%% file: main.tex
\documentclass[11pt]{amsart}
\pdfoutput=1
\usepackage{subfiles}
\input{preamble}

\title{Model categories for o-minimal geometry}
\author{Reid Barton \and Johan Commelin}

\begin{document}

\begin{abstract}
  We introduce a model category of spaces
  based on the definable sets of
  an o-minimal expansion of a real closed field.
  As a model category, it resembles the category of topological spaces,
  but its underlying category is a coherent topos.
  We will show in future work that its cofibrant objects
  are precisely the ``weak polytopes'' of Knebusch.
\end{abstract}

\maketitle

\subfile{tex/introduction}

\subfile{tex/omin}

\subfile{tex/morphism-omin}

\subfile{tex/categories}

\subfile{tex/proper-maps}

\subfile{tex/colimits}

\subfile{tex/realization}

\subfile{tex/homotopy}

\subfile{tex/topologies}

\subfile{tex/p_vs_d}

\subfile{tex/modelcat}

\subfile{tex/justleftproper}

\subfile{tex/related}

\appendix

\subfile{tex/cofibration}

\renewcommand*{\bibfont}{\small}
\printbibliography

\end{document}

%% file: preamble.tex
\usepackage[T1]{fontenc}
\usepackage[utf8]{inputenc}
\usepackage[english]{babel}

\usepackage{etoolbox}
\usepackage{ifthen}

\newbool{golf}
\booltrue{golf}
\newcommand{\golf}[2]{\ifbool{golf}{#1}{#2}}

\usepackage{microtype}
\usepackage{amsmath}
\usepackage{amssymb}
\usepackage{verbatim}
\usepackage{csquotes}
\usepackage{manfnt}

\usepackage{tikz-cd}

\usepackage[colorinlistoftodos]{todonotes}
\presetkeys{todonotes}{inline}{}
\makeatletter
\providecommand\@dotsep{5}
\def\listtodoname{List of Todos}
\def\listoftodos{\@starttoc{tdo}\listtodoname}
\makeatother

\usepackage{letltxmacro}
\LetLtxMacro{\oldtodo}{\todo}
\newcounter{todocounter}
\renewcommand{\todo}[2][]
{\stepcounter{todocounter}\oldtodo[#1]{\thetodocounter: #2}}

\newcommand{\lowpriotodo}[2][]{}

\usepackage{mathtools}

\usepackage{enumitem}

\setlist[enumerate,1]{label=\textup{(}\textit{\roman*}\textup{)},
                      ref  =\thetheorem.(\textit{\roman*})}

\usepackage{hyperref}
\usepackage[capitalise,nosort,nameinlink]{cleveref}

\crefname{enumi}{}{}
\crefname{enumii}{}{}
\crefname{enumiii}{}{}

\usepackage[backend=biber,doi=false,url=false,isbn=false,%
safeinputenc,sorting=nyt]{biblatex}
\bibliography{ominhtpy.bib}

\DeclareFieldFormat{postnote}{#1}
\DeclareFieldFormat{multipostnote}{#1}

\swapnumbers

\theoremstyle{definition}
\newtheorem{definition}{Definition}[section]
\newtheorem*{definition*}{Definition}
\newtheorem{construction}[definition]{Construction}
\newtheorem{convention}[definition]{Convention}
\newtheorem{example}[definition]{Example}
\newtheorem{nonexample}[definition]{Non-example}
\newtheorem{remark}[definition]{Remark}
\newtheorem{warning}[definition]{Warning}

\newtheorem{notecons}[definition]{Note for constructivists}
\newtheorem{speculation}[definition]{Speculation}

\theoremstyle{remark}           %
\newtheorem{nul}[definition]{}
\newenvironment*{nul*}{}{}      %

\theoremstyle{plain}
\newtheorem{proposition}[definition]{Proposition}
\newtheorem{lemma}[definition]{Lemma}
\newtheorem{theorem}[definition]{Theorem}
\newtheorem{corollary}[definition]{Corollary}

\newcommand\SetSymbol[1][]{\nonscript\:#1\vert\allowbreak\nonscript\:\mathopen{}}
\providecommand\given{} %
\DeclarePairedDelimiterX\set[1]\{\}{\renewcommand\given{\SetSymbol[\delimsize]}#1}

\newcommand{\NN}{\mathbb{N}}    %
\newcommand{\ZZ}{\mathbb{Z}}    %
\newcommand{\RR}{\mathbb{R}}    %
\newcommand{\CC}{\mathbb{C}}    %

\newcommand{\sS}{\mathcal{S}}   %

\newcommand{\DD}{\mathbb{D}}     %
\newcommand{\PP}{\mathbb{P}}     %

\newcommand{\tp}{\mathrm{top}}  %

\newcommand{\C}{C}              %

\newcommand{\op}{\mathrm{op}}   %

\newcommand{\Set}{\mathrm{Set}} %
\newcommand{\Top}{\mathrm{Top}} %
\newcommand{\CompHaus}{\mathrm{CompHaus}} %
\newcommand{\sSet}{\mathrm{sSet}} %
\DeclareMathOperator{\Sing}{Sing} %
\newcommand{\sSetfin}{{\sSet_{\mathrm{fin}}}}
\newcommand{\Topfin}{{\Top_{\mathrm{fcell}}}} %

\newcommand{\Hom}{\mathrm{Hom}} %
\newcommand{\id}{\mathrm{id}}   %

\newcommand{\Map}{\mathrm{Map}} %

\newcommand{\Shv}{\mathrm{Shv}} %

\newcommand{\PSh}{\mathrm{PSh}} %

\newcommand{\loc}{\mathrm{loc}} %

\newcommand{\Qalgre}{\tilde{\mathbb{Q}}}
\newcommand{\slin}{\mathrm{sl}} %
\newcommand{\sa}{\mathrm{sa}}   %
\newcommand{\an}{\mathrm{an}}   %
\newcommand{\anexp}{{\an,\exp}} %

\DeclareMathOperator*{\colim}{colim}

\newcommand{\yo}{\mathrm{y}}    %

\newcommand{\taucl}{{\tau_{\mathrm{cl}}}} %
\newcommand{\taupr}{{\tau_{\mathrm{pr}}}} %
\newcommand{\tauad}{{\tau_{\mathrm{ad}}}} %

\DeclareMathOperator{\Spec}{Spec}

\DeclareMathOperator{\ev}{ev}   %

\newcommand{\s}{\mathbf{s}}    %

\DeclareMathOperator{\sk}{sk}   %
\DeclareMathOperator{\Ex}{Ex}   %

\newcommand{\cof}{\mathrm{cof}} %
\DeclareMathOperator{\Ho}{Ho}   %
\newcommand{\lder}{\mathbb{L}}  %

\newcommand{\kerodon}[1]{%
 \cite[\href{https://kerodon.net/tag/#1}{Tag #1}]{kerodon}%
}

\definecolor{darkred}{rgb}{0.6,0,0}

\def\colon{:} %
\renewcommand{\emptyset}{\varnothing}
\def\subset{\subseteq}

%% file: tex/introduction.tex
\section{Introduction}

As models for the homotopy theory of spaces,
simplicial sets and topological spaces
have complementary advantages.
Simplicial sets arise naturally from algebraic constructions
such as the nerve of a category,
and they form a particularly well-behaved category.
Topological spaces are a familiar setting for
spaces of geometric origin
such as Lie groups and configuration spaces,
and their weak equivalences admit an elementary and direct description.
On the other hand,
topological spaces also exhibit certain phenomena
which are inconvenient in the context of the homotopy theory of spaces.
For example, a proof that $\pi_k(S^n) = 0$ for $k < n$
must take into account the existence of surjective continuous maps from $S^k$ to $S^n$.
Moreover, even a so-called convenient~category of topological spaces
(as in \cite{Steenrod67})
has rather poor exactness and compactness properties
compared to the category of simplicial sets.

In this work, we introduce a new model category of spaces
resembling that of topological spaces
but with better formal properties.
The basic idea is to replace topological spaces
by the definable sets of an \emph{o-minimal structure}\rlap{,}\footnote{%
    In this paper, ``o-minimal structure'' means by default
    an o-minimal expansion of a real closed field.
    We review the theory of o-minimal structures
    over the next several sections.
    }
a kind of geometric context for ``tame topology''.
Our main result is as follows.

\begin{theorem}[\cref{q-model-structure,quillen-eq}]
  \label{main-thm}
  Fix an o-minimal structure
  and let $\PP$ denote the category of closed and bounded definable sets
  and continuous definable functions,
  equipped with the Grothendieck topology generated by
  finite jointly surjective families.
  Then the geometric realization--$\Sing$ adjunction
  \[
    |{-}| : \sSet \rightleftarrows \Shv(\PP) : \Sing
  \]
  induces a model category structure on $\Shv(\PP)$
  in which every object is fibrant.
  Furthermore, this adjunction becomes a Quillen equivalence.
\end{theorem}

In this introduction, we first explain
the meaning of and the motivation for \cref{main-thm}.
After that, we discuss some central concepts of the paper
that do not appear explicitly in \cref{main-thm}.

\subsection*{O-minimal structures, and the homotopy theory of definable sets}

An o-minimal structure consists of, for each $n \ge 0$,
a chosen class of subsets of $R^n$ called the \emph{definable sets}.
Here, $R$ is a real closed field.
For simplicity, in this introduction we will take $R$ to be the real numbers $\RR$.

In the prototypical example of an o-minimal structure,
we declare the definable sets to be all the subsets of $\RR^n$
that can be expressed as boolean combinations
of sets of the form
$\{\,(x_1, \ldots, x_n) \in \RR^n \mid f(x_1, \ldots, x_n) \ge 0\,\}$
for polynomials $f \in \RR[x_1, \ldots, x_n]$.
These are the \emph{semialgebraic sets}
of real algebraic geometry.
They include, for example,
the standard simplices $|\Delta^n|$ and cubes $[0, 1]^n$
and the unit spheres $S^{n-1} \subset \RR^n$,
but not the set $\ZZ \subset \RR$
or the graph of the sine function.
A function between semialgebraic sets is called semialgebraic
if its graph is semialgebraic.
Semialgebraic sets are ``tame'';
for example, they have a well-behaved theory of dimension
which rules out the existence of
a surjective semialgebraic function from $S^k$ to $S^n$ if $k < n$.

Remarkably, it turns out that
the tameness of semialgebraic sets and functions
has little to do with the specific nature of polynomial inequalities.
Instead, it depends only on
basic structural properties of the class of semialgebraic sets
(such as closure under projection: the Tarski--Seidenberg theorem)
together with the fact that
the only semialgebraic subsets of $\RR^1$
are the finite unions of open intervals and single points.
This is the subject matter of o-minimality,
a branch of model theory which emerged in the 1980s
(\cite{vdD84}, and \cite{DefsetsI, DefsetsII, DefsetsIII}).
The theory of o-minimal structures
can be regarded as one realization of
the idea of ``topologie mod\'er\'ee'' of
Grothendieck's ``Esquisse d'un Programme''~\cite{Esquisse}
from around the same time,
which he proposed should form a new foundation for algebraic topology.
We include a summary of the relevant definitions and results of this theory
over the next few sections.
For now, the reader not already familiar with o-minimality
can continue to interpret ``definable set'' as meaning
a semialgebraic subset of some $\RR^n$.

We call a closed and bounded definable set a \emph{polytope},
and we write $\PP$ for the category of polytopes
and continuous definable functions between them;
a function is called definable if its graph is a definable set.
The polytopes play a role analogous to
that of the finite CW complexes or simplicial complexes in topology.
In fact, every polytope is definably homeomorphic to a finite simplicial complex
(the ``triangulation theorem'');
but the maps between polytopes depend on the choice of o-minimal structure.

There is an evident notion of homotopy between maps of $\PP$,
defined using the unit interval $[0, 1] \in \PP$ as cylinder object.
We define the homotopy category $\Ho \PP$ to be
the category with the same objects as $\PP$
and morphisms given by homotopy classes of maps.
By work of Delfs and Knebusch~\cite{DK85} and of Baro and Otero~\cite{BO10},
the category $\Ho \PP$ is known to be
equivalent to the homotopy category of finite CW complexes,
regardless of the choice of o-minimal structure or real closed field.

\subsection*{The category of sheaves}

In order to enlarge $\PP$ to a model for \emph{all} spaces,
we should perform some sort of cocompletion.
One possibility would be to take the free cocompletion of $\PP$,
that is, the presheaf category $\PSh(\PP)$.
An object of $\PSh(\PP)$ is something built out of
not just simplices (or cubes, or any other fixed family of shapes)
but rather all ``reasonable shapes'',
with the exact meaning of ``reasonable shape''
parameterized by the choice of o-minimal structure.
However, this construction ignores the ways in which
an object of~$\PP$ can already be built out of other objects of $\PP$.
For example,
gluing together two copies of the unit interval $[0, 1] \in \PP$ end-to-end
results in another interval $[0, 2] \in \PP$,
which is homeomorphic to the original interval $[0, 1]$
via the map $t \mapsto \frac12 t$.
This is what underlies the ability to concatenate paths in an object of $\PP$,
like in a topological space.
We would like to preserve this feature when cocompleting $\PP$,
in order to end up with a model category which resembles topological spaces.
This means that the pushout square in $\PP$
\[
  \begin{tikzcd}
    \{1\} \ar[r] \ar[d] & {[0, 1]} \ar[d] \\
    {[1, 2]} \ar[r] & {[0, 2]}
  \end{tikzcd}
  \tag{$*$}
\]
ought to be preserved in the cocompletion.

To construct a model category for spaces, which form an $\infty$-topos,
it is logical to use a 1-topos.
So, we are led to consider
the category of sheaves on $\PP$ for some topology.
Examples like ($*$) suggest taking the topology
in which the covering families of $X \in \PP$
are generated by finite \emph{closed} covers of~$X$.
Note that the topology generated by open covers
(more precisely, finite families of closed subsets whose interiors form a cover,
since the objects of $\PP$ behave like compact Hausdorff spaces)
is not fine enough for our purposes.

Topos-theoretic considerations indicate that
we should instead choose an even finer topology on $\PP$.
The category $\PP$ is a pretopos\rlap{,}%
\footnote{
  This observation seems not to have been made in the literature,
  though it follows fairly easily from the results of \cite{vdD98}.}
meaning it satisfies the Giraud pretopos axioms:
it has finite limits,
finite coproducts which are disjoint and pullback-stable,
and quotients of equivalence relations which are effective and pullback-stable.
On any pretopos $C$ there is the coherent topology,
generated by finite jointly effective epimorphic families.
The coherent topology is the finest subcanonical topology on $C$
generated by finite families,
and its category of sheaves is the initial topos
equipped with a pretopos morphism from $C$.
In the case of $\PP$,
the finite jointly effective epimorphic families
turn out to be just the finite families
that are jointly surjective on underlying sets.
For reasons to become apparent later,
we refer to this topology on $\PP$ as the \emph{proper} topology.

For the purpose of \cref{main-thm},
the main practical advantage of the proper topology
is its compatibility with the geometric realization
of a simplicial set.
We define the simplex functor $|\Delta^\bullet|_\PP : \Delta \to \PP$
by the same formulas as in topology.
Whatever (subcanonical) topology we choose on $\PP$,
its category of sheaves $\Shv(\PP)$
will be a cocomplete category which contains $\PP$ as a full subcategory,
so $\yo |\Delta^\bullet|_\PP : \Delta \to \PP \subset \Shv(\PP)$
extends to a colimit-preserving functor $|{-}|_{\Shv(\PP)} : \sSet \to \Shv(\PP)$.
On the other hand, while $\PP$ itself is not cocomplete,
it turns out to have enough colimits to construct
the geometric realization $|K|_\PP$ of any \emph{finite} simplicial set $K$.
Thus, there are functors
\[
  \begin{tikzcd}
    \sSetfin \ar[r, "|{-}|_\PP"] \ar[d] & \PP \ar[d, "\yo"] \\
    \sSet \ar[r, "|{-}|_{\Shv(\PP)}"] & \Shv(\PP)
  \end{tikzcd}
  \tag{$\dagger$}
\]
but in order for this square to commute (up to canonical isomorphism),
we have to give $\PP$ the proper topology
and not just the topology generated by closed covers.
The reason is that constructing the geometric realization $|K|_\PP$ of $K \in \sSetfin$
involves forming a series of pushouts
\[
  \begin{tikzcd}
    \coprod_{N_n K} |\partial \Delta^n|_\PP \ar[r] \ar[d] & {|{\sk^{n-1} K}|_\PP} \ar[d] \\
	  \coprod_{N_n K} |\Delta^n|_\PP \ar[r] & {|{\sk^n K}|_\PP}
  \end{tikzcd}
\]
in which the left map is a monomorphism,
but there is no reason that
the ``attaching maps'' $|\partial \Delta^n|_\PP \to |{\sk^{n-1} K}|_\PP$
of the nondegenerate $n$-simplices of~$K$
need be monomorphisms.
Pushouts of this form need not be preserved in $\Shv(\PP)$
when we give $\PP$ the topology generated by closed covers,
but they are preserved when we give $\PP$ the proper topology.

For the rest of the introduction we regard $\PP$ as equipped with the proper topology,
as in the statement of \cref{main-thm},
and write $|{-}|$ for the geometric realization valued in either $\PP$ or $\Shv(\PP)$.

\subsection*{The model category structure}

Once the basic properties of $\Shv(\PP)$ and geometric realization
have been established,
the proof of \cref{main-thm} is not difficult.
We transfer the standard model structure on simplicial sets
along the adjunction whose left adjoint is the bottom functor of ($\dagger$),
so that a map $f : X \to Y$ is a weak equivalence or a fibration of $\Shv(\PP)$
if and only if $\Sing f : \Sing X \to \Sing Y$ is one in $\sSet$.
The resulting weak factorization systems on $\Shv(\PP)$
are generated by the sets
\begin{align*}
  I & = \{\,|\partial \Delta^n| \to |\Delta^n| \mid n \ge 0\,\} \\
  & \cong \{\,S^{n-1} \to D^n \mid n \ge 0\,\}
\end{align*}
and
\begin{align*}
  J & = \{\,|\Lambda^n_i| \to |\Delta^n| \mid n > 0, 0 \le i \le n\,\} \\
  & \cong \{\,D^{n-1} \times \{0\} \to D^{n-1} \times [0, 1] \mid n > 0\,\}.
\end{align*}
Since each map of $J$ admits a retraction,
every object of $\Shv(\PP)$ is fibrant,
as in the model category of topological spaces.
From here, it is easy to check the conditions required
for the existence of the transferred model structure.
(Quillen originally used essentially the same method
to construct the model category structure on topological spaces in \cite{Qui67}.)

To verify that $|{-}| : \sSet \rightleftarrows \Shv(\PP) : \Sing$ is a Quillen equivalence,
it suffices to check that the unit map $\eta_K : K \to \Sing |K|$
is a weak equivalence for each $K \in \sSet$.
By expressing $K$ as a filtered colimit of finite simplicial sets,
we can reduce this to a statement about the top functor in the diagram~($\dagger$),
which we deduce from the result mentioned earlier that
the homotopy category of~$\PP$ is equivalent to that of finite CW complexes.

The resulting model category $\Shv(\PP)$
resembles the model category of topological spaces in many ways.
For example:
\begin{itemize}
\item
  As we saw above,
  its generating (acyclic) cofibrations
  can be described in the same way as for topological spaces.
  In particular, every object of $\Shv(\PP)$ is fibrant.
\item
  The model category $\Shv(\PP)$ is simplicial
  with the usual formulas for the simplicial mapping spaces, tensors and cotensors
  (a feature shared by cartesian closed versions of the category of topological spaces).
\item
  Objects of $\PP$, and in particular the spheres $S^n$,
  are cofibrant objects of $\Shv(\PP)$,
  and the unit interval $[0, 1] \in \PP$ serves as an interval object.
  Therefore, the homotopy groups $\pi_n$ can be defined
  in the same way as for topological spaces
  and are ``correct'' for all objects of $\Shv(\PP)$,
  in the sense that a map of $\Shv(\PP)$ is a weak equivalence
  if and only if it induces an isomorphism on all possible $\pi_n$
  (\cref{weq-iff-pi}).
\item
  The model category $\Shv(\PP)$ is right proper since every object is fibrant,
  and also left proper by \cref{left-proper}.
\end{itemize}
On the other hand:
\begin{itemize}
\item
  As a category, $\Shv(\PP)$ is a topos,
  in fact a coherent topos, whose coherent objects are exactly $\PP$.
  In particular, $\Shv(\PP)$ is locally \emph{finitely} presentable,
  since the topology on $\PP$ is generated by finite families.
  Likewise, the model category structure on $\Shv(\PP)$ is $\omega$-combinatorial
  in the sense that the generating (acyclic) cofibrations $I$ and $J$ described above
  have $\omega$-compact domains and codomains.
\item
  ``Gluing'' theorems describing the weak homotopy type of an object of $\Shv(\PP)$
  in terms of a cover by subobjects
  apply to arbitrary covers (in the topos-theoretic sense).
  In particular, they apply to covers of a polytope
  by finitely many \emph{closed} subsets.
  By contrast, in $\Top$ one typically has to use covers by open subsets.
  The hypothesis that the cover is open
  is typically used in conjunction with the Lebesgue number theorem
  and some kind of iterated subdivision.
  In the o-minimal setting,
  we rely on the triangulation theorem
  (and its variants such as normal triangulations \cite{BO10})
  in place of this style of argument,
  which is not available anyways if the field $R$ is nonarchimedean.

  These gluing theorems lie mostly outside the scope of this paper,
  but we do use a simple argument of this type (\cref{mono-mono-weq-pushout})
  in the proof that the model category $\Shv(\PP)$ is left proper.
\item
  Every cofibrant object of $\Shv(\PP)$ is an $I$-cell complex
  (not just a retract of an $I$-cell complex),
  and the full subcategory $\Shv(\PP)^\cof$ of cofibrant objects
  is closed under finite limits.
  (This fails for topological spaces
  since the equalizer of two continuous maps $f$, $g : [0, 1]^n \to [0, 1]$
  could be an arbitrary closed subset of $[0, 1]^n$.)
\end{itemize}

As a matter of fact, the category $\Shv(\PP)^\cof$
is already known in the literature on o-minimal structures:
it is equivalent to the category of \emph{weak polytopes}
defined in the semialgebraic case by Knebusch in~\cite{WSA}
using quite different methods
and generalized to the o-minimal setting in~\cite{Pi08}.
The advantage of the full model category $\Shv(\PP)$
is that it contains constructions
such as mapping spaces (exponentials) and loop spaces
which the category of weak polytopes lacks.
The definition of $\Shv(\PP)^\cof$ is also arguably
simpler than the original definition of weak polytopes,
and we hope that the methods of model category theory
can shed new light on this category.
The equivalence of $\Shv(\PP)^\cof$ to the category of weak polytopes
will be established in future work~\cite{paper3},
along with the claims in the last bullet point above.

The model category $\Shv(\PP)$ is parameterized by
the choice of a real closed field and an o-minimal structure.
In particular, the real closed field might be countable
(e.g., the real algebraic numbers $\Qalgre$),
or nonarchimedean.
It seems plausible that for suitable choices
(e.g., $\Qalgre$ with the semialgebraic structure,
which has an algorithm for quantifier elimination)
the model structure could be made constructive,
though we do not attempt to do so here.
See \labelcref{constructive} for further discussion.

\subsection*{General definable sets, and proper maps}

So far we have restricted our attention to the polytopes:
the definable sets which are closed and bounded.
This restriction is sensible in the context of algebraic topology,
where we generally build up objects from cells like simplices or disks
which are ``compact''---either topological spaces that are literally compact,
or the standard simplex $\Delta^n \in \sSet$
which we imagine as representing the compact topological space $|\Delta^n|$.
However, this restriction might not be appropriate in other settings.
For example, an affine algebraic variety $V$ over $R$
can be viewed as a semialgebraic set, which will rarely be bounded.
We could regard $V$ as defining a (non-representable) sheaf on $\PP$,
but this discards information about ``definability at infinity''
which we might prefer to retain.
For this reason, we will work for the most part
in the larger category $\DD$ of all definable sets and continuous definable functions.

Unlike $\PP$, the category $\DD$ is not a pretopos.
For example, let $X$ be a definable set with a nonclosed subset $W \subset X$,
and form the quotient of two copies of $X$
by the relation identifying the two copies of $W$.
This quotient exists in $\DD$
but it is the same as the quotient which identifies
the two copies of the closure $\overline W \supsetneq W$,
essentially because every object of $\DD$ is Hausdorff.
So, the original quotient cannot be effective.
(This situation does not arise in $\PP$
because $W$ would not be a polytope.)

There is however a restricted class of equivalence relations in $\DD$
which do have well-behaved quotients,
namely those which are componentwise \emph{proper}.
A map $p : X \to Y$ of $\DD$ is proper
if for each subset $K \subset Y$ which is a polytope,
the preimage $p^{-1}(K) \subset X$ is also a polytope.
The category $\DD$ together with the class of proper maps form a structure
which satisfies a restricted version of the Giraud pretopos axioms.
(Algebraic spaces and the class of \'etale maps
form another example of such a structure.)
It would take us too far afield to develop a general theory of such structures here,
but this notion is implicit in our treatment of
proper maps and quotients by proper equivalence relations in $\DD$.
From this perspective,
the reason that the theory simplifies for $\PP$
is that every map of $\PP$ is proper.

The topology we consider on $\DD$ is the one generated by
finitely jointly surjective families of proper maps
(hence the term ``proper topology'').
This is the analogue of the coherent topology on a pretopos
that takes into account the class of proper maps.

\subsection*{Characterization of sheaves for the proper topology}

The above definition of the proper topology
is convenient for topos-theoretic descent-type arguments,
but the corresponding sheaf condition
can be difficult to check
since we lack a simple description of a general surjective proper map.
Our main technical result gives
a second description of the sheaves for the proper topology,
akin to the characterization of the sheaves for the proper cdh-topology
of \cite{Voe10_unstable_motivic_homotopy}.

\begin{definition*}
  A \emph{distinguished square} in $\DD$ is a commutative square
  \[
    \begin{tikzcd}
      A \ar[r, "f"] \ar[d, "j"'] & B \ar[d, "j'"] \\
      X \ar[r, "f'"] & Y
    \end{tikzcd}
  \]
  satisfying the following equivalent conditions
  (\cref{closed-by-proper-iff-blowup}):
  \begin{enumerate}
  \item
    $j$ is a closed embedding,
    $f$ is proper,
    and the square is a pushout.
  \item
    $j'$ is a closed embedding,
    $f'$ is proper and induces an isomorphism from $X - j(A)$ to $Y - j'(B)$,
    and the square is a pullback.
  \end{enumerate}
\end{definition*}

\begin{proposition}[\labelcref{proper-sheaves-squares}]
  \label{proper-sheaf-iff}
  A functor $F : \DD^\op \to \Set$ is a sheaf for the proper topology
  if and only if $F(\emptyset) = {*}$
  and $F$ takes distinguished squares to pullbacks.
\end{proposition}

In algebraic geometry,
the definition of distinguished square (or abstract blow-up square)
is normally formulated using condition (\textit{ii}).
For our purposes, however, condition (\textit{i}) seems more useful
because it corresponds to the model category-theoretic idea
of building up a complex by attaching cells.
This is particularly true if we restrict attention to $\PP$,
in which case condition (\textit{i}) reduces to
``$j$ is a monomorphism and the square is a pushout''
(and it is somewhat unnatural to introduce the non-polytopes $X - j(A)$ and $Y - j'(B)$).
Henceforth, we refer to distinguished squares as ``closed-by-proper pushouts''.

The ``only if'' direction of \cref{proper-sheaf-iff}
can be rephrased by saying that
the Yoneda embedding $\yo : \DD \to \Shv(\DD)$
preserves the initial object and closed-by-proper pushouts;
this direction depends only on
the pretopos-like properties of $\DD$ along with its class of proper maps.
The ``if'' direction says that $\yo : \DD \to \Shv(\DD)$
is the universal cocompletion preserving these colimits,
and is an additional property of $\DD$.
We prove this direction using the same notion of a \emph{splitting sequence}
which appears in \cite{Voe10_unstable_motivic_homotopy}.

In this paper we do not really make use of the ``if'' direction
except to describe the model category $\Shv(\PP)$ in terms of a universal property
(\cref{model-category-universal}).
It will play a larger role in subsequent work.

\subsection*{Why o-minimality?}

From the perspective of a homotopy theorist,
the theory of o-minimal structures might seem like
a considerable amount of additional machinery.
Indeed, even proving that there are any o-minimal expansions of real closed fields at all
already requires the Tarski--Seidenberg theorem on projections of semialgebraic sets;
and the proof of the triangulation theorem (for instance) is quite intricate.
The reader might wonder whether
something simpler such as piecewise linearity could suffice.
There are two main reasons that we prefer to work in the o-minimal framework.
\begin{itemize}
\item
  Spaces of importance in homotopy theory that arise from geometry,
  such as the orthogonal and unitary groups,
  are frequently defined by nonlinear polynomial equations.
  These spaces generally have no natural representations as simplicial complexes
  but do exist naturally as objects in the o-minimal setting.
  (For example, $O(n)$ and $U(n)$ are group objects in $\PP$.)

  By working in a general o-minimal structure,
  we retain the option to enlarge the class of definable functions
  beyond the semialgebraic ones at our discretion,
  subject to the requirement that the structure remains o-minimal.
  For example, we can add the restrictions of
  the standard covering map $\RR \to S^1$, $t \mapsto \exp(2\pi i t)$
  to every compact interval.
  We list more examples (and non-examples) of o-minimal structures in the next section.
\item
  The piecewise linear category is too rigid
  to have a good supply of colimits.
  In particular, one cannot generally form a pushout
  which attaches a simplex along a non-injective attaching map\rlap{.}%
  \footnote{
    This limitation of the piecewise linear setting
    was also noted by Grothendieck
    in section~5 of~\cite{Esquisse}.
    }
  In an o-minimal expansion of a real closed field, however,
  the definability of multiplication provides enough ``flexibility''
  to build such pushouts.
  This is a theme which runs throughout the paper,
  eventually culminating in the fact that
  the geometric realization in $\Shv(\PP)$ of \emph{any} finite simplicial set
  is representable by an object of $\PP$.
  By expressing a general simplicial set as a filtered colimit of finite ones,
  we can thus reduce the question of the Quillen equivalence of $\sSet$ and $\Shv(\PP)$
  to a question about finite simplicial sets and $\PP$.
\end{itemize}
While we hope the reader will be interested in learning more about o-minimality,
we have tried to demonstrate here that
the theory of o-minimal structures
can be used effectively as a foundation
in a ``black-box'' fashion.

\subsection*{Organization of this paper}

In \cref{omin} we review the concept of an o-minimal structure,
recalling the main theorems that highlight its tameness properties.
We then define morphisms of o-minimal structures in \cref{morphism-omin}.
This concept is certainly not new,
but we are not aware of other discussions in the literature.
It combines the model-theoretic notions of expansions and reducts.
In \cref{cats} we introduce the categories $\DD$ and~$\PP$,
and prove some of their basic properties.

\Cref{proper-maps} reviews the class of proper morphisms,
which lie at the heart of this work.
They are the main ingredient for the proper topology,
which we use in the latter half of this article.
In \cref{colimits} we study colimits in~$\DD$ and~$\PP$;
properness is required for these to be well-behaved.
This section is a recasting of the last chapter of~\cite{vdD98}
in categorical language.

With \cref{geometric-realization} we start expounding
the homotopy theoretic aspects of~$\DD$ and~$\PP$.
We show that there is a geometric realization functor
from finite simplicial sets to~$\PP$.
In \cref{homotopy-d-p} we consider
the homotopy categories~$\Ho \DD$ and~$\Ho \PP$.
By a classical argument these categories are equivalent to each other,
and by results of Baro and Otero~\cite{BO10} and Delfs and Knebusch~\cite{DK85}
they are both equivalent to the homotopy category of finite CW complexes.
We then equip $\PP$ with the structure of a cofibration category%
\footnote{%
Roughly speaking \emph{cofibration categories}
are a variant of model categories
suitable for studying the homotopy theory of ``finite spaces'',
such as finite simplicial sets, compact Hausdorff spaces,
or objects of~$\PP$ (closed and bounded definable sets).
\Cref{cofibration-categories}
includes an introduction to cofibration categories.}
and deduce from the results mentioned above
that the geometric realization functor induces an equivalence of homotopy categories.

Having understood the homotopy category~$\Ho \PP$,
we turn our attention to categories of sheaves on~$\PP$,
a natural place for ``infinite spaces'', such as loop spaces.
In \cref{topologies} we study several Grothendieck topologies on~$\DD$ and~$\PP$.
The proper topology has the good property
that it makes the sheaf-valued geometric realization functor
agree with the representable one on all finite simplicial sets.
We study those categories of sheaves, and how they compare, in \cref{p_vs_d}.
Finally, we endow these categories of sheaves with
model category structures in \cref{modelcat}
and establish their basic model-categorical properties,
finishing with the proof that geometric realization is a Quillen equivalence
between $\sSet$ and $\Shv(\PP)$ or $\Shv(\DD)$.
In \cref{left-properness} we prove that
these model categories are left proper.

In \cref{related},
we touch upon some open ends
that should be treated in future work.
We also compare and contrast our results to other approaches in the literature.

\subsection*{Acknowledgments}

We thank
Mathieu Anel, Jonas Frey, Amador Martin-Pizarro and Inna Zakharevich
for helpful conversations related to this work.

The first author is supported by the Sloan Foundation (grant G-2018-10067).
The second author is supported by
DFG GK~1821 \emph{(Cohomological Methods in Geometry)}.

%% file: tex/omin.tex
\section{Review of o-minimal structures}
\label{omin}

\begin{nul*}
  In this section we give a brief introduction to o-minimality,
  with a particular focus on the tameness properties that it provides.
  The subject of o-minimal structures can be approached
  from geometry and from mathematical logic.
  We choose to present the basics of the theory from a geometric perspective.
  Some of the definitions are special cases
  of more general concepts in model theory,
  but here we will only appeal to model theory in a mild way.

  The book~\cite{vdD98} by van den Dries
  is an excellent introduction to the subject
  and includes most of the background theory
  that we need here
  (the main exception being the homotopical comparison results
  we will use in \cref{homotopy-d-p}).
  A rapid summary of the theory can also be found in
  the slightly earlier article~\cite{vdDM},
  particularly sections~2 and~4.
\end{nul*}

\begin{definition}[{\cite[1.2.1]{vdD98}}]
  \label{structure}
  A \emph{structure} on a set~$R$
  is a collection~$\sS = (\sS_n)_{n \in \NN}$
  that satisfies the following properties:
  \begin{enumerate}
    \item
      $\sS_n$ is a boolean subalgebra of the subsets of~$R^n$.
    \item
      If $X \in \sS_n$, then $R \times X$ and $X \times R$ are in~$\sS_{n+1}$.
    \item
      The set $\set{ (x_1, \ldots, x_n) \given x_1 = x_n }$ is in~$\sS_n$.
    \item
      If $X \in \sS_{n+1}$,
      and $\pi \colon R^{n+1} \to R^n$
      is the projection onto the first $n$ coordinates,
      then $\pi(X)$ is in~$\sS_n$.
  \end{enumerate}
  A set $X \subset R^n$ is called
  \emph{$\sS$-definable} or \emph{definable in~$\sS$}
  if $X \in \sS_n$.
  A function $f \colon X \to Y$ for $X \subset R^m$ and $Y \subset R^n$
  is \emph{$\sS$-definable} if its graph is $\sS$-definable.
  Usually the structure~$\sS$ will be clear from the context,
  and we simply speak of \emph{definable} sets and functions.
\end{definition}

\begin{nul}
  Let $\sS = (\sS_n)_{n \in \NN}$ be a structure on~$R$.
  For all $1 \le i$, $j \le n$,
  the set $\set{ (x_1, \dots, x_n) \given x_i = x_j }$ is definable.
  Any function $f : R^m \to R^n$ of the form
  $f(x_1, \ldots, x_m) = (x_{i(1)}, \ldots, x_{i(n)})$ is definable.
  The image or preimage of a definable set under a definable function
  is again definable.
  The composition of definable functions is definable.
  These basic facts can be found in~\cite[\S1.2]{vdD}.
\end{nul}

\begin{remark}
  \Cref{structure} may appear somewhat ad-hoc.
  As alluded to in the introductory paragraph of this section,
  it is a special case of a model-theoretic definition.
  Alternatively,
  it can also be rephrased using so-called hyperdoctrines~\cite{Lawvere_hyperdoctrines}.
  It is a sub-hyperdoctrine of the ``powerset hyperdoctrine on~$R$''
  whose category of ``types'' is the trivial Lawvere theory $\Set_{\mathrm{fin}}^\op$
  and whose category of ``attributes'' of type $X \in \Set_{\mathrm{fin}}^\op$
  is the poset of all subsets of~$R^X$.
\end{remark}

\begin{remark}
  \label{definable-iff-fol}
  Let $\sS$ be a structure on~$R$.
  Consider the signature $L_\sS$ with,
  for each $n \ge 0$,
  one $n$-ary relation symbol corresponding to
  each definable subset $X \subset R^n$,
  and one $n$-ary function symbol corresponding to
  each definable function $f : R^n \to R$.
  The set $R$ becomes a structure for $L_\sS$ (in the model-theoretic sense)
  in a tautological way.
  Then, for any first-order formula $\varphi$ in $n$ free variables
  over the signature $L_\sS$,
  the subset of $R^n$ that $\varphi$ cuts out is again definable.

  This observation is invaluable for verifying
  that a set is definable,
  as can be seen from \cref{closure-definable}.
  See~\S1.5 of~\cite{vdD98} for more details on this principle.
\end{remark}

\begin{nul}
  \label{eg-real-closed-field}
  We are mainly interested in structures on real closed fields.
  A field~$R$ is a \emph{real closed field}
  if one of the following equivalent conditions holds:
  \begin{enumerate}
    \item
      $R$ is not algebraically closed,
      but the field extension $R(\sqrt{-1})$ is algebraically closed.
    \item
      $R$ is \emph{elementarily equivalent} to the real numbers:
      every first-order formula in the language of fields is true in~$R$
      if and only if it is true in~$\RR$.
    \item
      $R$ is a totally ordered field,
      an element of $R$ is nonnegative if and only if
      it has a square root,
      and each polynomial $f \in R[X]$ of odd degree has a root in~$R$.
  \end{enumerate}
  Important examples of real closed fields include
  the real numbers~$\RR$
  and the algebraic real numbers~$\Qalgre$,
  which are both \emph{archimedean}:
  for every element~$x$,
  there exists a natural number~$n$ that is larger than~$x$.
  Nonarchimedean examples include
  the nonstandard real numbers from nonstandard analysis
  and the field of Puiseux series over a real closed field.

  The reader not familiar with real closed fields
  may assume that $R = \RR$
  (except in the proof of \cref{ho-p},
  which uses the real closed field $\Qalgre$).
\end{nul}

\begin{definition}[{\cite[1.3.2]{vdD98}}]
  \label{o-minimal-structure}
  Let $R$ be a set equipped with a dense linear order ${<}$ without endpoints,
  and let $\sS = (\sS_n)_{n \in \NN}$ be a structure on~$R$.
  Then $\sS$ is \emph{o-minimal} if it satisfies the following two properties:
  \begin{enumerate}
    \item
      The set $\set{ (x,y) \given x < y }$ is in~$\sS_2$.
    \item
      $\sS_1$ consists exactly of those sets that are finite unions
      of points and (possibly unbounded) open intervals in~$R^1 = R$.
  \end{enumerate}
  In this paper we assume, in addition to these basic axioms,
  that $R$ is an ordered ring
  and that the sets
  \[
    \set{(x,y,z) \given x + y = z}
    \quad\text{and}\quad
    \set{(x,y,z) \given x \cdot y = z}
  \]
  are definable.
  In this case $R$ must be a real closed field with its induced ordering
  \cite[Proposition~1.4.6]{vdD}.
  We say that $\sS$ \emph{extends} the real closed field~$R$.
\end{definition}

\begin{remark}
  \begin{enumerate}
    \item
      The `o' in `o-minimal' stands for `order'.
    \item
      Because of the finiteness assumption in the o-minimality condition,
      the set $\ZZ \subset \RR$ can never be definable
      in an o-minimal structure over~$\RR$.
      For the same reason,
      the functions
      $\cos : \RR \to \RR$ and $\exp \colon \CC \to \CC$ cannot be definable
      (where $\CC$ is identified with $\RR^2$).
      This might seem like a high price to pay
      for any benefits that o-minimality may bring.
      In \cref{real-examples} and \cref{exp-loc}
      we will see a partial work-around to this ``problem''.
    \item
      Let $f : X \to Y$ be a function between definable sets
      and $X_1$, \ldots, $X_k$ a cover of $X$
      by \emph{finitely many} definable sets.
      If the restriction of $f$ to each $X_i$ is definable,
      then $f$ is definable,
      since the graph of $f$ is the union of the graphs of the $f_i$.
      The finiteness hypothesis is necessary
      even if (say) $X = R = \RR$ and the $X_i$ are open intervals.
      For example, the floor function $\lfloor {-} \rfloor : \RR \to \RR$
      is definable in the semilinear structure of the next example
      when restricted to any bounded interval,
      but not on all of $\RR$.
  \end{enumerate}
\end{remark}

\begin{example}
  \label{omin-examples}
  Let $R$ be a real closed field.
  \begin{enumerate}
    \item
      A set $X \subset R^n$ is \emph{semilinear}
      if it belongs to the boolean algebra
      generated by affine half-spaces,
      that is,
      sets of the form
      $\set{ x \in R^n \given f(x) \ge 0 }$
      where
      $f : R^n \to R$ is an affine function.

      Semilinear sets form a structure on~$R$ that we denote by~$R_\slin$.
      The multiplication function ${-} \times {-} : R^2 \to R$ is not definable in~$R_\slin$
      and so by convention we do not refer to $R_\slin$ as an o-minimal structure.
      Definability of multiplication is actually
      not required for most of the results reviewed in this section,
      but it will come to play a crucial role later.
      From time to time we will use $R_\slin$ as a counterexample
      showing what goes wrong in its absence.
    \item
      A set $X \subset R^n$ is \emph{semialgebraic}
      if it belongs to the boolean algebra
      generated by sets of the form
      $\set{ x \in R^n \given f(x) \ge 0 }$
      with $f \in R[X_1, \dots, X_n]$.

      Semialgebraic sets form a structure on~$R$ that we denote by~$R_\sa$.
      The next theorem shows that this structure is o-minimal.
      Historically, the theory of o-minimality
      grew out of the observation that
      most of the pleasant properties of semialgebraic sets
      could in fact be proven
      using only the axioms of an o-minimal structure.
  \end{enumerate}
\end{example}

\begin{theorem}[Tarski--Seidenberg]
  Let $R$ be a real closed field.
  The projection of a semialgebraic set is again a semialgebraic set.
  Hence the structure~$R_\sa$ of semialgebraic sets
  is an o-minimal structure.
\end{theorem}

\begin{proof}
  Essentially by definition,
  the semialgebraic sets are those defined by
  quantifier-free formulas in the language of an ordered ring
  with parameters from~$R$.
  So, given such a formula $\varphi$ in $n+1$ variables,
  we need to find another quantifier-free formula $\psi$ in $n$ variables
  such that for every $y \in R^n$,
  \[
    \psi(y) \iff (\exists x \in R^{n+1}, \varphi(x) \land \pi(x) = y).
  \]
  This is an instance of
  the quantifier elimination theorem for real closed fields,
  which was first proven by Tarski~\cite{Tar51}.
  The remaining axioms for an o-minimal structure are easily verified.
\end{proof}

\begin{example}
  \label{real-examples}
  \begin{enumerate}
    \item
      Let $\RR_{\an}$ denote the smallest structure
      extending the real closed field~$\RR$
      which contains the graph of $f|_{[0,1]^m} : [0,1]^m \to \RR$
      for every real analytic function $f : \RR^m \to \RR$.
      Van den Dries showed in~\cite{Dries_1986_generalization_Tarski_Seidenberg}
      that $\RR_{\an}$ is an o-minimal structure.
    \item
      Let $\RR_{\exp}$ denote
      the smallest structure over~$\RR$
      extending the real closed field~$\RR$
      and containing the graph of $\exp \colon \RR \to \RR$.
      By a theorem of Wilkie~\cite{Wilkie_1996_Model_completeness}
      the structure $\RR_{\exp}$ is o-minimal.
    \item
      Van den Dries and Miller~\cite{Dries_Miller_1994_real_exponential_field}
      extended Wilkie's result (which was already announced in 1991),
      and showed that $\RR_{\anexp}$,
      the structure generated by $\RR_{\an}$ and~$\RR_{\exp}$, is o-minimal.
  \end{enumerate}
\end{example}

\begin{nul*}
  For the rest of this section,
  we fix a real closed field~$R$
  and an o-minimal structure~$\sS$ over~$R$.
  Whenever we say that a set or function is ``definable'',
  we mean definable in~$\sS$.
\end{nul*}

\begin{nul}
  \label{topology-basic}
  We now turn to topological matters.
  We equip $R$ with the order topology, each $R^n$ with the product topology,
  and each definable set $X \subset R^n$ with the subspace topology.
  We refer to the set $X$ equipped with this topology as
  the \emph{underlying topological space} of $X$,
  and denote it by $X_\tp$.
  The notions ``open'', ``closed'' and ``continuous''
  are interpreted relative to these topologies on definable sets.
  In practice, we are only interested in
  the \emph{definable} open and closed subsets of a definable set.
  However, the following fact means that
  the usual notion of continuity is still applicable.
\end{nul}

\begin{lemma}
  Let $f \colon X \to Y$ be a definable function
  between definable sets $X \subset R^m$ and $Y \subset R^n$.
  Then $f$ is continuous
  if and only if $f^{-1}(U)$ is open in $X$
  for every \emph{definable} open subset $U$ of $Y$.
  Similarly, we can also test continuity using definable closed subsets.
\end{lemma}

\begin{proof}
  The intersections of products of intervals with $Y$ are definable sets
  and form a base for the topology of $Y$.
\end{proof}

\begin{example}
  \label{closure-definable}
  If $X \subset R^n$ is definable
  then its closure $\overline X$ is also definable,
  using \cref{definable-iff-fol},
  since
  \begin{align*}
    & (y_1, \ldots, y_n) \in \overline X \iff \\
    & \qquad \forall \varepsilon > 0,
    \exists (x_1, \ldots, x_n) \in X,
    |x_1 - y_1| < \varepsilon \wedge \cdots \wedge |x_n - y_n| < \varepsilon.
  \end{align*}
\end{example}

\begin{nul}
  Likewise, the notions ``homeomorphism'', ``embedding'',
  and ``open (or closed, or dense) embedding'' from topology
  apply unchanged,
  except that we only apply them to definable functions.
  (Note that the inverse of a definable bijection
  is also definable.)
  Embeddings are characterized by a universal property
  relative to the map between underlying sets:
  they are the injective continuous definable functions $i : X \to Y$
  such that any continuous definable function $f : Z \to Y$
  factors (as a continuous definable function) through $i$
  if and only if the image of $f$ is contained in $i(X) \subset Y$.

  On the other hand, other notions like ``connected'' and ``compact''
  are not suitable in their usual forms.
  For instance, it is not hard to see that
  the underlying topological space of
  the interval $[0, 1] \subset R$
  is neither connected nor compact
  unless $R = \RR$.
  (In this text,
  we use the phrase ``$R = \RR$'' as shorthand for
  ``$R$ is Dedekind-complete''---or equivalently,
  (uniquely) isomorphic as an ordered field to the real numbers.)
\end{nul}

\begin{definition}[{\cite[Definition~1.3.5]{vdD}}]
  A definable set is \emph{definably connected}
  if it is nonempty
  and cannot be written as
  the disjoint union of two nonempty definable open subsets.
\end{definition}

\begin{definition}
  A subset $X \subset R^n$ is called \emph{bounded}
  if there exists an $r \in R$ such that $X \subset [-r, r]^n$.
  It is \emph{definably compact} if it is definable, closed and bounded.
  (There is also an important relative version of this notion, \emph{properness},
  which we discuss in \cref{proper-maps}.)
\end{definition}

\begin{remark}
  The notions of definable connectedness and definable compactness
  agree with the usual notions of connectedness and compactness
  in the case $R = \RR$.
  But even when $R \ne \RR$, a nonempty interval in $R$ is still definably connected.
  The key point is that
  any \emph{definable} nonempty, bounded-above subset of $R$
  has a least upper bound;
  this follows immediately from the definition of an o-minimal structure.
  Then the usual proof of connectedness of an interval goes through.
  Likewise, $[0, 1]$ is definably compact for every $R$.

  In fact, if $R$ is a countable real closed field,
  then the underlying topological spaces of
  $[0, 1]$ and $[0, 1] - \{1/2\}$ are homeomorphic by the back-and-forth method,
  since both are countable dense ordered sets with endpoints.
  But $[0,1]$ is definably connected and definably compact
  while $[0, 1] - \{1/2\}$ is neither.
  The homeomorphism between their underlying topological spaces
  produced by the back-and-forth method
  is necessarily non-definable.
\end{remark}

\begin{proposition}[{\cite[Proposition~6.1.10]{vdD98}}]
  \label{definably-compact-image}
  Let $f : X \to Y$ be a surjective continuous definable function.
  If $X$ is definably compact, then so is $Y$.
  In particular, the property of being definably compact
  is invariant under definable homeomorphism.
\end{proposition}

\begin{nul*}
  One of the tameness properties of
  definable sets in an o-minimal structure
  is a well-behaved notion of dimension.
  As an illustration of what \emph{cannot} happen,
  consider the graph $\Gamma \subset \RR^2$ of $x \mapsto \sin(\pi/x)$
  on the interval $(0,1)$.
  \[
    \begin{tikzpicture}
      \draw [domain=1:11,smooth,samples=201, line join=round] %
      plot [parametric, id=sin-inv-t-a] function {10/t, sin(3.1415926535*t)};
      \draw [domain=11:31.5,samples=411, line join=round] %
      plot [parametric, id=sin-inv-t-b] function {10/t, sin(3.1415926535*t)};
    \end{tikzpicture}
  \]
  The boundary $\overline \Gamma - \Gamma$ is the set
  $\set{ (0,y) \given y \in [-1,1] }$.
  For any reasonable homeomorphism-invariant definition of dimension,
  this would mean that $\Gamma$ and $\overline \Gamma - \Gamma$
  have the same dimension.
  We will see in \cref{dim-boundary-lt}
  that such pathologies do not occur in o-minimal geometry.
\end{nul*}

\begin{definition}
  \label{dimension}
  Let $X \subset R^n$ be a definable set.
  The dimension $\dim(X)$ of~$X$ is the maximal number~$k$
  such that there is an injective definable function
  $[0,1]^k \to X$.
  By convention, we say that $\dim(\varnothing) = -1$.
\end{definition}

\begin{remark}
  \cite{vdD98} uses a different definition of the dimension,
  which is better for establishing its basic properties.
  Once those have been established,
  it is easy to see that the dimension
  agrees with the definition we give above.
\end{remark}

\begin{lemma}[{\cite[Proposition~4.1.3]{vdD98}}]
  \label{dimension-facts}
  \begin{enumerate}
    \item
      If $X \subset Y \subset R^m$ are definable subsets,
      then $\dim(X) \le \dim(Y) \le m$.
    \item
      If $X \subset R^m$ and $Y \subset R^n$ are definable
      and there is a definable bijection between $X$ and~$Y$,
      then $\dim(X) = \dim(Y)$.
    \item
      \label{dim-union}
      If $X, Y \subset R^m$ are definable,
      then
      \[\dim(X \cup Y) = \max\{\dim(X), \dim(Y)\}.\]
  \end{enumerate}
\end{lemma}

\begin{lemma}[{\cite[Theorem~4.1.8]{vdD98}}]
  \label{dim-boundary-lt}
  Let $S \subset R^m$ be a nonempty definable set.
  Then $\dim(\overline S - S) < \dim(S)$.
  In particular, $\dim(\overline S) = \dim(S)$.
\end{lemma}

\begin{lemma}
  \label{dim-sdiff-interior-lt}
  Let $X \subset R^m$ be a nonempty definable set and $S \subset X$ be a definable subset,
  and let $U$ denote the interior of $S$ within $X$.
  Then $\dim(S - U) < \dim(X)$.
\end{lemma}

\begin{proof}
  If not, then we must have $\dim(S - U) = \dim(S) = \dim(X)$.
  This contradicts \cite[Corollary~4.1.9]{vdD98}.
\end{proof}

\begin{nul*}
  Another tameness property of o-minimal structures
  is the fact that an arbitrary definable function is ``almost'' continuous:
  it becomes continuous after a finite definable partition of its domain.
\end{nul*}

\begin{lemma}
  \label{continuous-on-partition}
  Let $f \colon X \to Y$ be a definable function between definable sets.
  Then there exists a finite partition of $X$
  into definable subsets~$X_i$
  such that $f|_{X_i} \colon X_i \to Y$ is continuous for each $i$.
\end{lemma}

\begin{proof}
  This follows from~\cite[6.2.5]{vdD98},
  setting $n = 0$ and $S = A = X$.
\end{proof}

\begin{remark}
  The reader might have expected the word ``continuous''
  to appear in \cref{dimension} and \cref{dimension-facts}.
  In fact, dimension is invariant even under arbitrary definable bijections,
  essentially because of \cref{continuous-on-partition}.
  For the same reason,
  the injective definable function in \cref{dimension}
  may be chosen to be continuous.
\end{remark}

\begin{nul*}
  Topological properties like continuity and closure
  can be detected using \emph{curves},
  which therefore play a role like that of sequences for metric spaces.
\end{nul*}

\begin{definition}
  \label{curve}
  Let $X \subset R^m$ be a definable set.
  A \emph{curve} is a continuous definable map $\gamma : (0, \varepsilon) \to X$
  for some $\varepsilon > 0$.
\end{definition}

\begin{remark}
  In \cite{vdD98} a curve is not assumed to be continuous,
  but on the other hand is assumed to be injective.
  This does not make much difference,
  because of the following fact.
\end{remark}

\begin{lemma}
  For any definable function $\gamma : (0, \varepsilon) \to R^n$
  there exists $0 < \varepsilon' < \varepsilon$ such that
  when restricted to $(0, \varepsilon')$,
  $\gamma$ is continuous and either constant or injective.
\end{lemma}

\begin{proof}
  This follows easily from the case $n = 1$,
  which is a consequence of the monotonicity theorem
  \cite[Theorem~3.1.2]{vdD98}.
\end{proof}

\begin{warning}
  We also make the left endpoint of a curve the endpoint of interest,
  and place it always at $0 \in R$,
  while \cite{vdD98} uses the right endpoint of an arbitrary interval
  $(a, b) \subset R$, \emph{including the case $b = +\infty$}.
  The choice of left or right endpoint is not important,
  but putting the endpoint at an actual element of $R$ (rather than $\infty$)
  makes a real difference in general.
  For example, in the semilinear structure
  there are no definable homeomorphisms
  between bounded intervals and unbounded intervals,
  and \cref{proper-iff-universally-closed} seems to be false.
  In our situation there is no problem
  with putting the endpoint at $0 \in R$
  because $(1, \infty)$, for example,
  is definably homeomorphic to $(0, 1)$ via $x \mapsto 1/x$;
  this is one place where we rely on the definability of multiplication.
\end{warning}

\begin{lemma}[Curve selection]
  \label{exists-curve-of-dense-embedding}
  Let $f \colon X \to Y$ be a definable dense embedding between
  definable sets $X \subset R^m$ and $Y \subset R^n$.
  For every point $y \in Y$,
  there exists a curve $\gamma \colon (0, \varepsilon) \to X$
  such that $f\gamma$ tends to~$y$ at~$0$.
\end{lemma}

\begin{proof}
  This follows from \cite[Corollary 6.1.5]{vdD98}.
\end{proof}

\begin{lemma}
  \label{continuous-iff-curve}
  Let $X \subset R^m$ and $Y \subset R^n$ be definable sets
  and let $f \colon X \to Y$ be a definable function.
  Let $x$ be a point of~$X$.
  Then $f$ is continuous at~$x$ if and only if
  for every curve $\gamma \colon (0, \varepsilon) \to X$
  with $\lim_{t \to 0} \gamma(t) = x$
  it is true that $\lim_{t \to 0} f(\gamma(t)) = f(x)$.
\end{lemma}

\begin{proof}
  The ``only if'' direction is obvious;
  for the ``if'' direction, apply \cite[Lemma~6.4.2]{vdD98}.
  We have to check the hypothesis for
  every (injective) definable function $\gamma : (0, \varepsilon) \to X$,
  not just the continuous ones;
  but we can make $\gamma$ continuous by reducing $\varepsilon$.
\end{proof}

\begin{nul*}
  Finally, we recall a few more results on o-minimal structures
  for later use.
\end{nul*}

\begin{theorem}[Definable choice]
  \label{definable-choice}
  Let $f \colon X \to Y$ be a definable surjection between definable sets.
  Then there exists a definable section $g \colon Y \to X$ of~$f$
  (that is, $f \circ g = \id$).
\end{theorem}

\begin{proof}
  Up to a slight reformulation, this is~\cite[Theorem~6.1.2]{vdD98}.
\end{proof}

\begin{lemma}[Shrinking lemma, {\cite[Lemma~6.3.6]{vdD98}}]
  \label{shrinking-lemma}
  Let $X \subset R^n$ be a definable set,
  covered by definable open subsets $U_1, \dots, U_k$.
  Then there exist, for $i = 1, \dots, k$,
  definable open subsets $V_i \subset X$
  with $\overline{V_i} \subset U_i$
  and such that $X = V_1 \cup \dots \cup V_k$.
\end{lemma}

\begin{nul}
  Another important result in o-minimal geometry is the triangulation theorem.
  This result says that there is a suitable notion of triangulation
  of a definable set in an o-minimal structure,
  and that every definable set can be triangulated in this sense.
  For details, we refer to chapter~8 of \cite{vdD}.
  We return to this topic in \cref{triangulation},
  at which point we have the language in place to give a precise statement.
\end{nul}

%% file: tex/morphism-omin.tex
\section{Morphisms of o-minimal structures}
\label{morphism-omin}

\begin{nul*}
  In this section we define morphisms of o-minimal structures.
  This terminology does not seem to appear elsewhere in the literature,
  although the basic idea is certainly not new.
  The relevant notions from model theory are
  \emph{expansions} and \emph{reducts} of structures.

  Our constructions in later sections
  will be functorial in the o-minimal structure,
  in a sense appropriate to the situation.
  The only place where we make essential use of this functoriality
  is in the proof that the homotopy theory of definable sets
  agrees with that of finite simplicial sets (\cref{ho-p}):
  via a zigzag of morphisms of o-minimal structures,
  we can compare the homotopy theory
  of a general o-minimal structure
  to that of $\RR_\sa$,
  where classical comparison theorems to topological spaces are available.
  Since verifying that these comparisons are equivalences of homotopy theories
  relies on difficult external results \cite{WSA,BO10},
  the reader might well treat \cref{ho-p} as a black box,
  in which case the material of this section could be considered optional.
\end{nul*}

\begin{definition}
  \label{omin-morphism}
  Let $R$ and $R'$ be two real closed fields,
  and let $\sS$ and $\sS'$
  be o-minimal structures expanding~$R$ and $R'$ respectively.
  A \emph{morphism of o-minimal structures} from $\sS$ to~$\sS'$
  is a pair $(\phi, \Phi)$
  consisting of a homomorphism of ordered rings $\phi \colon R \to R'$
  and a collection of functions
  $\Phi = (\Phi_n \colon \sS_n \to \sS'_n)_{n \in \NN}$
  satisfying the following conditions:
  \begin{enumerate}
  \item
    The functions $\Phi_n$ commute with
    the operations under which the definable sets of a structure are closed
    as described in \cref{structure}.
    Namely, the following equations hold whenever they make sense:
    \begin{align*}
      \Phi_n(\emptyset) & = \emptyset, \\
      \Phi_n(X \cup Y) & = \Phi_n(X) \cup \Phi_n(Y), \\
      \Phi_n(R^n - X) & = R'^n - \Phi_n(X), \\
      \Phi_{n+1}(R \times X) & = R' \times \Phi_n(X), \\
      \Phi_{n+1}(X \times R) & = \Phi_n(X) \times R', \\
      \Phi_n(\{\,(x_1, \ldots, x_n) \in R^n \mid x_1 = x_n\,\}) & =
      \{\,(x_1, \ldots, x_n) \in R'^n \mid x_1 = x_n\,\}, \\
      \Phi_n(\pi_R(X)) & = \pi_{R'}(\Phi_{n+1}(X))
    \end{align*}
    where in the last equation,
    $\pi_R : R^{n+1} \to R^n$ and $\pi_{R'} : R'^{n+1} \to R'^n$
    are the projections onto the first $n$ coordinates.
  \item
    For every $r \in R$, we have $\Phi_1(\{r\}) = \{\phi(r)\}$.
  \item
    $\Phi_2(\{\,(x, y) \in R^2 \mid x < y\,\})
    = \{\,(x, y) \in R'^2 \mid x < y\,\}$.
  \item
    Likewise,
    $\Phi_3$ takes the graphs of addition and multiplication on $R$
    to the graphs of the corresponding operations on $R'$.
  \end{enumerate}
  Because of condition (\textit{ii}),
  the homomorphism~$\phi$ is completely determined by~$\Phi_1$.
  For this reason,
  we will abuse notation and also write $\Phi$ for this induced function $R \to R'$,
  as well as for the entire morphism of o-minimal structures.
  We also suppress the subscripts on the individual operations $\Phi_n$.
\end{definition}

\begin{example}
  \label{morphism-omin-examples}
  \begin{enumerate}
    \item
      If $R$ is a real closed field,
      and $(R, \sS)$ and $(R, \sS')$ are two o-minimal structures with $\sS \subset \sS'$,
      then we obtain a morphism of o-minimal structures $(\id, \Phi)$,
      where $\Phi$ denotes the inclusion: $\Phi(X) = X$.
    \item
      Let $R \subset R'$ be a field extension of real closed fields.
      Then we get a morphism of o-minimal structures $\Phi : R_\sa \to R'_\sa$,
      defined as follows.
      If $X$ is a semialgebraic set over~$R$,
      and $\psi$ is a first-order formula
      in the language of ordered rings with constants from $R$ that describes~$X$,
      then we can interpret $\psi$ in~$R'$ to obtain a semialgebraic set over~$R'$;
      let $\Phi(X)$ be this set.
      (In the language of (semi)algebraic geometry:
      we can find polynomial inequalities that cut $X$ out of $R^n$,
      and proceed to view those polynomial inequalities over~$R'$.)

      The resulting semialgebraic set does not depend on the choice of~$\psi$,
      because the theory of real closed fields is \emph{model complete}:
      every embedding of real closed fields is an elementary embedding.
      Hence, if $\psi'$ is a second formula which also defines $X$,
      the statement $\forall x, \psi(x) \leftrightarrow \psi'(x)$ holds in $R$
      and therefore also in $R'$.
      Then it is also clear that $\Phi$ satisfies
      all the conditions of \cref{omin-morphism}
      because they correspond to operations on first-order formulas.
  \end{enumerate}
\end{example}

\begin{lemma}
  \label{Qsa-initial}
  The o-minimal structure $\Qalgre_\sa$
  of semialgebraic sets over the algebraic real numbers
  is initial in the category of o-minimal structures
  with morphisms as defined in \cref{omin-morphism}.
\end{lemma}

\begin{proof}
  Any o-minimal structure $(R, \sS)$
  admits a morphism from $\Qalgre_\sa$
  constructed by composing the two examples of \cref{morphism-omin-examples}.
  This morphism is unique because $\Qalgre$ is the initial real closed field
  and the semialgebraic structure is the smallest possible one,
  so the value of a morphism $\Phi : \Qalgre_\sa \to (R, \sS)$
  on any semialgebraic set $X$ is uniquely determined
  by the requirements on~$\Phi$.
\end{proof}

\begin{example}
  Let $\Phi : (R, \sS) \to (R', \sS')$ be any morphism of o-minimal structures.
  Then
  \[
    \Phi(\{\,(x, y) \in R^2 \mid 0 \le x \le y \le 1\,\})
    = \{\,(x, y) \in R'^2 \mid 0 \le x \le y \le 1\,\}
  \]
  because we can build the set $\{\,(x, y) \in R^2 \mid 0 \le x \le y \le 1\,\}$
  from the graphs of the relations $=$ and $\le$ and the sets $\{0\}$ and $\{1\}$
  by taking products with $R$, boolean operations, and projections,
  and $\Phi$ preserves all of these.
  ($\Phi(0) = 0$ and $\Phi(1) = 1$ because $\Phi$ is a ring homomorphism.)

  Note that if the field extension $R \subset R'$ is nontrivial,
  then the resulting set $\{\,(x, y) \in R'^2 \mid 0 \le x \le y \le 1\,\}$
  contains points which are not of the form $(\Phi(x), \Phi(y))$.
  Only when $X$ is finite is $\Phi(X)$ obtained by applying $\Phi$ ``pointwise''
  (and componentwise) to the elements of $X$.
\end{example}

\begin{nul}
  \label{first-order-prop}
  Fix $r$, $n_1$, \ldots, $n_r \ge 0$
  and let $P$ be a predicate of $r$ definable sets
  $X_1 \subset R^{n_1}$, \ldots, $X_r \subset R^{n_r}$
  in any o-minimal structure $(R, \sS)$
  with the following property:
  there exists a single sentence $\psi$
  in the language of an ordered ring
  augmented by relation symbols of arity $n_1$, \ldots, $n_r$,
  such that $P(X_1, \ldots, X_r)$ holds
  if and only if $\psi$ is valid in $R$
  when the $i$th relation symbol is interpreted as $X_i \subset R^{n_i}$.
  For example,
  \begin{itemize}
  \item
    with $r = 1$, $P(X)$ might be the predicate ``$X$ is nonempty''
    or the predicate ``$X \subset R^{n_1}$ is closed'';
  \item
    with $r = 3$ and $n_3 = n_1 + n_2$,
    $P(X, Y, Z)$ might be the predicate
    ``$Z$ is the graph of a continuous function from $X$ to $Y$''.
    (The function in question is automatically definable
    because we only let $X$, $Y$, $Z$ range over definable sets.)
  \end{itemize}
  We refer to such predicates as \emph{first-order properties}.

  Let $\Phi : (R, \sS) \to (R', \sS')$ be a morphism of o-minimal structures,
  and let $X_1 \subset R^{n_1}$, \ldots, $X_r \subset R^{n_r}$
  be $\sS$-definable sets.
  Then we have the following ``transfer principle'':
  \begin{itemize}
  \item[($*$)]
    For any first-order property $P$,
    $P(\Phi(X_1), \ldots, \Phi(X_r))$ holds
    if and only if $P(X_1, \ldots, X_r)$ does.
  \end{itemize}
  The proof is by induction on the structure of
  the formula $\psi$ defining $P$,
  using the fact that $\Phi$ preserves all the operations
  by which the interpretations of formulas
  are constructed from those of subformulas.
  As we will only need a finite number of instances of this transfer principle,
  which could be verified by hand,
  we leave the details to the reader.
\end{nul}

\begin{lemma}
  \label{omin-morphism-basic-props}
  Let $\Phi \colon (R,\sS) \to (R',\sS')$ be a morphism of o-minimal structures.
  \begin{enumerate}
  \item
    For every $\sS$-definable subset~$X$,
    we have $\Phi(X) = \varnothing \iff X = \varnothing$.
  \item
    For all $\sS$-definable subsets $X$ and $Y$,
    we have $\Phi(X) \subset \Phi(Y) \iff X \subset Y$.
  \item
    If $f \colon X \to Y$ is a definable function with graph $\Gamma_f$,
    then $\Phi(\Gamma_f)$ is the graph of a definable function
    from $\Phi(X)$ to $\Phi(Y)$.
  \end{enumerate}
  We denote this function by $\Phi f : \Phi X \to \Phi Y$.
  In the remaining items,
  all functions are assumed to be definable functions between definable sets.
  \begin{enumerate}[resume]
  \item If $f \colon X \to Y$ is continuous
    then $\Phi f \colon \Phi X \to \Phi Y$ is continuous.
  \item If $X$ is closed and bounded
    then $\Phi X$ is closed and bounded.
  \item If $f : X \to Y$ is injective
    then $\Phi f$ is injective.
  \item If $f : X \to Y$ is a homeomorphism
    then $\Phi f$ is a homeomorphism.
  \item If $f : X \to Y$ is an embedding
    then $\Phi f$ is an embedding.
  \item If $f_1 : X_1 \to Y$, \ldots, $f_n : X_n \to Y$ are jointly surjective,
    then $\Phi f_1$, \ldots, $\Phi f_n$ are jointly surjective.
  \item $\Phi$ preserves identities and composition.
  \end{enumerate}
\end{lemma}

\begin{proof}
  These all follow from the transfer principle.
  (Continuity can be described by a first-order formula
  using the $\varepsilon$-$\delta$ definition.
  A continuous bijection is a homeomorphism
  when the transpose of its graph is also the graph of a continuous function.
  An injective continuous function $f : X \to Y$ is an embedding
  when it is a homeomorphism onto its image.)
\end{proof}

%% file: tex/categories.tex
\section{\texorpdfstring
{The categories $\DD$ and $\PP$}
{The categories D and P}}
\label{cats}

\begin{definition}
  \label{categories-D-and-P}
  \begin{enumerate}
    \item
      We write $\DD = \DD(R,\sS)$ for the category of definable sets.
      More precisely, an object of $\DD$ is a pair $(n, X)$
      where $n \ge 0$ and $X$ is a definable set in~$R^n$.
      A morphism of $\DD$ is a continuous definable function.
    \item
      We write $\PP = \PP(R,\sS)$ for the full subcategory of~$\DD$
      on the definable sets which are definably compact
      (i.e., closed and bounded).
      We refer to the objects of $\PP$ as \emph{polytopes}.
  \end{enumerate}
\end{definition}

\begin{nul*}
  From this section onward, every mention of definability
  is with respect to a fixed o-minimal structure $(R, \sS)$,
  unless explicitly stated otherwise.
\end{nul*}

\begin{nul}
  By \cref{definably-compact-image},
  the property of being definably compact
  is invariant under definable homeomorphism.
  Therefore $\PP$ is a replete subcategory of~$\DD$:
  an object of $\DD$ that is isomorphic to an object in~$\PP$
  also belongs to~$\PP$.
\end{nul}

\begin{nul}
  We write $U : \DD \to \Set$ for the forgetful functor,
  and ${-}_\tp : \DD \to \Top$ for the functor
  sending a definable set $X$ to its underlying topological space $X_\tp$.
  We sometimes write $u : \PP \to \DD$ for the inclusion,
  but usually suppress it from the notation.
\end{nul}

\begin{nul}
  \label{limits-coproducts}
  The categories $\DD$ and~$\PP$ have all finite limits,
  which are computed ``as in $\Set$''.
  For example, the equalizer $E$ of two maps $f$, $g : X \to Y$
  is constructed as $E = \{\,x \in X \mid f(x) = g(x)\,\}$.
  If $X$ and $Y$ belong to $\PP$ (or even if just $X$ does)
  then $E$ again belongs to $\PP$.
  Finite products are also constructed in the obvious way.

  Similarly, there is no problem with coproducts.
  Let $X_1$, \dots, $X_m$ be a finite collection of definable sets.
  By enlarging their ambient spaces if necessary,
  we may assume each $X_i$ is given as a subset of the same $R^n$.
  Then we set
  \[
    X = \bigcup_{1 \le i \le m} \{i\} \times X_i \subset R^{1+n}.
  \]
  One easily verifies that the obvious maps $e_i : X_i \to X$
  make $X$ into the coproduct of the $X_i$ in $\DD$.
  If every $X_i$ is a polytope, then so is $X$,
  which is then also the coproduct of the $X_i$ in~$\PP$.

  For general finite colimits, however, the situation is much more subtle.
  We will discuss these in detail in \cref{colimits}.
\end{nul}

\begin{example}
  Let $G$ be a classical Lie group, for example $G = \mathrm{GL}(n, R)$.
  Then $G$ is naturally a group object in the category $\DD$.
  If $G$ is (definably) compact, for example $G = \mathrm{O}(n, R)$,
  then $G$ is a group object in $\PP$.
  These statements are valid for any real closed field $R$.
\end{example}

\begin{nul}
  \label{associated-functor}
  Let $\Phi \colon (R,\sS) \to (R',\sS')$ be a morphism of o-minimal structures.
  Then there is a functor
  \begin{align*}
    \DD(R,\sS) &\to \DD(R',\sS') \\
    X &\mapsto \Phi(X) \\
    f &\mapsto \Phi(f);
  \end{align*}
  see \cref{omin-morphism-basic-props} for more on~$\Phi(f)$.
  This functor restricts to a functor $\PP(R,\sS) \to \PP(R',\sS')$,
  since $\Phi(X)$ is a polytope if $X$ is.
  We call these functors the \emph{associated functors} of~$\Phi$,
  and denote them also by $\Phi$.
\end{nul}

%% file: tex/proper-maps.tex
\section{Proper maps}
\label{proper-maps}

\begin{nul*}
  Resuming our review of o-minimal structures,
  this section is devoted to an important class of morphisms of $\DD$,
  the \emph{proper} maps.
  Properness is a key hypothesis in the construction of well-behaved colimits,
  which we will discuss in the next section.
\end{nul*}

\begin{definition}[{\cite[Definition~6.4.4]{vdD98}}]
  \label{proper-map}
  Let $X$ and $Y$ be definable sets.
  A continuous definable function $f \colon X \to Y$ is \emph{proper}
  if the preimage of every definably compact subset of $Y$ is definably compact.
\end{definition}

\begin{proposition}
  \label{proper-iff-rr}
  If $R = \RR$, then a map $f : X \to Y$ of $\DD$ is proper
  if and only if it is proper as a map between the underlying topological spaces.
\end{proposition}

\begin{proof}
  See \cite[6.4.8 Exercise~3]{vdD98}.
\end{proof}

\begin{remark}
  For $R \ne \RR$, properness is not the same as
  properness of the map between underlying topological spaces,
  so for consistency we ought to write ``definably proper'' (as in \cite{vdD98}).
  However, this notion is so important that
  we prefer to use the shorter name for the definable version.
\end{remark}

\begin{lemma}
  \label{proper-of-iso}
  \label{proper-comp}
  \begin{enumerate}
    \item
      Isomorphisms are proper.
    \item
      The composition of proper morphisms is proper.
    \item
      A definable set $X$ is definably compact
      if and only if the unique map $X \to *$ is proper.
    \item
      Any morphism $f : X \to Y$ with $X$ definably compact is proper.
  \end{enumerate}
\end{lemma}

\begin{nul}
  In particular, every morphism of $\PP$ is proper.
  This simple observation is the reason that
  $\PP$ has better categorical properties than $\DD$.
  For example, $\PP$ turns out to be a pretopos while $\DD$ is not one.
  However, many parts of the theory
  (particularly \cref{colimits,topologies} here, and \cite{paper3})
  are perhaps surprisingly resilient to
  carrying this extra properness assumption where required.
\end{nul}

\begin{proposition}
  \label{proper-iff-curve}
  Let $f \colon X \to Y$ be a morphism of~$\DD$.
  Then $f$ is proper if and only if
  the following condition is satisfied:
  if $\gamma \colon (0, 1) \to X$ is a curve
  such that $f \circ \gamma \colon (0, 1) \to Y$ extends to $[0, 1)$,
  then $\gamma$ also extends to $[0, 1)$.
\end{proposition}

\begin{proof}
  This is \cite[Lemma 6.4.5~(1)]{vdD98},
  up to minor differences in the meaning of ``curve''.
  (Obviously, there is no difficulty if $\gamma$ is constant.)
\end{proof}

\begin{remark}
  \label{valuative-criterion}
  The condition in \cref{proper-iff-curve}
  can also be rephrased as saying that any commutative square
  \[
    \begin{tikzcd}
      {(0, 1)} \ar[r, "\gamma"] \ar[d] & X \ar[d, "f"] \\
      {[0, 1)} \ar[r] \ar[ru, dashed] & Y
    \end{tikzcd}
    \]
    admits a diagonal filler as shown by the dashed arrow.
    (The bottom triangle automatically commutes
    if the top one does,
    since definable sets are Hausdorff
    and so a curve defined on $(0, 1)$ has at most one extension to $[0, 1)$.
    For the same reason, the diagonal filler is unique if it exists.)

    In other words,
    this condition provides
    a ``definable valuative criterion of properness''.
\end{remark}

\begin{lemma}
  \label{proper-pullback}
  The pullback of a proper morphism is proper.
\end{lemma}

\begin{proof}
  This is an immediate consequence of \cref{valuative-criterion}.
\end{proof}

\begin{lemma}
  \label{proper-iff-pullback}
  Let $p \colon X \to Y$ be a morphism of~$\DD$.
  Then $p$ is proper if and only if
  for any polytope~$K$ and morphism $f \colon K \to Y$ of~$\DD$,
  the pullback $X \times_Y K$ is also a polytope.
\end{lemma}

\begin{proof}
  The ``if'' direction follows by applying the hypothesis
  to definably compact subsets $K$ of $Y$.
  On the other hand, suppose $p$ is proper
  and $f : K \to Y$ is any continuous definable map from a polytope $K$.
  Then the pullback $p' : X \times_Y K \to K$ of $p$ along $f$
  is again proper,
  so $X \times_Y K = p'^{-1}(K)$ is a polytope since $K$ is one.
\end{proof}

\begin{definition}
  A continuous definable map $p : X \to Y$ is \emph{definably closed}
  if the image under $p$ of any definable closed subset of $X$
  is closed in $Y$;
  and \emph{universally definably closed}
  if every pullback of $p$ is definably closed.
\end{definition}

\begin{warning}
  Unlike continuity,
  closedness of the map between underlying topological spaces
  cannot be checked on definable sets when $R \ne \RR$.
  See \cref{definably-closed-not-closed} for a counterexample.
\end{warning}

\begin{lemma}
  \label{closed-of-proper}
  Let $p \colon X \to Y$ be a morphism of~$\DD$.
  If $p$ is proper, then $p$ is definably closed.
\end{lemma}

\begin{proof}
  See \cite[Exercise~6.4.8]{vdD}.
  Hint: Use curve selection (\cref{exists-curve-of-dense-embedding}),
  definable choice (\cref{definable-choice}),
  and \cref{proper-iff-curve}.
\end{proof}

\begin{proposition}
  \label{proper-iff-universally-closed}
  Let $p \colon X \to Y$ be a morphism of~$\DD$.
  Then $p$ is proper if and only if $p$ is universally definably closed.
\end{proposition}

\begin{proof}
  Assume $p$ is proper.
  Since properness is stable under pullback,
  it suffices to show $p$ is definably closed,
  which we know by \cref{closed-of-proper}.

  Conversely, assume that $p$ is universally definably closed,
  and let $\gamma \colon (0,1) \to X$ be a curve.
  Assume that $p \circ \gamma$ extends to $\gamma' \colon [0,1) \to Y$.
  \[
    \begin{tikzcd}
      (0,1) \ar[d, hook] \ar[r, dashed, "\iota"'] \ar[rr, bend left=20, "\gamma"]
      & {[0,1) \times_Y X} \ar[d] \ar[r] & X \ar[d, "p"] \\
      {[0,1)} \ar[r, "\id"] & {[0,1)} \ar[r, "\gamma'"] & Y
    \end{tikzcd}
  \]
  We obtain an induced curve
  $\iota \colon (0,1) \to [0,1) \times_Y X$
  as in the diagram.
  Let $Z$ be the closure of the image of this curve.
  By assumption, the image of $Z$ under the projection
  $[0,1) \times_Y X \to [0,1)$ is closed, and it contains $(0,1)$.
  It is therefore equal to $[0,1)$,
  and we conclude that $\iota$ extends to $[0,1)$.
  Therefore $\gamma \colon (0,1) \to X$
  extends to $[0,1)$ as well,
  and we conclude that $p$ is proper by \cref{proper-iff-curve}.
\end{proof}

\begin{proposition}
  \label{closed-iff-injective-proper}
  A morphism of $\DD$ is a closed embedding if and only if
  it is injective and proper.
  In particular, a proper bijection is a homeomorphism.
\end{proposition}

\begin{proof}
  The ``only if'' direction is obvious.
  Suppose $f \colon X \to Y$ is injective and proper.
  Proper maps are definably closed,
  so the image $X' = f(X)$ is a (definable) closed subset of $Y$
  and the inverse function $g \colon X' \to X$ is continuous.
  Therefore $f$ is a definable homeomorphism onto its (closed) image $X'$.
\end{proof}

\begin{nul}
  In particular, the monomorphisms in $\PP$ are the closed embeddings.
\end{nul}

\begin{proposition}
  \label{proper-left-cancellation}
  If $f : X \to Y$ and $g : Y \to Z$ are maps of $\DD$ with $gf$ proper,
  then $f$ is proper.
\end{proposition}

\begin{proof}
  It suffices to find liftings in squares of the form
  \[
    \begin{tikzcd}
      (0, 1) \ar[r, "h"] \ar[d] & X \ar[d, "f"] \\
      {[0, 1)} \ar[r, "k"] \ar[ru, dashed] & Y
    \end{tikzcd}
    \tag{$*$}
  \]
  but any lift $[0, 1) \to X$ in the square
  \[
    \begin{tikzcd}
      (0, 1) \ar[r, "h"] \ar[d] & X \ar[d, "gf"] \\
      {[0, 1)} \ar[r, "gk"] \ar[ru, dashed] & Z
    \end{tikzcd}
  \]
  will do, because the bottom triangle of ($*$) will then automatically commute
  as noted in \cref{valuative-criterion}.
\end{proof}

\begin{proposition}
  \label{proper-right-cancellation}
  Let $(f_i : X_i \to Y)_{i = 1}^n$ be a finite family of maps of $\DD$
  which are jointly surjective,
  and let $g : Y \to Z$ be a map such that $gf_i$ is proper for each $i$.
  Then $g$ is proper.
\end{proposition}

\begin{proof}
  Since the hypotheses are preserved by pulling back along a map $Z' \to Z$,
  it suffices to check that $g$ is a definably closed map.
  If $V \subset Y$ is closed
  then $g(V) = \bigcup_{i=1}^n (g \circ f_i)(f_i^{-1}(V))$ is also closed
  because each $g \circ f_i$ is proper.
  (A remark: by \cref{proper-left-cancellation},
  the hypotheses also imply that each $f_i$ is proper.)
\end{proof}

\begin{nul*}
  We end this section with a proof that properness
  is a first order property of a morphism of~$\DD$,
  and is therefore preserved by morphisms of o-minimal structures.
  Our argument is essentially
  an explicit unwinding of the proof of \cite[Theorem~B.4]{LSA}.
  It uses the concept of a completion of a morphism,
  which we introduce next.
\end{nul*}

\begin{definition}[10.2.6 of~\cite{vdD98}]
  \label{completion}
  Let $f \colon X \to Y$ be a continuous definable map between
  definable sets $X \subset R^m$ and $Y \subset R^n$.
  A \emph{completion} of $f$ is a commutative square
  \[
    \begin{tikzcd}
      X \rar \ar[d, "f"'] & X' \dar{f'} \\
      Y \rar & Y'
    \end{tikzcd}
  \]
  where $X' \subset R^{m'}$ and $Y' \subset R^{n'}$
  are closed and bounded definable sets,
  $f'$ is a continuous definable function,
  and the horizontal maps are dense embeddings.
\end{definition}

\begin{construction}
  \label{standard-completion}
  Let $f \colon X \to Y$ be a continuous definable map between
  definable sets $X \subset R^m$ and $Y \subset R^n$.
  We will now construct a completion of~$f$,
  which we call the \emph{standard completion} of~$f$.

  Let $c \colon R \to (-2, 2)$ be the definable homeomorphism
  \[
    c(x) =
    \begin{cases}
      -2 - x^{-1}            &\text{if $x < -1$,} \\
      \phantom{-2 - {}}x     &\text{if $-1 \le x \le 1$,} \\
      \phantom{-}2 - x^{-1}  &\text{if $1 < x$.}
    \end{cases}
  \]
  By applying $c$ coordinatewise, we may assume that $X$ and~$Y$ are bounded.
  Let $\Gamma_f \subset R^m \times R^n$ be the graph of~$f$,
  and let $\overline{\Gamma_f}$ be the closure of $\Gamma_f$.
  It is closed and bounded,
  and contained in $\overline X \times \overline Y$,
  where $\overline X$ and $\overline Y$ are the closures of~$X$ and~$Y$ respectively.
  Then, consider the diagram
  \[
    \begin{tikzcd}
      X \ar[d, "f"'] \rar{\sim} & \Gamma_f \rar[hook] & \overline{\Gamma_f} \dar{\pi} \\
      Y \ar[hook]{rr} & & \overline Y
    \end{tikzcd}
  \]
  where $\pi \colon \overline{\Gamma_f} \to \overline Y$ denotes the natural projection.
  The horizontal maps are dense embeddings.
  Hence the diagram constitutes a completion of~$f$.
\end{construction}

\begin{lemma}
  \label{pullback-iff-set-pullback}
  Let
  \[
    \begin{tikzcd}
      X \ar[r] \ar[d] & X' \ar[d] \\
      Y \ar[r] & Y'
    \end{tikzcd}
  \]
  be a commutative square in $\DD$
  in which the horizontal morphisms are embeddings.
  Then the square is a pullback square in $\DD$
  if and only if the forgetful functor takes it to a pullback square of sets.
\end{lemma}

\begin{proof}
  The forgetful functor preserves pullbacks,
  and the converse direction follows from the universal property of embeddings.
\end{proof}

\begin{lemma}
  \label{proper-iff-completion}
  Let
  \[
    \begin{tikzcd}
      X \rar \ar[d, "f"'] & X' \dar{f'} \\
      Y \rar & Y'
    \end{tikzcd}
  \]
  be a completion square.
  Then $f$ is proper if and only if the square is a pullback.
\end{lemma}

\begin{proof}
  If the square is a pullback, then $f$ is proper, by \cref{proper-pullback}.

  Now assume that $f$ is proper.
  Identifying $X$ and $Y$ with subsets of $X'$ and $Y'$ respectively,
  it suffices to show $(f')^{-1}(Y) \subset X$.
  Suppose $x \in X'$ is a point whose image under $f'$ belongs to $Y$.
  By \cref{exists-curve-of-dense-embedding}
  there exists a curve $\gamma \colon (0,1) \to X$ that tends to~$x$ at~$0$.
  In particular $f\gamma$ tends to~$f'(y)$ at~$0$,
  and since $f$ is proper, this means that $x \in X$,
  by \cref{proper-iff-curve}.
  This concludes the proof.
\end{proof}

\begin{proposition}
  Properness is a first order property
  of a map of $\DD$
  (in the sense of \cref{first-order-prop})
  and so is preserved by morphisms of o-minimal structures.
\end{proposition}

\begin{proof}
  We apply \cref{proper-iff-completion}
  to the standard completion that we construct in \cref{standard-completion}.
  Let $f : X \to Y$ be a continuous definable map
  and assume that $X \subset R^m$ and $Y \subset R^n$ have already been made bounded.
  Then $f$ is proper if and only if
  for each point $(x, y) \in R^{m+n}$ belonging to the closure of the graph of $f$
  and with $y \in Y$,
  we have $x \in X$ and $y = f(x)$.
\end{proof}

%% file: tex/colimits.tex
\section{\texorpdfstring
{Colimits in $\DD$ and $\PP$}
{Colimits in D and P}}
\label{colimits}

\begin{nul*}
  We have already seen that $\DD$ has finite limits and finite coproducts,
  computed in the obvious way (\cref{limits-coproducts}),
  and that $\PP$ is closed under these constructions.
  For general finite colimits, the situation is considerably more subtle:
  not all pushouts exist in the category $\DD$
  (see \cref{closed-by-improper}),
  and some that do exist are ``bad''
  because they do not have the correct underlying set
  (see \cref{bad-colimits}).
  A class of colimits that exist and have the expected properties
  are the \emph{quotients by proper equivalence relations},
  constructed in \cite[\S10.2]{vdD98}.
  Here we present the results of that section
  in category-theoretic language.

  We begin with a characterization of coproducts in $\DD$
  which we will use to understand their interaction with pullbacks.
  (For coproducts this approach is not really necessary,
  but we will use the same method for quotients later.)
\end{nul*}

\begin{lemma}
  \label{coproduct-recognition}
  A diagram $(f_i : X_i \to Y)_{1 \le i \le n}$ is a coproduct in $\DD$
  if and only if the forgetful functor $U : \DD \to \Set$
  sends it to a coproduct diagram in $\Set$
  and each $f_i$ is a closed embedding.
\end{lemma}

\begin{proof}
  The ``only if'' direction is clear from the construction of finite coproducts.
  For the ``if'' direction,
  construct the induced map $f : \coprod_i X_i \to Y$.
  Using the ``only if'' direction,
  we see that $Uf$ is a bijection,
  and proper by \cref{proper-right-cancellation}
  because the inclusions $X_i \to \coprod_i X_i$ are jointly surjective.
  Then by \cref{closed-iff-injective-proper}, $f$ is an isomorphism
  and so the original diagram $(f_i : X_i \to Y)_{1 \le i \le n}$ is a coproduct.
\end{proof}

\begin{proposition}
  Finite coproducts in $\DD$ are disjoint and pullback-stable,
  and preserved by the forgetful functor $U : \DD \to \Set$.
\end{proposition}

\begin{proof}
  We already saw that $U$ preserves finite coproducts,
  and their disjointness follows.
  The pullback-stability of coproducts follows from the lemma
  together with the corresponding fact for $\Set$
  and the pullback-stability of closed embeddings.
\end{proof}

\begin{nul}
  Hence, the category $\DD$ is finitary extensive \cite{extensive_distrib_cat}.
  In particular, a square of the form
  \[
    \begin{tikzcd}
      X_i \ar[r] \ar[d, "f_i"'] & \coprod_i X_i \ar[d, "\coprod_i f_i"] \\
      Y_i \ar[r] & \coprod_i Y_i
    \end{tikzcd}
  \]
  is a pullback.
\end{nul}

\begin{proposition}
  Let $(f_i : X_i \to Y_i)_{i=1}^n$ be a family of maps of $\DD$.
  Then their coproduct $\coprod f_i : \coprod X_i \to \coprod Y_i$ is proper
  if and only if each $f_i$ is proper.
\end{proposition}

\begin{proof}
  Each $f_i$ is a pullback of $\coprod f_i : \coprod X_i \to \coprod Y_i$,
  which proves the ``only if'' direction.
  For ``if'', apply \cref{proper-right-cancellation}
  to the coproduct inclusions $X_i \to \coprod_i X_i$,
  noting that the composition $X_i \to \coprod_i X_i \to \coprod_i Y_i$
  equals $X_i \xrightarrow{f_i} Y_i \to \coprod_i Y_i$, which is proper.
\end{proof}

\begin{nul*}
  On the other hand, general coequalizers in $\DD$
  need not be so well-behaved even when they exist.
  In particular, they may not be preserved by the forgetful functor to $\Set$.
\end{nul*}

\begin{example}
  \label{bad-colimits}
  Take $X = Y = R = \RR$ and consider the maps $f : X \to Y$ and $g : X \to Y$
  given by $f(x) = x$ and $g(x) = x + 1$.
  In $\Set$, the coequalizer of $f$ and~$g$
  would be the circle $\RR/\ZZ$.
  But in $\DD$, the coequalizer of $f$ and~$g$ is the unique map $Y \to *$.
  This amounts to the claim that for any object~$Z$ of $\DD$,
  the only maps $h : Y \to Z$ such that $h f = h g$ are the constant maps.
  Indeed, suppose $h$ is such a map.
  Let $y \in Y$ and consider the set $\{\,y' \in Y \mid h(y') = h(y)\,\}$.
  This set contains $y$ and it is invariant under translation by $1$,
  so by o-minimality it can only be all of $Y = \RR$.

  The previous example involved objects of $\DD$.
  However, we can easily construct similar examples in $\PP$.
  For example, take $X = Y$ to be a circle, $f$ the identity
  and $g$ a rotation by an irrational multiple of $2\pi$.
  Or take $X = Y = [0, 1]$, $f$ the identity
  and $g : X \to Y$ defined by $g(x) = x/2$.
  In both cases it is not hard to see that the coequalizer of $f$ and $g$
  is again a point.

  At the end of this section (\cref{closed-by-improper}),
  we will give examples of finite colimits which do not exist in $\DD$.
  We do not know whether $\PP$ has all finite colimits.
\end{example}

\begin{nul*}
  The ``good'' coequalizers in $\DD$ are
  the quotients of proper equivalence relations.
  We review the relevant results from \cite[chapter~10]{vdD98}.
\end{nul*}

\begin{definition}[10.2.13 of~\cite{vdD98}]
  Let $X$ be a definable set.
  A \emph{definable equivalence relation} on $X$ is
  a definable subset $E \subset X \times X$
  which is (the graph of) an equivalence relation.
  A \emph{proper equivalence relation}
  is a definable equivalence relation $E$
  such that the two projections $E \subset X \times X \to X$
  are both proper morphisms.
\end{definition}

\begin{nul}
  Suppose for the moment that
  $E$ is any \emph{internal} equivalence relation on an object $X \in \DD$
  \cite[Definition~A.1.3.6]{Elephant1}.
  This means that $E$ is equipped with a monomorphism $i : E \to X \times X$
  and there are maps $r : X \to E$, $s : E \to E$, $t : E \times_X E \to E$
  which correspond to the axioms of an equivalence relation.
  In general, $i$ need not be an embedding;
  when it is, such an $E$ corresponds precisely to
  a definable equivalence relation as defined above.

  We write $p_1 : E \to X$ and $p_2 : E \to X$ for the two projections.
  They are interchanged by the isomorphism $s : E \to E$,
  hence isomorphic as morphisms of $\DD$.
  Therefore, if one of the projections is proper, so is the other.
  In this case, the map $i : E \to X \times X$ is also proper
  by \cref{proper-left-cancellation},
  since its composition with the projection to the first factor is $p_1$.
  Because $i$ is also a monomorphism, it is in fact a closed embedding.
  Hence, we conclude:
  \begin{itemize}
  \item
    A proper equivalence relation $E \subset X \times X$
    is the same as an internal equivalence relation $i : E \to X \times X$
    for which the compositions $E \to X \times X \to X$ are proper.
  \item
    If $E$ is a proper equivalence relation on $X$,
    then $E \subset X \times X$ is a closed subset.
  \end{itemize}
\end{nul}

\begin{definition}[10.2.2 and~10.2.3 of~\cite{vdD98}]
  \label{proper-quotient}
  Let $E$ be a definable equivalence relation on the definable set $X$.
  A \emph{definable quotient} of $X$ by $E$ is
  a definable set $Y$ together with a surjective continuous definable map $q : X \to Y$
  such that
  \begin{enumerate}
  \item[(a)] $q(x_1) = q(x_2)$ if and only if $(x_1, x_2) \in E$;
  \item[(b)] if $K$ is a definable subset of $Y$
    such that $q^{-1}(K)$ is closed in $X$,
    then $K$ is closed in $Y$.
  \end{enumerate}
  The quotient is called \emph{proper}
  if in addition the map $q$ is proper.
\end{definition}

\begin{nul}
  If $q : X \to Y$ is any proper surjective map,
  then the subset $E$ of $X \times X$ defined by
  $E = \{\,(x_1, x_2) \mid q(x_1) = q(x_2)\,\}$
  is a proper equivalence relation on $X$,
  since the projections $p_i : E \to X$
  are pullbacks of $q : X \to Y$.
\end{nul}

\begin{nul*}
  The converse also holds.
  This is the main result of \cite[\S10.2]{vdD98}.
\end{nul*}

\begin{theorem}[{\cite[Theorem 10.2.15]{vdD98}}]
  \label{exists-proper-quotient}
  Let $E$ be a proper equivalence relation on $X$.
  Then there exists a proper quotient of $X$ by $E$.
\end{theorem}

\begin{nul*}
  We now reformulate these results
  in categorical language.
\end{nul*}

\begin{proposition}
  \label{coequalizer-of-proper-surjective}
  Let $q : X \to Y$ be a proper surjective map of $\DD$.
  Then $q$ is the coequalizer of its kernel pair
  $E = X \times_Y X \rightrightarrows X$.
\end{proposition}

\begin{proof}
  Let $g : X \to Z$ be a morphism
  whose compositions with the two projections $E \rightrightarrows X$ are equal.
  Then $g$ factors uniquely through $q$ as a map of sets $h : Y \to Z$.
  By \cref{definable-choice}, the map $q : X \to Y$
  has a definable (but generally not continuous) section $s : Y \to X$.
  Then $h = gs$, so $h$ is definable.
  To check that $h$ is continuous, let $V \subset Z$ be a closed subset.
  Then $h^{-1}(V) = q(g^{-1}(V))$ is closed
  because a proper map is definably closed (\cref{closed-of-proper}).
\end{proof}

\begin{nul*}
  We now turn to the converse direction,
  that is, constructing a (proper) quotient of a proper equivalence relation.
  As we saw earlier, these can be identified with
  certain internal equivalence relations in $\DD$.
  In what follows, we denote such a proper equivalence relation
  by the notation $E \rightrightarrows X$.
\end{nul*}

\begin{proposition}
  \label{quotient-of-proper-eq}
  Let $E \rightrightarrows X$ be a proper equivalence relation on $X$.
  Then there is a coequalizer diagram
  $E \rightrightarrows X \xrightarrow{q} Y$,
  and $q$ is proper and surjective and has kernel pair $E \rightrightarrows X$.
  This coequalizer is preserved by the forgetful functor $U : \DD \to \Set$.
\end{proposition}

\begin{proof}
  Let $q : X \to Y$ be the proper quotient provided by \cref{exists-proper-quotient}.
  By definition,
  $q$ is proper and surjective and has kernel pair $E$,
  so we just need to check that $q$ is in fact the coequalizer of
  $E \rightrightarrows X$.
  Let $Z$ be any object of $\DD$
  and $g : X \to Z$ a map with $g \circ p_1 = g \circ p_2$;
  we must show $g$ descends to a unique map $h : Y \to Z$.
  The uniqueness follows from surjectivity of $q$.
  From condition (a) in \cref{proper-quotient},
  we see that
  $U$ takes $E \rightrightarrows X \xrightarrow{q} Y$ to a coequalizer
  and in particular there exists a function on the level of underlying sets
  $h : Y \to Z$ with $hq = g$.
  The function $h$ is definable by definable choice,
  as in the proof of \cref{coequalizer-of-proper-surjective}.
  It remains only to check that $h$ is continuous.
  Let $V \subset Z$ be a closed subset;
  then $q^{-1}(h^{-1}(V)) = g^{-1}(V) \subset X$ is closed,
  so by condition (b) in \cref{proper-quotient},
  $h^{-1}(V) \subset Y$ is closed as well.
\end{proof}

\begin{lemma}
  \label{proper-quotient-recognition}
  Let $E \rightrightarrows X$ be a proper equivalence relation on $X$.
  Then a diagram $E \rightrightarrows X \xrightarrow{q} Y$ is a coequalizer
  if and only if $U : \DD \to \Set$ takes it to a coequalizer
  and $q$ is proper.
  (Compare \cref{coproduct-recognition}.)
\end{lemma}

\begin{proof}
  The ``only if'' follows from the above result,
  since we know now that there exists a coequalizer with the claimed properties.
  Conversely, for the ``if'' direction,
  construct the coequalizer $E \rightrightarrows X \xrightarrow{q'} X/E$;
  then there is an induced map $r : X/E \to Y$
  which is bijective
  (since $U$ takes both
  $E \rightrightarrows X \xrightarrow{q} Y$ and
  $E \rightrightarrows X \xrightarrow{q'} X/E$
  to coequalizers)
  and proper
  (using \cref{proper-right-cancellation},
  because its composition with the proper surjective map $q'$
  is the proper map $q$).
  By \cref{closed-iff-injective-proper}, $r$ is an isomorphism.
\end{proof}

\begin{proposition}
  \label{proper-quotients-d}
  Proper equivalence relations in $\DD$ have effective quotients
  which are again proper and pullback-stable,
  and preserved by the forgetful functor $U : \DD \to \Set$.
\end{proposition}

\begin{proof}
  We just need to prove the pullback-stability;
  this follows from the preceding characterization of quotients,
  together with the pullback-stability of proper maps
  and of quotients of equivalence relations in $\Set$.
\end{proof}

\begin{nul}
  \label{p-pretopos}
  The full subcategory $\PP \subset \DD$
  is closed under these quotients,
  since if $q : X \to Y$ is surjective and $X \in \PP$
  then $Y \in \PP$ as well (\cref{definably-compact-image}).
  Furthermore, every internal equivalence relation in $\PP$
  is automatically proper.
  Since $\PP$ is also closed under finite limits and finite coproducts,
  we conclude that \emph{$\PP$ is a pretopos}
  \cite[Expos\'e~VI, Exercice~3.11]{SGA4}.

  Any pretopos has a \emph{coherent topology}
  generated by the finite families $(f_i : X_i \to Y)$
  for which the induced map $\coprod_i X_i \to Y$ is
  an effective epimorphism, that is, the quotient of its kernel pair.
  By \cref{coequalizer-of-proper-surjective,quotient-of-proper-eq}
  the effective epimorphisms in $\PP$ are precisely the (proper) surjective maps,
  and so we can describe the coherent topology on $\PP$ more directly
  as the one generated by finite jointly surjective families.
  We will meet this topology again in \cref{topologies}.

  The category $\DD$ is not a pretopos, because it is not balanced.
  For example, the map $\{0\} \amalg (0, 1] \to [0, 1]$
  is a monomorphism and an epimorphism but not an isomorphism.
  However, it does satisfy a version of the Giraud axioms for a pretopos
  provided appropriate properness conditions are added;
  this is the content of \cref{proper-quotients-d}.
  This perspective will be exploited more systematically in \cite{paper3}.
\end{nul}

\begin{nul}
  An important special class of quotients by proper equivalence relations
  are the \emph{closed-by-proper pushouts}.
  (In fact, \cite{vdD98} constructs this class of pushouts first
  along the way to constructing all quotients of proper equivalence relations.)

  We refer to these pushouts as ``closed-by-proper''
  even when working in the category $\PP$, in which all morphisms are proper.
  This is to avoid ambiguity with the ``closed-by-closed'' pushouts
  we consider later.
\end{nul}

\begin{proposition}
  \label{closed-by-proper-pushout}
  Let $j : A \to B$ and $f : A \to X$ be maps of $\DD$
  with $j$ a closed embedding and $f$ proper.
  Then there is a pushout square
  \[
    \begin{tikzcd}
      A \ar[r, "f"] \ar[d, "j"'] & X \ar[d, "j'"] \\
      B \ar[r, "f'"] & Y
    \end{tikzcd}
  \]
  in which $j'$ is a closed embedding and $f'$ is proper.
  This pushout is pullback-stable
  and preserved by the forgetful functor $U : \DD \to \Set$.
  The square is also a pullback,
  and $j'$ and $f'$ are jointly surjective.
\end{proposition}

\begin{proof}
  On $X \amalg B$ we consider the equivalence relation $E$
  generated by the relation $f(a) = j(a)$ for each $a \in A$.
  The equivalence relation generated by an \emph{arbitrary} relation
  need not be definable,
  because it may involve chains of generating relations of unbounded length.
  But in this case, the injectivity of $j : A \to B$
  implies that we need only consider chains of length at most two.
  Specifically, $E \subset (X \amalg B) \times (X \amalg B)$ is the image of
  a certain map
  \[
    (A \times_X A) \amalg A \amalg A \amalg X \amalg B \to (X \amalg B) \times (X \amalg B)
  \]
  whose summands encode the relations
  \[
    \mbox{$j(a) \sim j(a')$ when $f(a) = f(a')$}, \quad
    j(a) \sim f(a), \quad
    f(a) \sim j(a),
  \]
  \[
    x \sim x, \quad
    b \sim b.
  \]
  From this description it is also easy to check that
  $E$ is a \emph{proper} equivalence relation on $X \amalg B$,
  because $j$ and $f$ are proper.

  Now, let $q : X \amalg B \to Y$ be the quotient of this equivalence relation,
  and define $j'$ and $f'$ to be the compositions
  $X \to X \amalg B \xrightarrow{q} Y$ and
  $B \to X \amalg B \xrightarrow{q} Y$ respectively.
  Then $Y$ has the correct universal property to be
  the pushout of $j$ and $f$ in $\DD$.
  Furthermore, $j'$ and $f'$ are proper because $q$ is.
  The remaining properties
  (including the injectivity of $j'$, making it a closed embedding)
  follow from the corresponding properties
  for quotients by proper equivalence relations.
\end{proof}

\begin{nul}
  \label{closed-by-closed-pushout}
  An even more special case of quotients by proper equivalence relations
  are the \emph{closed-by-closed pushouts}.
  First, let $X \in \DD$ and let $V_1$, $V_2 \subset X$ be closed subsets
  whose union is $X$.
  Then the square
  \[
    \begin{tikzcd}
      V_1 \cap V_2 \ar[r] \ar[d] & V_1 \ar[d] \\
      V_2 \ar[r] & X
    \end{tikzcd}
    \tag{$*$}
  \]
  is a pushout.
  This is easily seen directly
  because continuous definable functions
  can be glued along finite closed covers.
  Furthermore, a map $f : X' \to X$
  induces a closed cover $V'_i = f^{-1}(V_i)$ of $X'$
  and a similar square for $V'_1$ and $V'_2$ which is also a pushout,
  so pushout squares of this form are pullback-stable
  (which we could also deduce from \cref{closed-by-proper-pushout}).

  Conversely, suppose
  $V_{12} \to V_1$ and $V_{12} \to V_2$
  are any two closed embeddings in $\DD$.
  By \cref{closed-by-proper-pushout},
  they have a pushout $X = V_1 \amalg_{V_{12}} V_2$
  in which the induced maps $V_1 \to X$ and $V_2 \to X$
  are also closed embeddings.
  Furthermore, this square is also a pullback square,
  so identifying $V_1$ and $V_2$ with closed subsets of $X$,
  we have $V_{12} = V_1 \cap V_2$.
  In summary, every closed-by-closed pushout
  is isomorphic to a square of the form ($*$).
\end{nul}

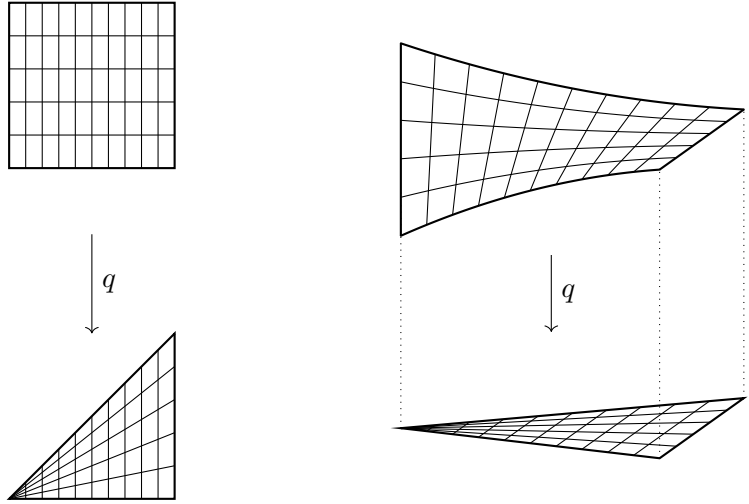
\begin{figure}
  \centering
  \begin{minipage}[c]{2in}
    \begin{tikzpicture}[scale=2.2]
      \foreach \v in {-0.6, -0.2, ..., 0.61} {
        \path [draw]
        (0, {0.5 * \v + 2.5}) -- (1, {0.5 * \v + 2.5});
      }
      \foreach \u in {0.1, 0.2, ..., 0.91} {
        \path [draw] (\u, 2) -- (\u, 3);
      }
      \path [draw, thick] (0, 2) -- (1, 2) -- (1, 3) -- (0, 3) -- cycle;

      \draw [->] (0.5, 1.6) -- (0.5, 1) node [right, midway] {$q$};

      \foreach \v in {-0.6, -0.2, ..., 0.61} {
        \path [draw]
        (0, 0) -- (1, {0.5 * \v + 0.5});
      }
      \foreach \u in {0.1, 0.2, ..., 0.91} {
        \path [draw] (\u, 0) -- (\u, \u);
      }
      \path [draw, thick] (0, 0) -- (1, 0) -- (1, 1) -- cycle;
    \end{tikzpicture}
  \end{minipage}
  \begin{minipage}[c]{2in}
    \begin{tikzpicture}[scale=4]
      \foreach \a/\b/\c/\h in {0.14/0.1/0.32/3} {
        \foreach \v in {-1, 1} {
          \draw [variable=\u,domain=0:1,smooth,thick,line cap=round]
          plot
          ({\u + \a * \u * \v},
           {\b * \u * \v + \c * (\h + \v * (1 - \u) * (1 - 0.5 * \u))});
        };
        \foreach \v in {-0.6, -0.2, ..., 0.61} {
          \draw [variable=\u,domain=0:1,smooth]
          plot
          ({\u + \a * \u * \v},
           {\b * \u * \v + \c * (\h + \v * (1 - \u) * (1 - 0.5 * \u))});
        };

        \foreach \u in {0, ..., 1} {
          \draw [variable=\v,domain=-1:1,smooth,thick]
          plot
          ({\u + \a * \u * \v},
           {\b * \u * \v + \c * (\h + \v * (1 - \u) * (1 - 0.5 * \u))});
        };
        \foreach \u in {0.1, 0.2, ..., 0.91} {
          \draw [variable=\v,domain=-1:1,smooth]
          plot
          ({\u + \a * \u * \v},
           {\b * \u * \v + \c * (\h + \v * (1 - \u) * (1 - 0.5 * \u))});
        };

        \foreach \v in {-0.6, -0.2, ..., 0.61} {
          \path [draw] (0, 0) -- (1 + \a * \v, \b * \v);
        }
        \foreach \u in {0.1, 0.2, ..., 0.91} {
          \path [draw] (\u + \a * \u, \b * \u) -- (\u - \a * \u, - \b * \u);
        }
        \path [draw, thick] (0, 0) -- (1+\a, \b) -- (1-\a, -\b) -- cycle;

        \begin{scope}[dotted]
          \path [draw] (0, 0) -- (0, {0.4 * (\h - 1)});
          \path [draw] (1 + \a, \b) -- (1 + \a, \b + \c * \h);
          \path [draw] (1 - \a, -\b) -- (1 - \a, -\b + \c * \h);
        \end{scope}

        \draw[->] (0.5, 1.8 * \c) -- (0.5, \c)
        node [right, midway] {$q$};
      }
    \end{tikzpicture}
  \end{minipage}
  \caption{Two views of the quotient of a square by its left edge.}
  \label{square-edge-quotient}
\end{figure}

\begin{example}
  \label{not-closed-by-closed-pushout}
  Here is an example of a closed-by-proper pushout in $\DD$
  which is not a closed-by-closed pushout.
  This example will recur in various guises in subsequent sections.

  Let $B = [0, 1]^2$ and $A = \{0\} \times [0, 1] \subset B$,
  and take $X = *$,
  with $j : A \to B$ the inclusion and $f : A \to X$ the unique map.
  If we can construct any (automatically proper) maps
  $q : B \to Y$ and $k : X \to Y$
  forming a square whose underlying sets form a pushout square in $\Set$,
  then we will have built a pushout in $\DD$ as well,
  by an argument similar to the one used to prove
  \cref{proper-quotient-recognition}.
  A simple choice is $Y = \{\,(x,y) \in R^2 \mid 0 \le y \le x \le 1\,\}$
  and $q(x,y) = (x, xy)$ with $k(*) = (0, 0)$.
  Indeed, $q$ is surjective
  and for $(x_1, y_1)$, $(x_2, y_2) \in B$,
  we have $q(x_1, y_1) = q(x_2, y_2)$
  if and only if $(x_1, y_1) = (x_2, y_2)$ or $x_1 = x_2 = 0$.

  \Cref{square-edge-quotient} shows
  two pictures of the map $q : B = [0,1]^2 \to Y$.
  On the left, we see that $q$ compresses
  each vertical segment $\{x\} \times [0,1]$
  to a vertical segment $\{x\} \times [0,x]$.
  On the right, we view $q$ as defining
  a family of spaces~$B_y$ indexed by $y \in Y$.
  This family consists of single points for all $y \ne (0, 0)$,
  but over the ``special point'' $y = (0, 0)$
  the family degenerates to an entire line segment.
  To obtain the domain of this model of $q$,
  take a flat rectangular strip and rotate its left edge by a quarter-twist.

  This example is a variation of the map $\pi$ from \cite[Example~8.2.10]{vdD98},
  which is the quotient of the square $[-1, 1]^2$
  by the bisecting segment $\{0\} \times [-1, 1]$
  and can be visualized by giving a rectangular strip an entire half-twist.
\end{example}

\begin{remark}
  \label{multiplication-required}
  Most of the theory of o-minimal topology
  we have introduced in previous sections
  is valid even when the o-minimal structure
  is only assumed to be an expansion of an ordered abelian group.
  However, \cref{exists-proper-quotient} on the existence of proper quotients
  relies in an essential way on the definability of \emph{multiplication}
  in the o-minimal structure.
  Using the above example, we illustrate what goes wrong in its absence.

  Recall the semilinear structure from \cref{omin-examples}:
  the smallest structure containing the constants of $R$
  and the graphs of the $<$ relation, addition,
  and scalar multiplication by each element of $R$.
  It is closely related to piecewise-linear geometry over $R$.
  In fact, \cite[8.2.14 Exercise~2]{vdD98} implies that
  any polytope $X$ in the semilinear structure \emph{is}
  (and not just is definably homeomorphic to)
  the underlying set $|K|$ of
  some finite geometric simplicial complex $K$ in $R$.
  If $f : X \to R^n$ is a continuous semilinear function,
  then by applying this fact to the graph of $f$,
  we see that $f$ is affine on each closed simplex of
  some complex $K$ with $|K| = X$.

  Now, working in the semilinear structure,
  consider the problem of finding
  a quotient of $B = [0,1]^2$ by its closed subset $A = \{0\} \times [0,1]$.
  Call a map $r : B \to Y$ a ``candidate''
  if $r$ is constant on the subset $A \subset B$,
  and for $\varepsilon > 0$ say that
  $r$ ``identifies points above $\varepsilon$''
  if there exist distinct points $(x_1, y_1)$, $(x_2, y_2) \in B$
  with $x_1 > \varepsilon$, $x_2 > \varepsilon$
  such that $r(x_1, y_1) = r(x_2, y_2)$.
  Then we claim:
  \begin{enumerate}
  \item[1.]
    For any $\varepsilon > 0$, there exists a candidate $r$
    which does not identify points above $\varepsilon$.
  \item[2.]
    For any candidate $r$, there exists $\varepsilon > 0$
    such that $r$ identifies points above $\varepsilon$.
  \end{enumerate}
  A quotient $q : B \to B/A$ is a candidate
  such that any candidate $r : B \to Y$ factors uniquely through $q$.
  But this is impossible in view of claims~1 and~2,
  because $q$ identifies points above some $\varepsilon > 0$,
  and then a candidate $r$ which does not identify points above $\varepsilon$
  cannot factor through $q$.

  Claim 1 is easily verified by writing down an explicit piecewise-linear map
  which collapses the triangle
  $\{\,(x, y) \mid y \ge x / \varepsilon\,\} \subset B$
  down to its diagonal edge $y = x / \varepsilon$
  and leaves the rest of $B$ fixed.
  For claim 2, we use the fact cited above that
  $r$ is affine on each simplex of a triangulation of the set $B \subset R^2$.
  Pick a point $p = (0, y)$ which is not a vertex of this triangulation.
  It must lie on the boundary of some 2-simplex $\sigma$,
  which must have two of its vertices on the edge $A = \{0\} \times [0, 1]$.
  Then $r$ maps these vertices to the same point,
  so because it is affine on $\sigma$,
  it also identifies some pairs of points with positive $x$-coordinates.

  Consequently, if we were to work with the semilinear structure $\sS = R_\slin$,
  \cref{closed-by-proper-pushout} would fail,
  and $\PP$ would not be a pretopos.
\end{remark}

\begin{nul}
  \label{pretopos-morphism}
  Let $\Phi : (R, \sS) \to (R', \sS')$ be a morphism of o-minimal structures.
  We claim that the induced functor $\Phi : \DD(R, \sS) \to \DD(R', \sS')$
  preserves finite limits, finite coproducts and
  quotients of proper equivalence relations.
  For finite limits and finite coproducts
  this is obvious by inspection of their constructions.
  Suppose $E \rightrightarrows X$ is a proper equivalence relation on $X$
  with quotient $q : X \to X/E$;
  then $\Phi E \rightrightarrows \Phi X$ is also
  a proper equivalence relation on $\Phi X$,
  and we must show that
  $\Phi E \rightrightarrows \Phi X \xrightarrow{\Phi q} \Phi (X/E)$
  is a coequalizer.
  The map $\Phi q$ is again proper and surjective
  and the square
  \[
    \begin{tikzcd}
      \Phi E \ar[r, "\Phi p_1"] \ar[d, "\Phi p_2"'] & \Phi X \ar[d, "\Phi q"] \\
      \Phi X \ar[r, "\Phi q"] & \Phi (X/E)
    \end{tikzcd}
  \]
  is a pullback,
  being the image under $\Phi$ of the corresponding square for $q$.
  Then by \cref{coequalizer-of-proper-surjective},
  $\Phi q$ is the coequalizer of its kernel pair
  $\Phi E \rightrightarrows \Phi X$.
  The functor $\Phi$ restricts to a functor $\Phi : \PP(R, \sS) \to \PP(R', \sS')$
  which is a pretopos morphism.
\end{nul}

\begin{nul}
  \label{quotients-top}
  In the preceding paragraph we considered
  the variation of $\DD$ in the o-minimal structure.
  When $R = \RR$, we can similarly compare $\DD$ to the category $\Top$.
  The functor ${-}_\tp : \DD \to \Top$ again
  preserves finite coproducts and finite limits by inspection.
  If $q : X \to Y$ is a proper surjective map of $\DD$
  then $q_\tp$ is also proper, by \cref{proper-iff-rr}.
  In particular, $q_\tp : X_\tp \to Y_\tp$ is closed
  and so $q_\tp$ induces a homeomorphism $X_\tp / Q \cong Y_\tp$,
  where $Q = (X \times_Y X)_\tp$ is
  the relation on $X_\tp$ defined by $Q(x, y)$ if $q(x) = q(y)$.
  (See \cite[I.5.3]{TG}.)
  Then $Y_\tp$ is the coequalizer of the maps
  $(X \times_Y X)_\tp \rightrightarrows X_\tp$,
  by the universal property of the quotient topology.
  Hence, if $q$ was the quotient of $X$
  by a proper equivalence relation $E \rightrightarrows X$,
  then the coequalizer $E \rightrightarrows X \to X/E$
  is preserved by the functor ${-}_\tp$.
  In particular, the restriction ${-}_\tp : \PP \to \CompHaus$
  is a pretopos morphism.

  As a special case, ${-}_\tp$ preserves closed-by-proper pushouts,
  because they can be constructed from
  coproducts and quotients by proper equivalence relations.
  (But there are other colimits in $\DD$ that it fails to preserve,
  like those of \cref{bad-colimits}.)
\end{nul}

\begin{warning}
  \label{definably-closed-not-closed}
  When $R \ne \RR$,
  a definably closed morphism $q : X \to Y$ of $\DD$
  is not necessarily closed as a map of underlying topological spaces,
  because $X$ has more closed subsets than closed definable subsets.
  (This was not a problem for continuity
  because continuity can be checked on a base for
  the topology of the target space.)
  For example, take $R = \Qalgre$,
  let $q : [0, 1]^2 \to [0, 1]$ be
  the projection onto the first coordinate,
  and pick $\alpha \in \RR - \Qalgre$ with $0 < \alpha < 1$.
  The topological space $[0, 1]^2_\tp = ([0, 1] \cap \Qalgre)^2$
  has a closed subset
  $V = \{\,(x, y) \mid x > |y - \alpha|\,\}$
  whose image $q(V) = (0, 1] \cap \Qalgre$ is not closed.
  Indeed, there are no points $(x, y) \in ([0, 1] \cap \Qalgre)^2$
  with $x = |y - \alpha|$,
  so the complement of $V$
  is the open set $\{\,(x, y) \mid x < |y - \alpha|\,\}$.
  But for every $x > 0$, there exists $y \in [0, 1] \cap \Qalgre$
  with $|y - \alpha| < x$.

  This example shows that the argument of \cref{quotients-top}
  relies on the hypothesis that $R = \RR$.
  However, in this example the map $q_\tp$ is clearly still a quotient map.
  We do not know in general whether
  $q_\tp$ is a quotient map when $q$ is proper and surjective and $R \ne \RR$.
\end{warning}

\begin{example}
  \label{closed-by-improper}
  We give an example of a closed embedding and an open embedding
  which admit no pushout in $\DD$,
  showing that in \cref{closed-by-proper-pushout},
  we cannot omit the hypothesis that $j$ is closed or that $f$ is proper.
  This example is closely related to
  Example~1.8 of \cite{Br87} (see the discussion after Lemma~2.6 there),
  but we give a self-contained proof of the nonexistence of the pushout.

  First, suppose
  \[
    \begin{tikzcd}
      A \ar[r, "f"] \ar[d, "j"'] & X \ar[d, "j'"] \\
      B \ar[r, "f'"] & Y
    \end{tikzcd}
  \]
  is any pushout in $\DD$ whatsoever.
  We claim that the maps $j'$ and $f'$ must be jointly surjective.
  Let $y \in Y$.
  Then (regarding $Y$ as embedded in some~$R^n$)
  the function $g = d(y, -)$ is a continuous definable map $g : Y \to R$
  with $g^{-1}(\{0\}) = \{y\}$.
  If $y$ is not in the image of either $f'$ or $j'$,
  then the restrictions $gf'$ and $gj'$ factor through the subspace $(0, \infty) \subset R$,
  and then by the universal property of the pushout,
  $g$ must factor through $(0, \infty)$ as well, which is impossible.

  Next, suppose $Z$ is a definable set
  and $B$ and $X$ are definable subsets of $Z$ with $B \cup X = Z$.
  Let $A = B \cap X$, and suppose there is a pushout $Y$ as above
  with $f$ and $j$ the inclusion maps.
  Then, there is an induced comparison map $k : Y \to Z$.
  On the other hand the underlying set $UZ$
  is the pushout of the diagram $UB \leftarrow UA \to UX$
  and so there is an induced map $l : UZ \to UY$.
  Using the fact that $f'$ and $j'$ are jointly surjective,
  it is not hard to see that $Uk$ and $l$ are mutually inverse,
  so set-theoretically we may identify $Y$ with~$Z$.
  However, $Y$ could be embedded in some $R^n$ in a different way,
  giving it a different ``definable topology''.
  For example, if $B$ and $X$ are disjoint subsets of $Z$
  then $Y = B \amalg X$ is generally not homeomorphic to $Z = B \cup X$.

  Now, in $Z = \{(0, 0)\} \cup (0, 1) \times [0, 1)$,
  take
  \[
    B = [0, 1) \times \{0\}, \quad
    X = (0, 1) \times [0, 1), \quad
    A = (0, 1) \times \{0\}.
  \]
  We claim the inclusion maps $f : A \to X$ and $j : A \to B$ admit no pushout.
  Suppose $Y$ is a pushout, as above.
  We will probe the topology of $Y$ using curves.

  \begin{figure}
    \[
      \tikz{
        \draw (0, 0) node[left] {$A = {}$};
        \draw (0, -1) -- (2, -1);
        \filldraw[fill=white] (0, -1) circle (0.05);
        \filldraw[fill=white] (2, -1) circle (0.05);

        \draw[->] (2.5, 0) -- node[above] {$f$} (3.5, 0);

        \draw (5, 0) node[left] {$X = {}$};
        \draw (5, -1) -- (7, -1);
        \draw[dashed] (7, -1) -- (7, 1) -- (5, 1) -- (5, -1);
        \draw (6, 0) node {$U$};
        \filldraw[fill=white] (5, -1) circle (0.05);
        \filldraw[fill=white] (7, -1) circle (0.05);

        \draw[->] (1, -1.5) -- node[left] {$j$} (1, -2.5);
        \draw[->] (6, -1.5) -- node[right] {$j'$} (6, -2.5);

        \draw (0, -4) node[left] {$B = {}$};
        \draw (0, -5) -- (2, -5);
        \filldraw[fill=black] (0, -5) circle (0.05) node[below] {$b_0$};
        \filldraw[fill=white] (2, -5) circle (0.05);

        \draw[->] (2.5, -4) -- node[above] {$f'$} (3.5, -4);

        \draw (5, -4) node[left] {$Y = {}$};
        \draw (5, -5) .. controls (5.75, -4.5) and (6.25, -5.5) ..
        node[below]{$L$} (7, -5);
        \draw[dashed] (7, -5) .. controls (7.5, -2.5) and (5.5, -2.5) .. (5, -5);
        \draw (6.25, -4) node {$W$};
        \filldraw[fill=black] (5, -5) circle (0.05) node[below] {$y_0$};
        \filldraw[fill=white] (7, -5) circle (0.05);
      }
    \]
    \caption{A nonexistent pushout in $\DD$.}
  \end{figure}

  Let $b_0 \in B$ denote the point $(0, 0)$
  and $y_0 \in Y$ its image under $f'$,
  and let $L \subset Y$ denote the image of $A$
  and $W \subset Y$ the image of $U := X - A$.
  Then $f'(B) = \{y_0\} \cup L$ and $j'(X) = Y - \{y_0\} = L \cup W$.
  We claim that $y_0$ belongs to the closure of $W$.
  We prove this in a series of two claims.
  \begin{enumerate}
  \item
    $y_0$ belongs to the closure of $L$.

    Consider a curve $\gamma : [0, \varepsilon] \to B$
    with $\gamma(0) = b_0$ and $\gamma(t) \ne b_0$ for $t > 0$.
    The image of $\gamma$ under $f'$ is a curve $\gamma'$
    with $\gamma'(0) = y_0$ and $\gamma'(t) \in L$ for $t > 0$,
    and therefore $y_0 \in \overline{L}$.
  \item
    $L$ is contained in the closure of $W$.

    Each point of $L$ is the image under $f'$ of a point $(x, 0) \in X$.
    Then we choose a curve $\gamma : [0, \varepsilon] \to X$
    with $\gamma(0) = (x, 0)$ and $\gamma(t) \in U$ for $t > 0$.
    By the same argument as before, the original point of $L$
    belongs to the closure of $W$.
  \end{enumerate}

  By curve selection, then,
  we may choose a curve $\gamma : [0, \varepsilon] \to Y$
  such that $\gamma(0) = y_0$ and $\gamma(t) \in W$ for $t > 0$.
  The map $j'|_U : U \to W$ is a definable bijection,
  so after reducing $\varepsilon$ if necessary,
  we can factor the restriction of $\gamma$ to $(0, \varepsilon]$ through $U$,
  obtaining a curve $\tilde \gamma : (0, \varepsilon] \to U \subset X$
  such that $j' \circ \tilde \gamma$ converges to $y_0$ at $0$.
  The image of $\tilde \gamma$ is closed in $X$:
  this could only fail if $\tilde \gamma$ is completable at $0$ to a point of $X$,
  but that would contradict the fact that $y_0$ is not in the image of $j'$.

  Now, $A$ and the image of $\tilde \gamma$
  are disjoint closed subsets of $X$,
  so we can choose a continuous definable function $h : X \to R$
  which restricts to $0$ on $A$ and to $1$ on the image of $\tilde \gamma$.
  By the universal property of $Y$,
  there is a function $h' : Y \to R$ whose restriction to $X$ is $h$
  and whose restriction to $B$ is $0$.
  In particular $h'(y_0) = 0$.
  But this is a contradiction,
  because we can also compute $h'(y_0)$ as the limit of $h'$ along the curve $\gamma$,
  and by construction $h'$ takes the value $1$ at all points $\gamma(t)$ with $t > 0$.
\end{example}

\begin{example}
  A similar sort of argument shows that
  there is no quotient $X/A$ in $\DD$ for, say,
  $X = R^2$ and $A \subset X$ the $x$-axis;
  this does not contradict \cref{closed-by-proper-pushout}
  because the map $A \to *$ is not proper.
  Here the idea is that if $q : X \to X/A$ is such a quotient
  and $g : X/A \to [0, \infty)$ is the distance in $X/A$ to the quotient basepoint,
  then $(x, y) \mapsto xg(q(x, y))$ cannot descend to a continuous function on $X/A$.
  The point is that the function $x$ is unbounded
  on points of $X$ ``near'' (but not in) $A$.
  More precisely, we can find curves $\gamma$
  with $\gamma(0) \in A$ and $\gamma(t) \in X - A$ for $t > 0$
  on which the continuous function~$x$ takes values as large as we wish.
  This could not occur if $A$ was definably compact.
\end{example}

\begin{remark}
  \label{no-infinite-colimits}
  Finally, while the small categories $\DD$ or $\PP$
  clearly cannot admit general infinite colimits,
  we caution the reader that
  even some diagrams in $\DD$ which may appear to be infinite colimits
  actually are not,
  because of the finitary nature of definability.
  For example, suppose that $R = \RR$.
  Then as a set, and even as a topological space,
  $\RR$ is the directed union of its open subspaces $(-N, N)$ for $N \ge 1$.
  However, the diagram
  \[
    (-1, 1) \to (-2, 2) \to \cdots \to \RR
  \]
  is not a colimit diagram in $\DD$,
  because we can find a function $f : \RR \to [0, 1]$
  which is continuous and definable when restricted to each $(-N, N)$,
  but not globally definable on $\RR$.
  For example, $f(x)$ could be the distance from $x$ to the nearest integer.
\end{remark}

%% file: tex/realization.tex
\section{The geometric realization of a finite simplicial set}
\label{geometric-realization}

\begin{nul}
  \label{definable-simplex}
  The classical geometric realization functor
  $|{-}|_\Top : \sSet \to \Top$
  is the (essentially unique) colimit-preserving extension
  of the ``topological simplex functor''
  $|\Delta^\bullet|_\Top : \Delta \to \Top$.
  It can be constructed explicitly by the formula
  \[
    |K|_\Top = \colim_{\sigma : \Delta^n \to K} |\Delta^n|_\Top
    \tag{$*$}
  \]
  and characterized as the left adjoint to the functor
  $\Sing^\Top : \Top \to \sSet$
  defined by
  \[
    (\Sing^\Top X)_n = \Hom_\Top(|\Delta^n|_\Top, X).
  \]

  In the o-minimal setting, we also have a ``definable simplex functor''
  $|\Delta^\bullet|_\DD : \Delta \to \DD$.
  On objects, it can be defined as
  \[
    |\Delta^n|_\DD =
    \{\,(t_1, \ldots, t_n) \mid 0 \le t_1 \le \cdots \le t_n \le 1\,\}
    \subset R^n.
  \]
  The face and degeneracy maps are given by the usual formulas,
  which involve only copying and deleting variables
  and the constants $0$ and $1$.
  The $i$th vertex $e^n_i \in |\Delta^n|_\DD$,
  which is the image of the map $|\Delta^0|_\DD \to |\Delta^n|_\DD$
  induced by the map $[0] \to [n]$ picking out $i \in [n]$,
  is the vertex with the last $i$ coordinates $t_{n-i+1}$, \ldots, $t_n$ equal to $1$,
  and $|\Delta^n|_\DD$ is the convex hull of $e^n_0$, \ldots, $e^n_n$ in $R^n$.
  We identify $|\Delta^1|_\DD$ with the unit interval $[0, 1] \subset R$,
  with $e^1_0 = 0$ and $e^1_1 = 1$.

  Each set $|\Delta^n|_\DD$ is closed and bounded,
  so we may also regard $|\Delta^\bullet|_\DD$ as taking values in $\PP$,
  in which case we denote it by $|\Delta^\bullet|_\PP$.
  When $R = \RR$, the composition
  $\Delta \xrightarrow{|\Delta^\bullet|_\DD} \DD \xrightarrow{{-}_\tp} \Top$
  equals $|\Delta^\bullet|_\Top$.

  As in the topological case,
  the functor $|\Delta^\bullet|_\DD : \Delta \to \DD$
  induces a functor $\Sing^\DD : \DD \to \sSet$
  defined by
  \[
    (\Sing^\DD X)_n = \Hom_\DD(|\Delta^n|_\DD, X).
  \]
  However, as $\DD$ is not cocomplete,
  we cannot expect this functor $\Sing^\DD$ to have a left adjoint.
  In fact, any definable set has only finitely many path components,
  so we cannot reasonably hope to form
  a $\DD$-valued geometric realization of
  a simplicial set with infinitely many connected components.
\end{nul}

\begin{nul}
  \label{realization-outline}
  The goal of this section is to show that
  there is a well-behaved $\DD$-valued geometric realization functor $|{-}|_\DD$
  defined on the full subcategory of \emph{finite} simplicial sets.
  This functor $|{-}|_\DD : \sSetfin \to \DD$ is not a left adjoint,
  because $\Sing^\DD$ does not take values in finite simplicial sets.
  Instead, $|{-}|_\DD$ is a left adjoint of $\Sing^\DD$
  \emph{relative} to the inclusion $\sSetfin \subset \sSet$,
  as we will explain in \cref{relative-adjoint}.

  To construct the $\DD$-valued geometric realization
  of a finite simplicial set~$K$,
  we should form some sort of colimit of the objects $|\Delta^n|_\DD$
  indexed by the simplices $\Delta^n \xrightarrow{\sigma} K$ of $K$.
  However, as we saw in the previous section,
  only certain kinds of colimits in $\DD$ are well-behaved,
  so some care is needed.
  Roughly speaking, we will proceed as follows.
  \begin{enumerate}
  \item
    By definition,
    the geometric realization of a single simplex $\Delta^n$
    is the object $|\Delta^n|_\DD$.
  \item
    Suppose $K$ is a simplicial subset of $\Delta^n$ for some $n \ge 0$;
    we call such a $K$ (together with its embedding into a $\Delta^n$)
    a \emph{simplicial complex}.
    Then $|K|_\DD$ can be constructed as the union of
    the corresponding faces of $|\Delta^n|_\DD$.
    We prove this by induction on the number of simplices contained in $K$,
    using (\textit{i}) and the fact that
    the union of two closed subsets of $|\Delta^n|_\DD$
    can be built as a closed-by-closed pushout.
  \item
    Step (\textit{ii}) applies in particular to
    the simplicial complexes $\partial \Delta^n \subset \Delta^n$
    for $n \ge 0$,
    so we can identify $|\partial \Delta^n|_\DD$
    with a closed subset of $|\Delta^n|_\DD$
    (namely, the union of its proper faces).
    Then, we can build up a general finite simplicial set $K$
    by attaching a simplex along its boundary
    a finite number of times,
    and so we can construct the geometric realization $|K|_\DD$
    using a corresponding sequence of closed-by-proper pushouts in $\DD$.
  \end{enumerate}

  This approach to constructing the topological geometric realization functor
  is described in detail in \kerodon{001X}.
  We borrow much of its terminology and notation,
  and omit some details from proofs which have analogues there.
\end{nul}

\begin{nul*}
  To make our goal more precise, we begin with the following definitions.
\end{nul*}

\begin{definition}
  \label{exhibits-as-realization}
  Let $\C$ be a category equipped with a functor $Q : \Delta \to \C$.
  We define the functor $\Sing^Q : \C \to \sSet$
  by the formula $(\Sing^Q X)_n = \Hom_C(Q[n], X)$.

  Let $K$ be a simplicial set and $Y$ an object of $\C$.
  Then we say a map $u : K \to \Sing^Q Y$
  \emph{exhibits $Y$ as a $Q$-realization of $K$}
  if, for each object~$Z$ of~$\C$,
  the composition
  \[
    \Hom_\C(Y, Z) \xrightarrow{\Sing^Q} \Hom_\sSet(\Sing^Q Y, \Sing^Q Z)
    \xrightarrow{- \circ u} \Hom_\sSet(K, \Sing^Q Z)
  \]
  is an isomorphism.
  More informally, we sometimes say that
  $Y$ ``is'' the $Q$-realization of $K$,
  with the map $u : K \to \Sing^Q Y$ implicit.
  We say $K$ \emph{admits a $Q$-realization}
  if such a pair $(Y, u)$ exists.

  Suppose now that $F : \C \to \C'$ is
  a functor to a second category $\C'$.
  Then there is a natural transformation
  $\Sing^Q ({-}) \to \Sing^{F \circ Q} F({-})$
  whose action on simplices is given by the action of $F$ on morphisms.
  Let $K$ be a simplicial set which admits a $Q$-realization $(Y, u)$.
  We say that $F$ \emph{preserves the $Q$-realization of $K$}
  if the composite $K \xrightarrow{u} \Sing^Q Y \to \Sing^{F \circ Q} FY$
  exhibits $FY$ as the $(F \circ Q)$-realization of $K$.
\end{definition}

\begin{nul}
  \label{realization-as-colimit}
  Informally, $Y$ is a $Q$-realization of $K$
  when $Y$ represents the functor $\Hom_\sSet(K, \Sing^Q(-))$.
  We can also reformulate this condition as follows.
  Giving a map $u : K \to \Sing^Q Y$
  amounts to specifying for each simplex $\sigma : \Delta^n \to K$ of $K$
  a ``characteristic map'' $u_\sigma : Q[n] \to Y$
  in a compatible way.
  The compatibility conditions amount to the statement that
  the maps $u_\sigma$ form a cocone with vertex $Y$
  on the diagram indexed on the category of simplices of $K$
  \kerodon{00X4}
  sending a simplex $([n], \sigma : \Delta^n \to K)$
  to the object $Q[n]$ of $\C$.
  Then $u$ exhibits $Y$ as the $Q$-realization of $K$
  precisely when this cocone is a colimit cocone.
  A functor $F : C \to C'$ preserves the $Q$-realization $u$
  when it preserves this colimit.
\end{nul}

\begin{nul}
  We speak of ``geometric realizations''
  when $Q$ is one of the standard simplex functors
  taking values in $\DD$ (or $\PP \subset \DD$) or $\Top$.
  In this case we also write the target category
  as a subscript or superscript in place of $Q$.

  From the outline in \cref{realization-outline}
  one can see that the $\DD$-valued geometric realization of a finite simplicial set $K$
  will in fact always be an object of $\PP$.
  Since $\PP \subset \DD$ is a full subcategory,
  it follows immediately from \cref{exhibits-as-realization}
  that the $\DD$-valued geometric realization of $K$
  is also a $\PP$-valued geometric realization of $K$
  (i.e., a realization with respect to the functor $|\Delta^\bullet|_\PP : \Delta \to \PP$).
  Thus, in the end it makes no real difference
  whether we take $\PP$ or $\DD$ as the target category of geometric realization.
\end{nul}

\begin{nul}
  \label{simplex-realization}
  The object $K = \Delta^n$ always admits a $Q$-realization for any $Q : \Delta \to C$:
  we take $Y = Q[n]$
  and $u : \Delta^n \to \Sing^Q Q[n]$
  the map corresponding to
  the identity element of $(\Sing^Q Q[n])_n = \Hom(Q[n], Q[n])$.
  Furthermore, this $Q$-realization is preserved by any functor.
  Thus, if we denote the functor $Q$ by notation like $|\Delta^\bullet|_\C$,
  there is no conflict between $|\Delta^n|_\C = Q[n]$
  and the $Q$-realization $|\Delta^n|_Q$.
\end{nul}

\begin{example}
  \label{realization-example-explicit}
  The simplicial set $\Lambda^2_1 = \Delta^1 \amalg_{\Delta^0} \Delta^1$
  consists of two $1$-simplices glued end-to-end.
  Let $Y = [0, 2] \subset R$
  and let $u : \Lambda^2_1 \to \Sing^\DD Y$ be the map
  sending the first $\Delta^1$ of $\Lambda^2_1$
  to the inclusion map $|\Delta^1|_\DD = [0, 1] \to [0, 2]$,
  and the second $\Delta^1$ to the map $[0, 1] \to [0, 2]$
  given by $t \mapsto t+1$.
  Then $u$ exhibits $Y$ as the $\DD$-valued geometric realization of $\Lambda^2_1$,
  because a continuous definable map from $Y = [0, 2]$ to some object $Z \in \DD$
  is uniquely determined by two continuous definable maps
  $[0, 1] \to Z$ and $[1, 2] \to Z$ that agree at the point~$1$,
  hence corresponds uniquely via $u$ to a map $\Lambda^2_1 \to \Sing^\DD Z$.

  Note that, just like in $\Top$,
  the $\DD$-valued geometric realization of $\Lambda^2_1$
  is (definably) homeomorphic to a single interval.
  This is what lies behind the ability to concatenate paths by the usual formula.
  In general, finitary elementary constructions of this kind in the homotopy theory of spaces
  work the same way in the definable setting.
\end{example}

\begin{nonexample}
  Take $R = \RR$
  and let $K$ be the simplicial set
  whose $0$-simplices are the integers
  and with a nondegenerate $1$-simplex from $n$ to $n+1$
  for each integer $n$.
  Then
  it may appear that (as in topology) the map $u : K \to \Sing^\DD \RR$
  defined as follows exhibits $Y = \RR$ as a $\DD$-valued realization of $K$:
  \begin{itemize}
  \item
    $u$ sends the $0$-simplex $n$
    to the map $|\Delta^0|_\DD = * \to \RR$ with value $n$.
  \item
    $u$ sends the $1$-simplex from $n$ to $n+1$
    to the map $|\Delta^1|_\DD = [0, 1] \to \RR$ sending $t$ to $n+t$.
  \end{itemize}
  However, this fails because
  $\RR$ does not have enough definable continuous functions.
  In \cref{exhibits-as-realization}, take $Z = \RR$.
  If $u$ exhibited $\RR$ as a $\DD$-valued realization of $K$,
  then every map $f : K \to \Sing^\DD \RR$ would arise from a map $g : \RR \to \RR$
  by restriction to the intervals $[n, n+1]$ for each integer $n$.
  But the function $f$ sending each $1$-simplex of $K$
  to the map $\psi : [0, 1] \to \RR$, $\psi(t) = \min \{t, 1-t\}$
  does not arise in this way from any definable continuous function $g : \RR \to \RR$,
  since $g$ would have to have infinitely many isolated zeros.
  (Compare \cref{no-infinite-colimits}.)

  More generally, it seems likely that
  no infinite simplicial set has a $\DD$-valued (or $\PP$-valued) geometric realization,
  but we have not proved this.
\end{nonexample}

\begin{nul}[Functoriality]
  \label{realization-functoriality}
  The construction outlined in \cref{realization-outline}
  is not obviously functorial in the finite simplicial set $K$.
  However, the universal property that defines a $Q$-realization of $K$
  also makes it automatically functorial in~$K$.
  Before we turn to the construction of geometric realizations,
  we first recall this functoriality of $Q$-realizations in general
  and their description as a relative left adjoint to the functor $\Sing^Q$.

  Fix a category $C$ equipped with a functor $Q : \Delta \to C$.
  The $Q$-realization of a simplicial set $K$,
  when it exists,
  is essentially unique by the Yoneda lemma.
  Furthermore, suppose $f : K \to K'$
  is a map between two simplicial sets
  which admit $Q$-realizations $u : K \to \Sing^Q Y$ and $u' : K' \to \Sing^Q Y'$.
  Then there is a unique map $g : Y \to Y'$ such that the diagram
  \[
    \begin{tikzcd}
      K \ar[r, "u"] \ar[d, "f"'] & \Sing^Q Y \ar[d, "\Sing^Q g"] \\
      K' \ar[r, "u'"] & \Sing^Q Y'
    \end{tikzcd}
  \]
  commutes,
  obtained by applying the defining property of $Y'$ to the map $u' f$.
  We call $g$ the $Q$-realization of $f$
  (with respect to the chosen $Q$-realizations $u$ and $u'$ of $K$ and $K'$).

  Suppose then that $S \subset \sSet$
  is a full subcategory consisting of objects
  each of which admits some $Q$-realization.
  Choose a $Q$-realization $u_K : K \to \Sing^Q Y_K$ of each $K \in S$,
  and define $|K|_Q$ to be the object $Y_K$.
  This defines the action on objects of a functor $|{-}|_Q : S \to C$
  whose sends a morphism $f : K \to K'$
  to its $Q$-realization $g : Y_K \to Y_{K'}$ as described in the previous paragraph.
  We call this the $Q$-realization functor (on $S$).
  The maps $u_K : K \to \Sing^Q Y_K$ then assemble into
  a natural transformation $u : J \to \Sing^Q |{-}|_Q$,
  where $J : S \to \sSet$ denotes
  the inclusion of the full subcategory $S$.

  The $Q$-realization functor can be characterized
  using the language of \emph{relative adjoints} \cite{Ulmer}.
\end{nul}

\begin{definition}
  \label{relative-adjoint}
  Let $A$, $B$, and $C$ be categories and $J : A \to B$ a functor.
  (For us, the functor $J$ will always be fully faithful.)
  We say that a functor $F : A \to C$ is \emph{left adjoint relative to $J$}
  to a functor $U : C \to B$
  when we are given a specified isomorphism
  \[
    \Hom_C(FX, Y) \cong \Hom_B(JX, UY)
    \tag{$\dagger$}
  \]
  natural in $X \in A$ and $Y \in C$.
\end{definition}

\begin{nul}
  When $A = B$ and $J$ is the identity functor,
  an adjunction relative to $J$ is the same as an ordinary adjunction $F \dashv U$.
  In the general case,
  the right adjoint $U$ of a relative adjunction
  determines the left adjoint $F$ uniquely,
  but not vice versa;
  this is because ($\dagger$) determines all maps out of objects of the form $FX$,
  but only determines the maps into objects of the form $UY$
  from objects in the image of $J$.
\end{nul}

\begin{nul}
  The preceding discussion is related to relative adjunctions as follows.
  Fix $C$, $Q$, and $S \subset \sSet$ as above.
  Then the functor $\Sing^Q : C \to \sSet$ has
  a left adjoint relative to $J : S \to \sSet$
  if and only if each object of $S$ admits a $Q$-realization,
  and in this case the relative left adjoint is
  the $Q$-realization functor $|{-}|_Q : S \to C$.
\end{nul}

\begin{example}
  Suppose the category $C$ is cocomplete (e.g., $C = \Top$).
  Then every simplicial set admits a $Q$-realization,
  so we may take $S$ to be all of $\sSet$,
  and the above notions reduce to the usual ones:
  the $Q$-realization functor $|{-}|_Q$ is an ordinary left adjoint to $\Sing^Q$.
  Any colimit-preserving functor also preserves all $Q$-realizations
  (by \cref{realization-as-colimit}).
\end{example}

\begin{nul*}
  Now, we turn to the construction of
  geometric realizations of finite simplicial sets.
  First, we make precise the construction appearing in
  step (\textit{ii}) of \cref{realization-outline}.
\end{nul*}

\begin{construction}
  \label{complex-realization-construction}
  Let $n \ge 0$.
  A face $\sigma$ of $\Delta^n$ contains some nonempty subset
  $V(\sigma) \subset \{0, \ldots, n\}$ of the vertices of $\Delta^n$.
  We let $|\sigma|_\DD$ denote the convex hull of the vertices $e^n_i$ of $|\Delta^n|_\DD$
  (see \cref{definable-simplex})
  as $i$ ranges over $V(\sigma)$.
  If the face $\sigma$ is the image of the injective map $\varphi : \Delta^k \to \Delta^n$
  (with $k + 1 = |V(\sigma)|$)
  then $|\sigma|_\DD$ is also the image of
  the induced map $|\Delta^\varphi|_\DD : |\Delta^k|_\DD \to |\Delta^n|_\DD$.
  The restriction of
  the canonical map $u_{\Delta^n} : \Delta^n \to \Sing^\DD |\Delta^n|_\DD$ to $\sigma$
  factors through the simplicial subset
  $\Sing^\DD |\sigma|_\DD \subset \Sing^\DD |\Delta^n|_\DD$,
  forming a square
  \[
    \begin{tikzcd}
      \sigma \ar[d] \ar[r, "u_\sigma"] & \Sing^\DD |\sigma|_\DD \ar[d] \\
      \Delta^n \ar[r, "u_{\Delta^n}"] & \Sing^\DD |\Delta^n|_\DD
    \end{tikzcd}
  \]
  and $u_\sigma$ exhibits $|\sigma|_\DD$ as
  the geometric realization of $\sigma$
  because $u_\sigma$ is isomorphic to the map
  $u_{\Delta^k} : \Delta^k \to \Sing^\DD |\Delta^k|_\DD$.

  Now let $K$ be a \emph{simplicial complex},
  that is,
  a simplicial subset of $\Delta^n$ for some $n \ge 0$.
  Then $K$ is the union of those faces of $\Delta^n$ which it contains.
  We let $|K|_\DD$ denote the union of the subsets $|\sigma|_\DD$
  as $\sigma$ ranges over the faces of~$\Delta^n$ contained in~$K$.
  The restriction of $u_{\Delta^n}$ to $K$
  again factors through the simplicial subset $\Sing^\DD |K|_\DD$,
  forming the square below.
  \[
    \begin{tikzcd}
      K \ar[d] \ar[r, "u_K"] & \Sing^\DD |K|_\DD \ar[d] \\
      \Delta^n \ar[r, "u_{\Delta^n}"] & \Sing^\DD |\Delta^n|_\DD
    \end{tikzcd}
  \]
  The notations $|K|_\DD$ and $u_K$ will be justified
  by \cref{simplicial-complex-realization}.
  Note that $|K|_\DD$ is a closed subset of $|\Delta^n|_\DD$,
  and so belongs to $\PP$.
\end{construction}

\begin{example}
  Suppose $K = \Lambda^2_1 \subset \Delta^2$.
  Then $K$ is the union of the two edges of $\Delta^2$ that contain $1$
  (along with the three vertices of $\Delta^2$, which will not affect the result).
  The set $|K|_\DD \subset |\Delta^2|_\DD \subset R^2$
  is the union of the segment from $e^2_0 = (0, 0)$ to $e^2_1 = (0, 1)$
  and the segment from $e^2_1 = (0, 1)$ to $e^2_2 = (1, 1)$,
  The map $u_K : K \to \Sing^\DD |K|_\DD$
  sends each $1$-simplex of $K$ to
  the evident corresponding map from $|\Delta^1|_\DD = [0, 1]$ into $|K|_\DD$.
  The map $(x_1, x_2) \mapsto x_1 + x_2$
  defines a homeomorphism from $|K|_\DD$ to $[0, 2]$
  and so by \cref{realization-example-explicit},
  $u_K$ does exhibit $|K|_\DD$ as the geometric realization of $K$.
\end{example}

\begin{nul*}
  We now come to the main results of this section.
  We can summarize these results by the slogan:
  while the piecewise-linear setting
  admits geometric realizations of finite simplicial complexes,
  the semialgebraic (or more generally o-minimal) setting
  admits geometric realizations of \emph{all} finite simplicial sets.
\end{nul*}

\begin{proposition}
  \label{simplicial-complex-realization}
  Let $K \subset \Delta^n$ be a simplicial complex.
  Then, with the notation of \cref{complex-realization-construction},
  the map $u_K : K \to \Sing^\DD |K|_\DD$
  exhibits $|K|_\DD$ as
  the $\DD$-valued (hence also $\PP$-valued) geometric realization of $K$.
  Moreover, as a $\PP$-valued geometric realization,
  this realization is preserved by any functor $F : \PP \to C'$
  satisfying the following two conditions:
  \begin{enumerate}
  \item[(a)] $F$ sends $\emptyset \in \PP$ to an initial object of $C'$.
  \item[(b$'$)] $F$ sends closed-by-closed pushouts in $\PP$ to pushouts in $C'$.
  \end{enumerate}
  If $K \subset K' \subset \Delta^n$,
  then the geometric realization of the inclusion $K \subset K'$
  is the inclusion $|K|_\DD \subset |K'|_\DD$
  between subsets of $|\Delta^n|_\DD$.
\end{proposition}

\begin{proposition}
  \label{finite-realizations}
  \label{realization-of-monomorphism}
  Any finite simplical set admits a $\PP$-valued realization.
  Moreover, these realizations are preserved by any functor $F : \PP \to C'$
  satisfying the following two conditions:
  \begin{enumerate}
  \item[(a)] $F$ sends $\emptyset \in \PP$ to an initial object of $C'$.
  \item[(b)] $F$ sends closed-by-proper pushouts in $\PP$ to pushouts in $C'$.
  \end{enumerate}
  In particular, these realizations are also $\DD$-valued realizations.
  The geometric realization of a monomorphism of finite simplicial sets
  is a closed embedding.
\end{proposition}

\begin{nul*}
  The proofs are based on the following general fact,
  which follows easily from \cref{exhibits-as-realization}.
  It is a minor extension of \kerodon{0023}.
\end{nul*}

\begin{lemma}
  \label{realization-colimit}
  Fix a functor $Q : \Delta \to C$
  and suppose $(K_i)_{i \in I}$ is a diagram of simplicial sets
  each of which admits a $Q$-realization $u_i : K_i \to \Sing^Q Y_i$.
  Suppose the induced diagram $(Y_i)_{i \in I}$ (as in \cref{realization-functoriality})
  has a colimit $Y$ in $C$, and write $K = \colim_{i \in I} K_i$.
  Then the induced map $u : K \to \Sing^Q Y$
  exhibits $Y$ as the $Q$-realization of $K$.
  This $Q$-realization is preserved by any functor $F : C \to C'$
  which preserves each $Q$-realization $u_i$
  as well as the colimit $Y = \colim_{i \in I} Y_i$.
\end{lemma}

\begin{nul}
  \label{realization-empty}
  Taking $I = \emptyset$,
  we see that the $Q$-realization of $\emptyset \in \sSet$
  is the initial object of $C$ (assuming $C$ has one),
  and this $Q$-realization is preserved by any functor which preserves the initial object.
\end{nul}

\begin{proof}[Proof of \cref{simplicial-complex-realization}]
  Fix $n \ge 0$ and let $\mathcal{T}$ denote
  the poset of simplicial subsets of $\Delta^n$,
  ordered by inclusion.
  The assignment $K \mapsto |K|_\DD$ of \cref{complex-realization-construction}
  defines a functor from $\mathcal{T}$ to the poset of closed subsets of $|\Delta^n|_\DD$.
  For each morphism $K \subset K'$ of $\mathcal{T}$, the square
  \[
    \begin{tikzcd}
      K \ar[d] \ar[r, "u_K"] & \Sing^\DD |K|_\DD \ar[d] \\
      K' \ar[r, "u_{K'}"] & \Sing^\DD |K'|_\DD
    \end{tikzcd}
    \tag{$*$}
  \]
  commutes, because both $u_K$ and $u_{K'}$ are restrictions of $u_{\Delta^n}$.
  We will prove the following statement for each $K \in \mathcal{T}$
  by induction on the number of faces of $\Delta^n$ contained in $K$.
  \begin{itemize}
    \item[($\dagger$)]
    $u_K$ exhibits $|K|_\DD$ as the $\PP$-valued geometric realization of~$K$,
    and any functor $F : \PP \to \C$ satisfying conditions (a) and (b$'$)
    preserves this geometric realization.
  \end{itemize}
  In particular, we can take $F$ to be the inclusion of $\PP$ in $\DD$
  to conclude that $|K|_\DD$ is also the $\DD$-valued geometric realization of~$K$.
  The last part of the lemma follows from the commutativity of ($*$).

  The base case is \cref{realization-empty}.
  Suppose $K$ is nonempty,
  and write $K = \sigma \cup K'$
  with $\sigma$ a maximal face of $\Delta^n$ contained in $K$
  and $K' \subsetneq K$.
  Then also $|K|_\DD = |\sigma|_\DD \cup |K'|_\DD$,
  and so the squares below are pushouts.
  \[
    \begin{tikzcd}
      \sigma \cap K' \ar[r] \ar[d] & K' \ar[d] &
      {|\sigma \cap K'|_\DD} \ar[r] \ar[d] & {|K'|_\DD} \ar[d] \\
      \sigma \ar[r] & K &
      {|\sigma|_\DD} \ar[r] & {|K|_\DD}
    \end{tikzcd}
  \]
  We know ($\dagger$) holds for $\sigma \cap K'$ and for $K'$
  because they have fewer simplices than $K$,
  and also holds for $\sigma$ by \cref{simplex-realization}.
  The pushout in $\PP$ is a closed-by-closed pushout,
  and so we can apply \cref{realization-colimit}.
\end{proof}

\begin{proof}[Proof of \cref{finite-realizations}]
  Let $K$ be a finite simplicial set.
  Recall that the \emph{$n$-skeleton} $\sk^n K$ of $K$
  is the simplicial subset of $K$
  generated by the nondegenerate $k$-simplices of~$K$ for $k \le n$,
  with $\sk^{-1} K = \emptyset$.
  We will prove by induction on $n \ge -1$ that
  $\sk^n K$ admits a $\PP$-valued realization
  preserved by any functor satisfying conditions (a) and (b).
  Then the same holds for $K$
  as we have $K = \sk^n K$ for sufficiently large $n$
  (since $K$ has only finitely many nondegenerate simplices).

  The base case $n = -1$ is \cref{realization-empty}.
  For $n \ge 0$, there is a pushout square in $\sSet$
  \[
    \begin{tikzcd}
      \coprod \partial \Delta^n \ar[r] \ar[d] & \sk^{n-1} K \ar[d] \\
      \coprod \Delta^n \ar[r] & \sk^n K
    \end{tikzcd}
  \]
  where the coproducts on the left
  range over the nondegenerate $n$-simplices of~$K$,
  of which there are only finitely many.
  Construct $|{\sk^n K}|_\PP$ as the pushout
  \[
    \begin{tikzcd}
      {|{\coprod \partial \Delta^n}|_\PP} \ar[r] \ar[d] &
      {|{\sk^{n-1} K}|_\PP} \ar[d] \\
      {|{\coprod \Delta^n}|_\PP \ar[r]} &
      {|{\sk^n K}|_\PP}
    \end{tikzcd}
  \]
  in $\PP$.
  We have already constructed the other three $\PP$-valued realizations
  in this diagram:
  $|{\coprod \partial \Delta^n}|_\PP$ and $|{\coprod \Delta^n}|_\PP$
  because they are simplicial complexes,
  and $|{\sk^{n-1} K}|_\PP$ by the inductive hypothesis.
  Furthermore, $\coprod \partial \Delta^n$ is a subcomplex of $\coprod \Delta^n$,
  so the vertical map on the left is a closed inclusion
  and so this pushout in $\PP$ exists.
  Thus, the pushout of this diagram really is
  the $\PP$-valued realization of $\sk^n K$
  by \cref{realization-colimit}.
  Moreover, suppose $F$ is a functor satisfying conditions (a) and (b).
  Then $F$ also satisfies condition (b$'$) of \cref{simplicial-complex-realization},
  so $F$ preserves the realizations
  $|{\coprod \partial \Delta^n}|_\PP$ and $|{\coprod \Delta^n}|_\PP$,
  and it preserves the realization $|{\sk^{n-1} K}|_\PP$ by the inductive hypothesis.
  Since $F$ preserves closed-by-proper pushouts,
  it also preserves the realization $|{\sk^n K}|_\PP$.
  This completes the inductive argument.

  Similarly,
  suppose $j : K \to K'$ is a monomorphism between finite simplicial sets.
  Then $j$ can be expressed as a composition of finitely many
  pushouts of maps of the form $\partial \Delta^n \to \Delta^n$,
  so it suffices to consider a map which is a single pushout
  \[
    \begin{tikzcd}
      \partial \Delta^n \ar[r] \ar[d] & K \ar[d] \\
      \Delta^n \ar[r] & K'
    \end{tikzcd}
  \]
  but then we can construct $|K'|_\PP$ as the closed-by-proper pushout
  \[
    \begin{tikzcd}
      {|\partial \Delta^n|_\PP} \ar[r] \ar[d] & {|K|_\PP} \ar[d] \\
      {|\Delta^n|_\PP} \ar[r] & {|K'|_\PP}
    \end{tikzcd}
  \]
  and so $|K|_\PP \to |K'|_\PP$,
  a pushout of the closed embedding $|\partial \Delta^n|_\PP \to |\Delta^n|_\PP$,
  is also a closed embedding.
\end{proof}

\begin{nul}
  Combining \cref{realization-functoriality} and \cref{finite-realizations},
  we deduce the existence of a geometric realization functor $|{-}|_\PP : \sSetfin \to \PP$
  which is left adjoint to $\Sing^\PP : \PP \to \sSet$
  relative to the inclusion $\sSetfin \subset \sSet$.
  Since $|K|_\PP$ is also a $\DD$-valued geometric realization of $K$,
  we will generally omit the subscript
  indicating the target category of geometric realization
  when it is $\PP$ or $\DD$.
\end{nul}

\begin{example}
  \label{non-complex-realization}
  Let $K$ be the simplicial set $\Delta^2/\Delta^{\{0,1\}}$.
  Then $|K|$ can be constructed as the following pushout.
  \[
    \begin{tikzcd}
      {|\Delta^{\{0,1\}}|} \ar[r] \ar[d] & * \ar[d] \\
      {|\Delta^2|} \ar[r] & {|K|}
    \end{tikzcd}
  \]
  Up to definable homeomorphism,
  this square agrees with the one of \cref{not-closed-by-closed-pushout}.
  The existence of this pushout relies on
  the definability of multiplication in the o-minimal structure.

  We will later see an example of a functor
  satisfying conditions ($a$) and ($b'$) of \cref{simplicial-complex-realization}
  but not condition ($b$) of \cref{finite-realizations},
  and which does not preserve the geometric realization of $K$.
\end{example}

\begin{proposition}
  \label{horn-inclusion-retract}
  The geometric realization of a horn inclusion $\Lambda^n_i \to \Delta^n$
  ($n \ge 1$, $0 \le i \le n$) admits a retraction.
\end{proposition}

\begin{proof}
  Let $w$ denote the barycenter of the face opposite
  vertex $e^n_i$ of $|\Delta^n|$.
  By \cref{simplicial-complex-realization},
  the $\PP$-realization of the inclusion $\Lambda^n_i \subset \Delta^n$
  may be identified with the inclusion of $|\Lambda^n_i|$ in $|\Delta^n|$,
  where $|\Lambda^n_i|$ is the union of the convex hulls of
  all subsets of $\{e^n_0, \ldots, e^n_n\}$ of size $n$ containing $e^n_i$.
  Then $|\Delta^n|$ retracts onto~$|\Lambda^n_i|$
  by projection in the direction from $w$ to $e^n_i$;
  this is a continuous piecewise linear function with rational coefficients
  and therefore definable.
\end{proof}

\begin{corollary}
  For any object $X$ of $\DD$, $\Sing^\DD X$ is a Kan complex.
\end{corollary}

\begin{proof}
  Using the relative adjunction between $|{-}|$ and $\Sing^\DD$,
  a lifting problem
  \[
    \begin{tikzcd}
      \Lambda^n_i \ar[r] \ar[d] & \Sing^\DD X \\
      \Delta^n \ar[ru, dashed]
    \end{tikzcd}
  \]
  can be converted into the lifting problem
  \[
    \begin{tikzcd}
      {|\Lambda^n_i|} \ar[r] \ar[d] & X \\
      {|\Delta^n|} \ar[ru, dashed]
    \end{tikzcd}
  \]
  but any such lifting problem has a solution
  obtained using a retraction $|\Delta^n| \to |\Lambda^n_i|$
  provided by \cref{horn-inclusion-retract}.
\end{proof}

\begin{corollary}
  Let $\Phi : (R, \sS) \to (R', \sS')$ be a morphism of o-minimal structures.
  Then the diagram
  \[
    \begin{tikzcd}
      \sSetfin \ar[r, "|{-}|"] \ar[rd, "|{-}|"'] & \PP(R, \sS) \ar[d, "\Phi"] \\
      & \PP(R', \sS')
    \end{tikzcd}
  \]
  commutes (up to a canonical natural isomorphism).
\end{corollary}

\begin{proof}
  The definition of the functor $|\Delta^\bullet| : \Delta \to \PP$
  is uniform in the o-minimal structure,
  so $\Phi \circ |\Delta^\bullet|_{\PP(R, \sS)}$
  is canonically isomorphic to $|\Delta^\bullet|_{\PP(R', \sS')}$;
  and $\Phi$ preserves geometric realizations
  because it preserves the initial object and closed-by-proper pushouts
  (by \cref{pretopos-morphism}).
\end{proof}

\begin{corollary}
  Suppose $R = \RR$.
  Then for each finite simplicial set $K$,
  the underlying space $|K|_\tp$ of $|K|$
  is canonically homeomorphic to
  the usual geometric realization $|K|_\Top$ of $K$.
\end{corollary}

\begin{proof}
  The functor ${-}_\tp : \PP \to \Top$
  satisfies the hypotheses of \cref{finite-realizations}
  by \cref{quotients-top},
  and the underlying space of $|\Delta^n|_\PP$
  is $|\Delta^n|_\Top$ by construction.
\end{proof}

\begin{proposition}
  \label{realization-pushouts}
  $|{-}| : \sSetfin \to \PP$ preserves the initial object
  and pushouts of monomorphisms.
\end{proposition}

\begin{proof}
  This follows from \cref{realization-colimit}
  once we check that $|{-}|$ takes these diagrams
  (the empty diagram, and spans with one leg a monomorphism)
  to diagrams which admit a colimit in $\PP$.
  For the former this is obvious
  and for the latter this follows from
  \cref{realization-of-monomorphism} and \cref{closed-by-proper-pushout}.
\end{proof}

\begin{nul}
  \label{finite-realization-finite-limits}
  Geometric realization also preserves finite limits.
  This will be easy to prove later (\cref{geom-real-lex}),
  once we have more tools.
\end{nul}

\begin{nul*}
  Finally, we formulate a version of the triangulation theorem \cite[Theorem~8.2.9]{vdD}
  using the language of geometric realization.
  It provides a converse to the last part of
  \cref{simplicial-complex-realization}.
\end{nul*}

\begin{proposition}
  \label{triangulation}
  Let $f : X \to Y$ be a closed embedding in $\PP$.
  Then $f$ is isomorphic to a map of the form $|K'| \to |K|$
  for $K' \subset K \subset \Delta^n$ an inclusion of simplicial complexes.
\end{proposition}

\begin{proof}
  We regard $X$ as a closed subset of $Y \subset R^N$.
  By \cite[Theorem~8.2.9]{vdD},
  there is a geometric simplicial complex $L$
  and a homeomorphism between $Y$ and the underlying set of $L$
  which takes $X$ to a subcomplex $L'$ of $L$.
  Note that in \cite{vdD}, a simplicial complex is allowed to be ``open'':
  a face of a simplex in the complex might not belong to the complex.
  However, since $Y$ is a closed and bounded definable set,
  the simplicial complex $L$ must be ``closed''
  by the comment after \cite[Definition~8.2.3]{vdD}.
  Likewise, the image of $X$ must be the underlying set
  of a closed simplicial subcomplex $L'$ of $L$.

  Now, choose any ordering of the vertices of $L$
  and let $K$ be the corresponding simplicial complex
  as a simplicial subset of $\Delta^n$ (for a suitably chosen $n$),
  and $K' \subset K$ the subcomplex corresponding to $L' \subset L$.
  By \cref{simplicial-complex-realization},
  $|K|$ is the underlying set of a geometric simplicial complex
  (a subcomplex of $|\Delta^n|$)
  of the same combinatorial type as $L$.
  Then by \cite[\S8.1.6]{vdD} there is an induced definable homeomorphism
  between $|K|$ and the underlying set of $L$,
  and hence between $|K|$ and $Y$.
  This homeomorphism takes $|K'|$,
  which corresponds to the subcomplex $L'$,
  to $X$.
\end{proof}

\begin{remark}
  In particular, the functor $|{-}| : \sSetfin \to \PP$
  is essentially surjective.
  Consequently,
  for any morphism $\Phi : (R, \sS) \to (R', \sS')$ of o-minimal structures,
  the induced functor $\Phi : \PP(R, \sS) \to \PP(R', \sS')$
  is also essentially surjective.
  In a sense, then, enlarging the o-minimal structure
  does not add any new objects;
  it only adds more morphisms between the existing objects.
\end{remark}

\begin{warning}
  In \cref{triangulation},
  suppose $Y$ originated as
  the geometric realization of a simplicial complex $L$,
  and $X$ is a general definable closed subset of $Y$.
  A priori, the triangulation $K$ of $Y$ produced by \cref{triangulation}
  might have an underlying simplicial complex
  with no particular combinatorial relation to $L$.
  For example, if $L$ consists of a single simplex,
  it is not obvious that $K$ is (weakly) contractible as a simplicial set.
  We know only that the realizations of $K$ and $L$
  are related by a definable homeomorphism in the given o-minimal structure.

  In fact, we will see in the next section
  that the homotopy type of a simplicial complex $K$
  is indeed determined by its geometric realization (in any o-minimal structure).
  However, this requires additional inputs beyond
  the triangulation theorem of \cite{vdD98}.
\end{warning}

%% file: tex/homotopy.tex
\section{\texorpdfstring
{The homotopy theory of $\DD$ and $\PP$}
{The homotopy theory of D and P}}
\label{homotopy-d-p}

\begin{definition}
  Two parallel maps $f_0$, $f_1 : A \to X$ of $\DD$ are \emph{homotopic}
  if there exists a map $H : A \times [0, 1] \to X$
  such that $H(a, 0) = f_0(a)$ and $H(a, 1) = f_1(a)$ for all $a \in A$.
  In this case we write $f_0 \simeq f_1$.
\end{definition}

\begin{nul}
  \label{ho-dd-pp}
  Homotopy defines an equivalence relation on $\Hom_\DD(A, X)$
  for each pair of objects $A$ and $X$ of $\DD$,
  and these equivalence relations are compatible with composition.
  Thus, we may define the \emph{classical homotopy category}
  $\Ho \DD$ by $(\Ho \DD)(A, X) = \DD(A, X)/{\simeq}$.
  Similarly, we define $\Ho \PP$ by $(\Ho \PP)(A, X) = \PP(A, X)/{\simeq}$;
  $\Ho \PP$ is a full subcategory of $\Ho \DD$.
\end{nul}

\begin{nul}
  \label{ho-p-eq-ho-d}
  Every object of $\DD$ is homotopy equivalent to an object of $\PP$
  and so the inclusion $\PP \to \DD$ induces
  an equivalence $\Ho \PP \xrightarrow{\sim} \Ho \DD$.
  In fact, for any object~$X$ of~$\DD$
  there is a closed subspace $A \subset X$ which is a polytope
  such that the inclusion $A \to X$ is a homotopy equivalence.
  To construct such an $A$,
  first triangulate $X$ by an \emph{open simplicial complex}
  \cite[Theorem 8.2.9]{vdD98}:
  roughly speaking, a finite complex $K$ built from open simplices
  in which a face of a simplex of $K$ is not required to belong to $K$.
  Subdivide $K$ once to produce a new open complex $K'$,
  and then let $A$ be the \emph{core} of $K'$:
  the union of those open simplices of $K'$
  all of whose faces also belong to $K'$.
  Then there is a strong deformation retract of $X$ onto
  the subspace $A$, which is a finite (closed) simplicial complex.
  We refer the reader to \cite[III.1]{LSA} for more details.

  Therefore, in the remainder of this section
  we discuss only the homotopy theory of $\PP$.
\end{nul}

\begin{nul}
  The main result of this section is that
  the functor
  \[
    \sSetfin \xrightarrow{|{-}|} \PP
  \]
  induces an equivalence of homotopy categories
  \[
    \Ho \sSetfin \xrightarrow{\sim} \Ho \PP
  \]
  where $\Ho \sSetfin$ denotes
  the category obtained by inverting
  the weak equivalences between finite simplicial sets.
  By \cref{ho-ssetfin-ff}, $\Ho \sSetfin$ is also
  the full subcategory of the usual homotopy category $\Ho \sSet$
  whose objects are the finite simplicial sets.

  We rely on previous work of several authors \cite{BO10,DK85}
  which shows that $\Ho \PP$
  is equivalent to the homotopy category of finite CW complexes.
  Our contribution in this section is to reformulate this result
  in a way more suitable for extension to the setting of model categories.
  In particular, the model categories we introduce in \cref{modelcat}
  always come with a direct left Quillen functor from $\sSet$,
  but in general are related to $\Top$ only by a zigzag of Quillen adjunctions.
\end{nul}

\begin{nul*}
  In this section we use the framework of \emph{cofibration categories},
  a model for homotopy theories which admit finite (homotopy) colimits.
  \cref{cofibration-categories} contains background material
  on cofibration categories.

  Our first goal is to equip $\PP$ with the structure of a cofibration category.
  This is easy to do using the following key fact.
\end{nul*}

\begin{proposition}
  \label{pp-hep}
  Every monomorphism of $\PP$ satisfies the homotopy extension property.
  That is, if $j : A \to B$ is a closed embedding
  and $f : B \to X$ and $h : A \times [0, 1] \to X$ are maps
  with $h(a, 0) = f(j(a))$ for $a \in A$,
  then there exists $H : B \times [0, 1] \to X$
  with $H(j(a), t) = h(a, t)$ and $H(b, 0) = f(b)$.
\end{proposition}

\begin{proof}
  This is proven for
  closed embeddings between semialgebraic sets in \cite[Theorem 5.1]{DK84},
  even when $A$ and $B$ are not necessarily polytopes.
  For polytopes, we sketch a direct proof as follows.
  Any closed embedding in $\PP$ can be triangulated,
  so it can be expressed as a composition of pushouts
  of cell inclusions $j_n : |\partial \Delta^n| \to |\Delta^n|$.
  The class of maps with the homotopy extension property
  with respect to a fixed cylinder functor (here $- \times [0, 1]$)
  is closed under composition
  and under pushouts preserved by the cylinder functor
  \cite[I.2.6 and~I.2.7]{KampsPorter}.
  So it suffices to show that
  each $j_n$ satisfies the homotopy extension property,
  and this can be done by writing down an explicit formula
  for a retraction
  $|\Delta^n| \times [0, 1] \to
  |\partial \Delta^n| \times [0, 1] \cup |\Delta^n| \times \{0\}$.
\end{proof}

\begin{proposition}
  $\PP$ has the structure of a cofibration category
  in which a morphism is a cofibration if it is a monomorphism
  and a weak equivalence if it is a homotopy equivalence.
  In this cofibration category, every object is cofibrant and fibrant.
\end{proposition}

\begin{proof}
  We check the axioms of an $I$-category in \cref{cc-of-icat},
  using the cylinder functor $I$ given by $IA = A \times [0, 1]$
  with $p(a, t) = a$ and $i_\varepsilon(a) = (a, \varepsilon)$ for $\varepsilon \in \{0, 1\}$.
  (I1) is clear.
  (I2) follows from \cref{closed-by-proper-pushout};
  $I$ preserves pushouts of monomorphisms
  because such pushouts are pullback-stable.
  (I3) follows from \cref{pp-hep}.
  (I4) follows from identifying $B \amalg_{A,0} IA \amalg_{A,1} B$
  with the closed subset $\{\,(b, t) \mid b \in A \vee t = 0 \vee t = 1\,\}$
  of $B \times [0, 1]$.
  For (I5), take $T : IIA \to IIA$ defined by $T((a, t), t') = ((a, t'), t)$.
\end{proof}

\begin{nul}
  By \cref{cc-ho-cof-fib},
  since every object of $\PP$ is both cofibrant and fibrant,
  we can identify $\Ho \PP$ as defined in \cref{ho-dd-pp}
  with the category obtained by inverting the weak equivalences of $\PP$,
  as in \cref{ho-cat}.
  For brevity, we denote the Hom-sets $\Ho C(A, B)$ of the category $\Ho C$
  by $[A, B]_C$
  when they can be computed directly
  as homotopy classes of maps from $A$ to $B$,
  such as when $C = \PP$ or when $A$ and $B$ are cell complexes in $\Top$.
\end{nul}

\begin{nul}
  Recall that the morphisms of cofibration categories
  are the exact functors (\cref{cofibration-category}):
  the functors which preserve cofibrations, acyclic cofibrations,
  the initial object, and pushouts of cofibrations.
  By \cref{realization-of-monomorphism,realization-pushouts},
  the geometric realization functor $|{-}| : \sSetfin \to \PP$
  preserves all of these
  except (possibly) acyclic cofibrations.
  Preservation of acyclic cofibrations can be checked directly using the fact that
  every acyclic cofibration of $\sSetfin$
  can be written as a retract of a finite composition of pushouts
  of horn inclusions $\Lambda^n_i \to \Delta^n$
  \cite[Proposition 5.13]{BaSc2}.
  We will also obtain this
  as a consequence of the proof below.
\end{nul}

\begin{proposition}
  \label{ho-p}
  The functor $|{-}| : \sSetfin \to \PP$ is exact
  and induces an equivalence of homotopy categories.
\end{proposition}

\begin{proof}
  We denote the o-minimal structure in question by $(R, \sS)$
  and recall the following facts which, taken together,
  relate the homotopy theory of $\PP$ to that of $\Top$.
  \begin{enumerate}
  \item
    By \cite[Theorem 3.1]{BO10},
    the (non-full) inclusion $J : \PP(R_\sa) \to \PP(R, \sS)$
    induces isomorphisms
    \[
      [A, B]_{\PP(R_\sa)} \xrightarrow{\sim} [JA, JB]_{\PP(R,\sS)}.
    \]
  \item
    By \cite[Theorem III.3.1]{DK85},
    for any extension of real closed fields $R' \subset R''$,
    the associated functor $\Phi : \PP(R'_\sa) \to \PP(R''_\sa)$
    of \cref{associated-functor}
    induces isomorphisms
    \[
      [A, B]_{\PP(R'_\sa)} \xrightarrow{\sim} [\Phi A, \Phi B]_{\PP(R''_\sa)}.
    \]
  \item
    By \cite[Theorem III.5.1]{DK85},
    the functor $U : \PP(\RR_\sa) \to \Top$
    sending a polytope to its underlying topological space
    induces isomorphisms
    \[
      [A, B]_{\PP(\RR_\sa)} \xrightarrow{\sim} [UA, UB]_\Top.
    \]
  \end{enumerate}
  Now, consider the following diagram,
  in which each functor out of $\sSetfin$ is given by geometric realization,
  and the remaining ones are the comparison functors described above.
  This diagram commutes up to natural isomorphism
  because the functors $J$, $\Phi$, $U$ preserve
  initial objects and closed-by-proper pushouts.
  \[
    \begin{tikzcd}
      \sSetfin \ar[rd] \ar[rrd] \ar[rrrd, bend left=10]
      \ar[rdd] \ar[rddd, bend right=10] \\
      & \PP(\Qalgre_\sa) \ar[r, "\mathit{(ii)}"'] \ar[d, "\mathit{(ii)}"]
      & \PP(\RR_\sa) \ar[r, "\mathit{(iii)}"'] & \Top \\
      & \PP(R_\sa) \ar[d, "\mathit{(i)}"] \\
      & \PP(R, \sS)
    \end{tikzcd}
    \tag{$*$}
  \]
  Let $K$ and $L$ denote objects of $\sSetfin$.
  Then there is a zigzag of isomorphisms
  \begin{align*}
    (\Ho \sSetfin)(K, L)
    & \xrightarrow{\sim} (\Ho \sSet)(K, L)
    \xrightarrow{\sim} [|K|, |L|]_\Top \\
    & \xleftarrow{\sim} [|K|, |L|]_{\PP(\RR_\sa)}
    \xleftarrow{\sim} [|K|, |L|]_{\PP(\Qalgre_\sa)} \\
    & \xrightarrow{\sim} [|K|, |L|]_{\PP(R_\sa)}
    \xrightarrow{\sim} [|K|, |L|]_{\PP(R,\sS)}
  \end{align*}
  natural in $K$ and $L$.
  In the first step, we used \cref{ho-ssetfin-ff}.
  The second isomorphism is classical,
  and the rest are given by facts (\textit{i}), (\textit{ii}), (\textit{iii})
  as indicated in the diagram.

  A priori, we do not know that this isomorphism
  is induced by a functor $\Ho \sSetfin \to \Ho \PP(R, \sS)$.
  However, if $i : L \to L'$ is a weak equivalence in $\sSetfin$,
  then from the naturality square
  \[
    \begin{tikzcd}
      \Ho \sSetfin(K, L) \ar[r, "\sim"] \ar[d, "\sim"']
      & \Ho \sSetfin(K, L') \ar[d, "\sim"] \\
      {[|K|, |L|]_{\PP(R, \sS)}} \ar[r] & {[|K|, |L'|]_{\PP(R, \sS)}}
    \end{tikzcd}
  \]
  we conclude that
  the induced map $[|K|, |L|]_{\PP(R, \sS)} \to [|K|, |L'|]_{\PP(R, \sS)}$
  is an isomorphism.
  Applying this to $K = L$ and $K = L'$,
  it follows that $|i| : |L| \to |L'|$ is a homotopy equivalence.
  Therefore, $|{-}| : \sSetfin \to \PP(R, \sS)$ preserves weak equivalences,
  so we have completed the verification that it is an exact functor.

  The remaining morphisms out of $\sSetfin$ in ($*$) also preserve weak equivalences,
  since we could specialize $(R, \sS)$ to the other o-minimal structures in the diagram.
  Thus, they descend uniquely to functors on the level of homotopy categories
  and so we obtain an induced diagram
  \[
    \begin{tikzcd}
      \Ho \sSetfin \ar[rd] \ar[rrd] \ar[rrrd, bend left=10]
      \ar[rdd] \ar[rddd, bend right=10] \\
      & \Ho \PP(\Qalgre_\sa) \ar[r] \ar[d] & \Ho \PP(\RR_\sa) \ar[r] & \Ho \Top \\
      & \Ho \PP(R_\sa) \ar[d] \\
      & \Ho \PP(R, \sS)
    \end{tikzcd}
  \]
  in which each functor out of $\Ho \sSetfin$ is given by $\Ho |{-}|$.
  The horizontal and vertical functors in this diagram are fully faithful,
  as is $\Ho |{-}| : \Ho \sSetfin \to \Ho \Top$.
  It follows that $\Ho |{-}| : \Ho \sSetfin \to \Ho \PP(R, \sS)$ is fully faithful,
  and it is essentially surjective
  because $|{-}| : \sSetfin \to \PP(R, \sS)$ is already essentially surjective
  (every polytope can be triangulated).
\end{proof}

\begin{nul*}
  For future use, we extract the following consequence.
\end{nul*}

\begin{corollary}
  \label{eta-weq-of-finite}
  For each $K \in \sSetfin$,
  the unit map $\eta_K : K \to \Sing |K|_\PP$
  is a weak homotopy equivalence of simplicial sets.
\end{corollary}

\begin{proof}
  The hypotheses of \cref{cc-unit-weq} are satisfied:
  $|{-}| : \sSetfin \to \PP$ is an exact functor
  inducing an equivalence of homotopy categories
  and every object of $\PP$ is fibrant.
\end{proof}

\begin{remark}
  \label{pretopoi-for-spaces}
  We have shown that the category $\PP$ is a pretopos
  which also has a cofibration category structure
  that models the homotopy theory of finite CW complexes.
  Other categories of the same type include:
  \begin{itemize}
  \item
    $\sSetfin$, the category of finite simplicial sets.
    The cofibration category structure is the one
    obtained by restricting the standard model structure on $\sSet$.
    In the cofibration category $\sSetfin$,
    every object is cofibrant,
    but there are very few fibrant objects.
  \item
    $\CompHaus$, the category of compact Hausdorff topological spaces.
    (We could impose some cardinality restriction
    to make $\CompHaus$ a small category, if desired.)
    The cofibration category structure is again the one
    obtained by restricting the standard model structure on $\Top$.
    The main point is that $\CompHaus \subset \Top$
    is closed under pushouts \emph{of monomorphisms},
    so in particular under pushouts of cofibrations.
    The usual mapping cylinder construction
    produces the required factorizations.
    In the cofibration category $\CompHaus$,
    every object is fibrant but not every object is cofibrant.
  \end{itemize}
  When $R = \RR$, there are pretopos morphisms
  \[
    \sSetfin \xrightarrow{|{-}|} \PP \xrightarrow{{-}_\tp} \CompHaus
  \]
  which are also exact functors of cofibration categories,
  and factor the topological geometric realization functor.
  The pretopos $\PP$ sits between $\sSetfin$ and $\CompHaus$
  in a ``Goldilocksian'' way:
  it has enough maps between objects
  (in particular, retractions for ``open cylinder'' inclusions)
  to make every object fibrant,
  but also few enough maps between objects
  to allow every monomorphism to be a cofibration.
  Much of what we prove here and in \cite{paper3}
  is a formal consequence of these properties of $\PP$,
  and could be generalized to any pretopos
  (equipped with an interval, to define homotopy)
  with the same properties.
  However, we do not know of any other examples!

  The functor ${-}_\tp : \PP \to \CompHaus$ is faithful but not full.
  This is essential:
  we cannot cut down $\CompHaus$ to a category with the same properties as $\PP$
  merely by removing objects,
  because the equalizer of two continuous maps between ``nice'' spaces
  may be a space which is not ``nice''.
  We must remove morphisms as well.
\end{remark}

%% file: tex/topologies.tex
\section{\texorpdfstring
{Grothendieck topologies on $\DD$ and $\PP$}
{Grothendieck topologies on D and P}}
\label{topologies}

\begin{nul*}
  The categories $\DD$ and $\PP$ model
  the homotopy theory of finite CW complexes.
  However, many basic constructions in homotopy theory
  leave the finite world,
  such as the formation of loop spaces.
  Our aim is to enlarge $\DD$ and $\PP$ to model categories
  which model the full homotopy theory of spaces.
  In doing so, we would like to preserve certain ``well-behaved'' colimits
  while freely adjoining other colimits
  which are badly-behaved or nonexistent in $\DD$ or $\PP$.
  Sheaf theory provides a general mechanism for achieving this.

  In this section we introduce several Grothendieck topologies
  on the categories $\DD$ and $\PP$.
  We will eventually be interested mainly in the \emph{proper topology}
  (\cref{proper-topology}).
\end{nul*}

\begin{definition}
  \label{open-topology}
  The \emph{open topology} on $\DD$ is the Grothendieck topology
  generated by the pretopology in which
  the covering families of an object $X \in \DD$
  are the \emph{finite} covers $\{U_i \to X\}_{i \in I}$
  by definable open subsets $U_i$ of $X$.

  The open topology on $\PP$ is defined similarly,
  using finite families $\{V_i \to X\}_{i \in I}$
  of definable closed subsets $V_i$ of $X$ whose interiors cover~$X$.
\end{definition}

\begin{proposition}
  \label{open-topology-subcanonical}
  The open topologies on $\DD$ and on $\PP$ are subcanonical.
\end{proposition}

\begin{proof}
  This follows from the fact that, for any definable sets $X$ and $Y$
  and any finite family $\{V_i \to X\}_{i \in I}$ of definable subsets of $X$
  whose interiors cover~$X$,
  an arbitrary function from $X$ to $Y$ is continuous and definable
  if and only if its restriction to each $V_i$ is so.
  (The ``glued'' function is definable because
  its graph is the union of the graphs of the original functions.)
\end{proof}

\begin{remark}
  The category $\DD$ with the open topology
  is a definable analogue of the big Zariski site in algebraic geometry.
  However,
  only \emph{definable} open subsets can possibly appear in a covering family in $\DD$,
  by the way $\DD$ was defined.
  Moreover, we need to restrict to \emph{finite} covering families
  to obtain a subcanonical topology,
  because of examples like the one in \cref{no-infinite-colimits}.
  For these reasons, we don't use the name ``Zariski topology''.

  One way to make the analogy to algebraic geometry more exact
  is to replace each object $X \in \DD$
  by the coherent (or spectral) space $\widetilde X$
  which corresponds via Stone duality
  to the distributive lattice of definable open subsets of $X$,
  as described in \cite{EJP06}.
  Then definable open subsets of $X$
  correspond to quasicompact open subsets of $\widetilde X$.
  In particular, $\widetilde X$ itself is always quasicompact
  (regardless of whether $X$ was definably compact),
  so it makes no difference whether we restrict to finite covers of $\widetilde X$.
  For our purposes we can make do with the original definable set $X$,
  as long as we declare only finite covers to be generating covering families.
\end{remark}

\begin{nul}
  \label{covers}
  Recall that
  a family of morphisms $\{f_\alpha : Y_\alpha \to X\}_{\alpha \in A}$
  is said to \emph{cover} $X$ in a Grothendieck topology $\tau$
  if the sieve it generates belongs to $\tau$.
  When $\tau$ is the topology generated by a given Grothendieck pretopology,
  this means that there is a generating covering family of $X$
  each of whose morphisms factors through $f_\alpha$ for some $\alpha \in A$.
  This is also equivalent to the condition that
  $\coprod_{\alpha \in A} \yo Y_\alpha \to \yo X$ is an epimorphism
  in the category of sheaves for the topology $\tau$.
\end{nul}

\begin{proposition}
  \label{open-topology-covers}
  A family of morphisms $\{f_\alpha : Y_\alpha \to X\}_{\alpha \in A}$
  covers $X$ in the open topology (on $\DD$ or on $\PP$)
  if and only if there exists a finite set of
  definable open subsets $U_i \subset X$ covering $X$
  such that, for each $i$,
  some $f_{\alpha_i}$ splits after pullback along the inclusion $U_i \to X$.
\end{proposition}

\begin{proof}
  In $\DD$, the family $\{f_\alpha : Y_\alpha \to X\}_{\alpha \in A}$ covers $X$
  if and only if the sieve it generates contains the sieve generated by
  a finite definable open cover $\{U_i \to X\}_{i \in I}$.
  This means for each $i \in I$ there exists $\alpha_i \in A$ such that
  the inclusion $U_i \to X$ factors through $f_{\alpha_i} : Y_{\alpha_i} \to X$,
  or equivalently, such that $f_\alpha$ splits after pullback along $U_i \to X$.

  In $\PP$, we conclude similarly that
  the family $\{f_\alpha : Y_\alpha \to X\}_{\alpha \in A}$ covers $X$
  if and only if there exists
  a finite family $\{V_i \to X\}_{i \in I}$ of definable closed subsets of $X$
  such that for each $i$,
  some $f_{\alpha_i}$ splits after pullback along $V_i \to X$.
  If so, then setting $U_i$ to be the interior of $V_i$,
  the sets $U_i$ form a definable open cover of $X$
  satisfying the conclusion of the lemma.
  Conversely, suppose $\{U_i \to X\}_{i \in I}$ is
  any finite cover of $X$ by definable open subsets
  such that for each $i$,
  some $f_{\alpha_i}$ splits after pullback along $U_i \to X$.
  By the shrinking lemma (\cref{shrinking-lemma}),
  we may find definable closed subsets $V_i$ of $X$ with $V_i \subset U_i$
  whose interiors still cover $X$;
  then for each $i$,
  the same $f_{\alpha_i}$ still splits after pullback along $V_i \to X$.
\end{proof}

\begin{nul*}
  The open topology does not suit our aims
  because it only allows gluing along open subsets
  and for this reason its category of sheaves fails to resemble $\Top$,
  as we explain next.
\end{nul*}

\begin{nul}
  \label{too-few-open-colimits}
  In the category of sheaves (on $\DD$ or on $\PP$) for the open topology,
  let $J$ be the sheaf defined by the pushout square shown below.

  \[
    \begin{tikzcd}
      \yo \{\frac 12\} \ar[r] \ar[d] & \yo {[0, \frac 12]} \ar[d] \\
      \yo {[\frac 12, 1]} \ar[r] & J
    \end{tikzcd}
  \]
  The inclusion maps $[0, \frac 12] \to [0, 1]$ and $[\frac 12, 1] \to [0, 1]$
  induce a morphism $\varphi : J \to \yo [0, 1]$.
  However, $\varphi$ is not an isomorphism.
  If it was one, the map $\yo [0, \frac 12] \amalg \yo [\frac 12, 1] \to \yo [0, 1]$
  would be an epimorphism,
  so the two maps $[0, \frac 12] \to [0, 1]$ and $[\frac 12, 1] \to [0, 1]$
  would cover $[0, 1]$ in the open topology.
  But there is no finite open cover of $[0, 1]$ over each set of which
  one of these two maps splits.

  In $\Top$, and also in $\PP$ and in $\DD$, the square
  \[
    \begin{tikzcd}
      \{\frac 12\} \ar[r] \ar[d] & {[0, \frac 12]} \ar[d] \\
      {[\frac 12, 1]} \ar[r] & {[0, 1]}
    \end{tikzcd}
  \]
  \emph{is} a pushout,
  and this fact is what allows paths to be composed by the usual formula.
  We would like to preserve this feature in a category of sheaves.
  To do so, we need to use a finer topology,
  one for which the Yoneda embedding preserves pushout squares of this type.
\end{nul}

\begin{definition}
  \label{closed-topology}
  The \emph{closed topology} on $\DD$ or on $\PP$ is the Grothendieck topology
  generated by the pretopology in which
  the covering families of an object $X$
  are the finite covers $\{V_i \to X\}_{i \in I}$ by definable closed subsets $V_i$ of $X$.
  We denote the closed topology by $\taucl$.
\end{definition}

\begin{proposition}
  \label{closed-topology-subcanonical}
  The closed topology (on $\DD$ or on $\PP$) is subcanonical.
\end{proposition}

\begin{proof}
  As \cref{open-topology-subcanonical},
  using the fact that a function is continuous
  if its restriction to each member of a finite cover by closed subsets is.
\end{proof}

\begin{proposition}
  \label{closed-topology-covers}
  A family of morphisms $\{f_\alpha : Y_\alpha \to X\}_{\alpha \in A}$
  covers $X$ in the closed topology (on $\DD$ or on $\PP$)
  if and only if there exist a finite set $I$
  and definable closed subsets $V_i \subset X$ covering $X$
  such that, for each $i$,
  some $f_{\alpha_i}$ splits after pullback along the inclusion $V_i \to X$.
\end{proposition}

\begin{proof}
  This follows from the definition of the closed topology
  together with~\cref{covers}.
\end{proof}

\begin{nul}
  Recall from \cref{closed-by-closed-pushout} that
  if the definable set $X$ is covered by
  two definable closed subsets $V_1$ and $V_2$,
  then the square
  \[
    \begin{tikzcd}
      V_1 \cap V_2 \ar[r] \ar[d] & V_1 \ar[d] \\
      V_2 \ar[r] & X
    \end{tikzcd}
  \]
  is a pushout (in $\DD$, or in $\PP$ if $X$ belongs to $\PP$),
  as well as a pullback;
  squares of this form are the closed-by-closed pushouts.
\end{nul}

\begin{proposition}
  \label{closed-sheaves-squares}
  A presheaf $F$ (on $\DD$ or on $\PP$) is a sheaf for the closed topology
  if and only if it satisfies both of the following conditions:
  \begin{itemize}
  \item[(a)] $F(\emptyset) = *$.
  \item[(b$'$)] $F$ sends closed-by-closed pushouts to pullback squares.
  \end{itemize}
  In particular, the Yoneda embedding (valued in sheaves for the closed topology)
  preserves the initial object
  and sends closed-by-closed pushouts to pushouts.
\end{proposition}

\begin{proof}
  By \cref{closed-by-closed-pushout},
  condition (b$'$) is equivalent to
  the sheaf condition for all covers of an object $X$
  by two definable closed subsets $V_1$ and $V_2$.
  Then the first claim follows by an inductive argument.

  The second claim follows from the ``only if'' direction,
  because the Yoneda embedding $\yo$ preserves a given colimit
  if and only if $\Hom(\yo -, F)$ sends it to a limit for every sheaf $F$;
  and by the Yoneda lemma,
  this functor is isomorphic to $F$ itself.
\end{proof}

\begin{proposition}
  \label{closed-finer-than-open}
  The closed topology (on $\DD$ or on $\PP$)
  is strictly finer than the open topology.
\end{proposition}

\begin{proof}
  On $\PP$, any generating cover for the open topology
  is also a generating cover for the closed topology.
  On $\DD$, we can use the shrinking lemma (\cref{shrinking-lemma})
  to refine any finite cover of a definable set $X$
  by definable open subsets $\{U_i\}_{i \in I}$
  to a cover $\{V_i\}_{i \in I}$ by definable closed subsets
  with $V_i \subset U_i$ for each $i$.
  The closed topology is not equal to the open topology
  because they have different sheaves
  by \cref{too-few-open-colimits,closed-sheaves-squares}.
\end{proof}

\begin{notecons}
  \label{constructive}
  In classical mathematics, the real numbers $\RR$ satisfy both of
  the following properties.
  \begin{enumerate}
  \item
    If $a < a' < b < b'$ are real numbers,
    then $(a, b) \cup (a', b') = (a, b')$.
  \item
    If $a < b < c$ are real numbers,
    then $[a, b] \cup [b, c] = [a, c]$.
  \end{enumerate}
  Constructively, property (\textit{i}) is valid
  but property (\textit{ii}) is not.
  The problem is that the real numbers need not satisfy
  the totality axiom $\vdash x \le y \vee y \le x$.
  That means that the second diagram of \cref{too-few-open-colimits}
  is not necessarily a pushout in $\Top$ (or even in $\Set$).
  As a result,
  a constructive development of the homotopy theory of topological spaces along traditional lines
  would run into the same sort of problems that we would have with
  sheaves on $\DD$ for the open topology,
  such as the inability to concatenate paths by the usual formula.

  One workaround for this issue is to work with
  paths defined on $[0, 1]$ which are constant on
  open neighborhoods of the endpoints.
  Such paths \emph{can} be concatenated constructively,
  using property (\textit{i}).
  More generally, any kind of gluing must be replaced by
  gluing along open subsets.
  This is a departure from classical algebraic topology,
  in which one builds cell complexes by attaching cells along their boundaries.

  Here, we propose an alternative approach
  to building a constructive homotopy theory of spaces
  which better resembles the classical one.
  Namely, we dispense with the real numbers entirely
  and work instead in an o-minimal structure
  of a strongly constructive nature:
  \begin{itemize}
  \item
    $R$ satisfies the trichotomy law
    $\vdash x < y \vee x = y \vee x > y$.
  \item
    There is an algorithm\footnote{
      Internally to the constructive logic,
      an ``algorithm'' is just a function.
      Externally, in order to give an example of such a function,
      we must provide an algorithm to compute it.
    }
    which accepts any formula in one free variable
    in the language of an ordered ring with parameters from $R$,
    and outputs a description of the subset of $R$ it defines
    as the union of a list of singletons and (possibly infinite) open intervals.
  \end{itemize}
  In particular, we can compute the sum $a + b$ for given $a$, $b \in R$
  by applying this algorithm to the formula $x = a + b$
  in the free variable $x$,
  and likewise for the other field operations.
  For example, the initial o-minimal structure $\Qalgre_\sa$
  is of this type
  because the real algebraic numbers form
  a computable field satisfying the trichotomy law
  and the theory of real closed fields
  has an algorithm for quantifier elimination.

  Although we have not checked all details,
  it appears that the theory of o-minimal geometry
  becomes constructively valid under these hypotheses.
  (Many proofs in \cite{vdD98} do invoke
  case analysis on the truth of certain propositions,
  but these propositions are (or are equivalent to)
  closed formulas in the language of the o-minimal structure
  and so are decidable by the second hypothesis.)
  Since $R$ satisfies trichotomy,
  the second square of \cref{too-few-open-colimits}
  will be a pushout in $\DD$,
  and we can preserve this pushout
  by taking sheaves for the closed topology.

  In this work we do not restrict ourselves to constructive methods,
  but we also know of no particular obstruction
  to building a constructive version of the theory developed here.
\end{notecons}

\begin{nul}
  The category $\Shv(\DD, \taucl)$ of sheaves on $\DD$ for the closed topology
  is of course cocomplete,
  so there is a geometric realization functor
  $|{-}| : \sSet \to \Shv(\DD, \taucl)$
  which preserves colimits
  and fits into the square below.
  \[
    \begin{tikzcd}
      \Delta \ar[r, "|\Delta^\bullet|"] \ar[d, "\yo"'] & \DD \ar[d, "\yo"] \\
      \sSet \ar[r, "|{-}|"] & \Shv(\DD, \taucl)
    \end{tikzcd}
  \]
  However, the square
  \[
    \begin{tikzcd}
      \sSetfin \ar[r, "|{-}|"] \ar[d] & \DD \ar[d, "\yo"] \\
      \sSet \ar[r, "|{-}|"] & \Shv(\DD, \taucl)
    \end{tikzcd}
  \]
  does not commute, as we show in the example below.
  It does commute when restricted to simplicial \emph{complexes}
  by \cref{simplicial-complex-realization}
  because $\yo : \DD \to \Shv(\DD, \taucl)$
  preserves the initial object and closed-by-closed pushouts
  (\cref{closed-sheaves-squares}).
  It would be preferable to arrange that
  the Yoneda embedding preserves
  the geometric realizations of all finite simplicial sets,
  and we can achieve this using the proper topology
  which we will introduce below.
\end{nul}

\begin{example}
  \label{quotient-not-closed-cover}
  Let $K$ be the simplicial set $\Delta^2 / \Delta^{\{0,1\}}$
  of \cref{non-complex-realization}.
  We claim that the $\Shv(\DD, \taucl)$-valued geometric realization of $K$
  does not agree with $\yo |K|_\DD$.
  In other words, the diagram
  \[
    \begin{tikzcd}
      \yo |\Delta^{\{0,1\}}| \ar[r] \ar[d] & \yo |\Delta^0| \ar[d] \\
      \yo |\Delta^2| \ar[r] & \yo |K|
    \end{tikzcd}
  \]
  is not a pushout in $\Shv(\DD, \taucl)$.
  It suffices to check that the bottom map
  $\yo |\Delta^2| \to \yo |K|$ is not an epimorphism in $\Shv(\DD, \taucl)$,
  or equivalently that
  the quotient map $p : |\Delta^2| \to |K|$
  is not a cover in $\DD$ for the closed topology.
  (Up to definable homeomorphism,
  $p$ is the same as the map $q$ of \cref{not-closed-by-closed-pushout}.)
  By \cref{closed-topology-covers},
  this amounts to the statement that there is no finite closed cover
  $(V_i)_{i \in I}$ of $|K|$
  over each member of which $p$ splits.

  Suppose for the sake of contradiction that
  $(V_i)_{i \in I}$ is a finite closed cover of $|K|$
  and, for each $i \in I$, $g_i : V_i \to |\Delta^2|$
  is a section of $p$ over $V_i$.
  Let $W_i$ be the image of $g_i$ and set $W = \bigcup_{i \in I} W_i$.
  Since $p$ is a homeomorphism
  away from the subset $|\Delta^{\{0, 1\}}| \subset |\Delta^2|$,
  every point outside $|\Delta^{\{0, 1\}}|$ must belong to $W$.
  But $W$ is closed, so $W$ must contain all of $|\Delta^{\{0, 1\}}|$ as well.
  This is a contradiction since
  each $W_i$ can contribute only one point to $W \cap |\Delta^{\{0, 1\}}|$.

  The failure of the Yoneda embedding $\yo$ to preserve
  the geometric realization of $K$
  is due to the fact that $\yo$
  preserves closed-by-closed pushouts
  but not closed-by-proper pushouts.
  The same argument would apply if we replaced $\DD$ by $\PP$ throughout.
\end{example}

\begin{definition}
  \label{proper-topology}
  The \emph{proper topology} on $\DD$ or $\PP$ is the Grothendieck topology
  generated by the pretopology in which
  the covering families of an object $X$
  are the finite families of proper maps $\{Y_i \to X\}_{i \in I}$
  which are jointly surjective.
  We denote the proper topology by $\taupr$.
\end{definition}

\begin{convention}
  Starting in the next section, we will always assume that
  the categories $\PP$ and $\DD$ are equipped with the proper topology
  unless otherwise specified.
\end{convention}

\begin{proposition}
  The proper topology is subcanonical,
  and strictly finer than the closed topology.
\end{proposition}

\begin{proof}
  If $q : Y \to X$ is a proper surjective map
  then it is the coequalizer of its kernel pair
  by \cref{coequalizer-of-proper-surjective}.
  This implies that the proper topology is subcanonical.

  Closed embeddings are proper maps
  so every generating covering family for the closed topology
  is a covering family for the proper topology.
  On the other hand,
  \cref{quotient-not-closed-cover} gave an example of a proper surjection
  which is not a cover for the closed topology.
\end{proof}

\begin{nul}
  In $\PP$, the condition that the maps be proper is redundant,
  so the covering families are just the finite jointly surjective families.
  Thus, the proper topology on $\PP$ is also the coherent topology
  (described in \cref{p-pretopos})
  and so $\PP$ can be recovered as the pretopos of coherent objects
  of $\Shv(\PP, \taupr)$.
  In $\DD$, a general finite jointly surjective family is not a cover.
  For instance, the map $\{0\} \amalg (0, 1] \to [0, 1]$ of $\DD$ is surjective
  but cannot be a cover in any subcanonical topology on $\DD$,
  since a definable function $[0, 1] \to R$
  whose restriction to $\{0\} \amalg (0, 1]$ is continuous
  need not itself be continuous.

  One could define the (not subcanonical) \emph{constructible topology} on $\DD$
  to have all finite jointly surjective families as covers.
  In the category of sheaves on $\DD$ for the constructible topology,
  the maps between the sheafifications of representable objects $\yo X$ and $\yo Y$
  are \emph{all} the definable (but not necessarily continuous) maps from $X$ to $Y$.
  This follows from \cref{continuous-on-partition}.
  So, one might as well instead work with
  the category $\mathrm{Def}$ of definable sets and all definable functions
  from the start.
  We will not make use of the constructible topology here.
\end{nul}

We now come to our main technical result,
a second characterization of the sheaves for the proper topology.

\begin{proposition}
  \label{proper-sheaves-squares}
  A presheaf $F$ (on $\DD$ or on $\PP$) is a sheaf for the proper topology
  if and only if it satisfies both of the following conditions:
  \begin{itemize}
  \item[(a)] $F(\emptyset) = *$.
  \item[(b)] $F$ sends closed-by-proper pushouts to pullback squares.
  \end{itemize}
  In particular, the Yoneda embedding (valued in sheaves for the proper topology)
  preserves the initial object
  and sends closed-by-proper pushouts to pushouts.
\end{proposition}

\begin{remark}
  \Cref{proper-sheaves-squares} is analogous to
  the characterization of sheaves for the ``proper cdh-structure''
  used in algebraic geometry \cite{Voe10_unstable_motivic_homotopy}.
  (This analogy can be made precise: see \cref{closed-by-proper-iff-blowup}.)
  The proof follows the same general outline,
  using the notion of splitting sequence (\cref{splitting-sequence}),
  and could be fit into Voevodsky's framework of
  \emph{cd-structures}~\cite{Voe10_homotopy_theory_decomposed}.
  In order to keep the exposition self-contained,
  we do not explicitly use the notion of cd-structure here.

  In algebraic geometry,
  the covers $p : Y \to X$ of the proper cdh-topology
  have to be ``completely decomposed'' (``cd'')
  in the sense that every point of $X$ has a preimage in $Y$
  with the same residue field.
  In our setting, there is no analogue of this condition
  because any point of any definable set has the same ``residue field'' (namely $R$).
\end{remark}

\begin{nul*}
  We defer the proof of \cref{proper-sheaves-squares} to the end of this section.
\end{nul*}

\begin{remark}
  \label{topologies-agree}
  Let $C$ be any pretopos.
  The following are equivalent for a presheaf $F$ on $C$:
  \begin{itemize}
  \item
    $F(\emptyset) = *$ and $F$ takes pushouts of monomorphisms to pullback squares.
  \item
    $F$ is a sheaf for the topology $\tauad$ generated by
    the empty cover of $\emptyset$
    together with the family $\{B \xrightarrow{f'} Y, X \xrightarrow{j'} Y\}$
    for each pushout
    \[
      \begin{tikzcd}
        A \ar[r, "f"] \ar[d, "j"'] & X \ar[d, "j'"] \\
        B \ar[r, "f'"] & Y
      \end{tikzcd}
    \]
    of a monomorphism $j$ by an arbitrary map $f$.
  \end{itemize}
  This is not a tautology because the map $f'$ need not be a monomorphism.
  It can be shown using the theory of cd-structures,
  or by the proof of in Theorem~4.1 of \cite{GL12}.

  In general, the topology $\tauad$ might be coarser than the coherent topology.
  For example, if $G$ is a nontrivial finite group
  and $C$ is the pretopos of finite $G$-sets,
  then a family of maps is a $\tauad$-covering
  if and only if the maps are jointly surjective on $H$-fixed points
  for every subgroup $H$ of $G$,
  whereas the coverings for the coherent topology
  are just the families that are jointly surjective on the underlying sets.
  For $G = \ZZ/2$, this is precisely the difference between
  the Nisnevich and the \'etale topologies on $\Spec \RR$.
  On the other hand, it is not too hard to see that
  $\tauad$ equals the coherent topology for the pretopos of finite simplicial sets,
  using the inductive construction of a finite simplicial set by attaching simplices,
  which are projective objects.

  The essential content of \cref{proper-sheaves-squares} is that
  $\tauad$ equals the coherent topology for the pretopos $\PP$.
  The point here is that
  this property can be formulated purely in terms of $\PP$ as a category,
  without making reference to its definition.
  One could give an analogous formulation for $\DD$
  by also taking into account the class of proper maps.
\end{remark}

\begin{nul}
  The topologies we have introduced are all functorial
  in the o-minimal structure, in the following sense.
  Let $\Phi : (R, \sS) \to (R', \sS')$ be a morphism of o-minimal structures.
  Then $\Phi$ determines functors
  \[
    \Phi : \PP(R, \sS) \to \PP(R', \sS'), \quad
    \Phi : \DD(R, \sS) \to \DD(R', \sS')
  \]
  and for each topology
  $\tau \in \{\mathrm{open}, \mathrm{closed}, \mathrm{proper}\}$,
  the functors $\Phi$ preserve generating covering families for $\tau$.
  This is because $\Phi$ preserves
  open embeddings, closed embeddings, properness,
  finite coproducts, and surjectivity.
  Consequently, we obtain geometric morphisms
  \[
    \Phi_! : \Shv(\PP(R, \sS), \tau) \rightleftarrows
    \Shv(\PP(R', \sS'), \tau) : \Phi^*,
  \]
  \[
    \Phi_! : \Shv(\DD(R, \sS), \tau) \rightleftarrows
    \Shv(\DD(R', \sS'), \tau) : \Phi^*,
  \]
  in which the left adjoints $\Phi_!$ are determined
  by the condition that $\Phi_!(\yo K) = \yo (\Phi K)$
  for $K \in \PP(R, \sS)$ (respectively, $K \in \DD(R, \sS)$).

  Our notation $\Phi_! \dashv \Phi^*$ for these geometric morphisms
  differs from that used in topos theory or algebraic geometry
  because our morphisms of o-minimal structures go in the ``algebraic'' direction:
  the morphism $\Phi$ determines a map of rings $\Phi : R \to R'$.
  The operation $\Phi_!$ corresponds to
  pullback of sheaves along the induced map $\Spec \Phi : \Spec R' \to \Spec R$.
\end{nul}

\begin{nul}
  \label{topological-realization}
  Suppose that $R = \RR$.
  In that case, we showed in \cref{quotients-top} that
  if $f : X \to Y$ is a proper surjective map,
  then $f_\tp : X_\tp \to Y_\tp$ is a quotient map in $\Top$.
  We claim that the colimit-preserving extension of
  ${-}_\tp : \DD \to \Top$ to $\PSh(\DD)$
  descends to the category $\Shv(\DD, \taupr)$,
  producing a commutative diagram
  \[
    \begin{tikzcd}
      \DD \ar[r, "\yo"] \ar[rd, "\yo"'] &
      \PSh(\DD) \ar[r] \ar[d] &
      \Top \\
      & \Shv(\DD, \taupr) \ar[ru]
    \end{tikzcd}
  \]
  in which the top composition is ${-}_\tp$
  and the unmarked arrows are left adjoints.
  Since ${-}_\tp$ commutes with finite coproducts,
  it suffices to check that when $f : X \to Y$ is a proper surjective map,
  the diagram $(X \times_Y X)_\tp \rightrightarrows X_\tp \to Y_\tp$
  is a coequalizer in $\Top$;
  and this follows from the aforementioned fact
  that $X_\tp \to Y_\tp$ is a quotient map.

  We denote this induced functor by $|{-}|_\tp : \Shv(\DD, \taupr) \to \Top$
  and refer to it as ``topological realization''.
  Of course, we can also make the same definition with $\DD$ replaced by $\PP$
  or $\taupr$ replaced by the open or closed topology.
\end{nul}

\begin{remark}
  \Cref{proper-sheaves-squares}
  also yields a characterization of the category of sheaves
  by a universal property:
  it is the free cocompletion subject to the condition
  that the colimits mentioned in the last sentence are preserved.
  We will discuss this further at \cref{model-category-universal}.
\end{remark}

\begin{nul*}
  The remainder of this section is devoted to
  proving \cref{proper-sheaves-squares}.
  We start with the easier direction.
\end{nul*}

\begin{lemma}
  \label{proper-squares-of-sheaf}
  If $F$ is a sheaf for the proper topology (on $\DD$ or $\PP$)
  then $F$ satisfies the two conditions of \cref{proper-sheaves-squares}.
\end{lemma}

\begin{proof}
  $F(\emptyset) = *$ is equivalent to
  $F$ satisfying the sheaf condition for the empty cover of $\emptyset$.
  For condition (b), suppose given a distinguished proper square
  \[
    \begin{tikzcd}
      A \ar[r, "f"] \ar[d, "j"'] & A' \ar[d, "j'"] \\
      B \ar[r, "g"] & B'
    \end{tikzcd}
  \]
  and sections $\alpha' \in FA'$ and $\beta \in FB$
  such that $f^* \alpha' = j^* \beta$.
  We must show that there exists a unique section $\beta' \in FB'$
  with $g^* \beta' = \beta$ and $(j')^* \beta' = \alpha'$.
  Uniqueness follows from applying the sheaf property of $F$
  to the proper covering family $\{A' \xrightarrow{j'} B', B \xrightarrow{g} B'\}$ of $B'$.
  Existence will follow as well once we check that,
  extending the above diagram to
  \[
    \begin{tikzcd}
      & A \times_{A'} A \ar[r, shift left, "\pi_1"] \ar[r, shift right, "\pi_2"']
      \ar[d, "j \times j"']
      & A \ar[r, "f"] \ar[d, "j"'] & A' \ar[d, "j'"] \\
      B \ar[r, "\Delta"]
      & B \times_{B'} B \ar[r, shift left, "\pi_1"] \ar[r, shift right, "\pi_2"']
      & B \ar[r, "g"] & B'
    \end{tikzcd}
  \]
  we have $\pi_1^* \beta = \pi_2^* \beta$,
  so that $\alpha'$ and $\beta$ constitute a compatible family of sections
  for this covering family.
  (The corresponding condition for $\alpha'$ is automatic
  because $j' : A' \to B'$ is a monomorphism.)

  We claim that the maps
  $j \times j : A \times_{A'} A \to B \times_{B'} B$
  and $\Delta : B \to B \times_{B'} B$
  form a closed cover of $B \times_{B'} B$.
  Indeed, they are closed embeddings (because $j$ is one)
  and they cover $B \times_{B'} B$ set-theoretically
  because, as a set,
  $B'$ is obtained from $B$ by identifying $j(a)$ and $j(a')$
  if $a$ and $a'$ have the same image under $f$.
  So, by the sheaf condition for this family,
  two sections of $F$ on $B \times_{B'} B$ are equal
  if they agree after pullback along $j \times j$ and along $\Delta$.
  This is the case for the sections $(\pi_1)^* \beta$ and $(\pi_2)^* \beta$
  because the sections
  \[
    (j \times j)^* (\pi_i)^* \beta
    = (\pi_i)^* j^* \beta
    = (\pi_i)^* f^* \alpha'
  \]
  are equal (because $\pi_1 f = \pi_2 f$),
  as are the sections $\Delta^* (\pi_i)^* \beta = \beta$.
\end{proof}

\begin{nul*}
  For the reverse direction,
  we will use the following auxiliary definition.
  (Compare \cite[Definition~2.15]{Voe10_unstable_motivic_homotopy}.)
\end{nul*}

\begin{definition}
  \label{splitting-sequence}
  Let $f : Y \to X$ be a surjective map of $\DD$.
  A \emph{splitting sequence} for $f$
  is a finite filtration
  $\emptyset = Z_{n+1} \subset Z_n \subset \cdots \subset Z_0 = X$
  of $X$ by closed subsets
  such that for each $0 \le i \le n$,
  the pullback of $f$ to $Z_i - Z_{i+1} \subset X$ admits
  a (continuous definable) section.
  Here $n$ is the \emph{length} of the splitting sequence.
  (We allow $n = -1$ if $X = \emptyset$.)
\end{definition}

\begin{proposition}
  \label{splitting-seq-of-surjective}
  Every surjective map $f : Y \to X$ of $\DD$ admits a splitting sequence;
  in fact, one of length at most $\dim(X)$.
  (By convention, we define $\dim(\emptyset) = -1$.)
\end{proposition}

\begin{proof}
  We prove the statement by induction on $d = \dim(X)$.
  When $d = -1$, $X$ is empty and we can take
  the unique splitting sequence of length $-1$.
  Suppose $d \ge 0$ and let $f : Y \to X$ be a surjective map.
  By \cref{definable-choice}, $f$ has
  a definable (but not necessarily continuous) section $g : X \to Y$.
  By \cref{continuous-on-partition},
  $X$ admits a partition into finitely many definable subsets $(S_i)_{i \in I}$
  such that $g$ is continuous when restricted to each $S_i$.
  Let $U_i$ be the interior (in $X$) of $S_i$ for each $i$, and
  let $U$ be the union of all the $U_i$.
  Then $g$ is a continuous definable section of $f$ on $U$.
  Set $Z_1 = X - U$.
  Then $Z_1$ is the union of the finitely many sets $S_i - U_i$
  and so
  $\dim(Z_1) = \max_{i \in I} \dim(S_i - U_i) < \dim(X)$
  by \cref{dim-union,dim-sdiff-interior-lt}.
  By the inductive hypothesis,
  we can find a splitting sequence for the pullback of $f$ to $Z_1$
  of length at most $\dim(Z_1) \le d-1$.
  Then appending $Z_0 = X$ to this sequence
  yields a splitting sequence for $f$ of length at most $d$.
\end{proof}

\begin{proposition}
  \label{section-of-splitting}
  Let $f : Y \to X$ be a proper map
  and $U \subset X$ an open subset with closed complement $Z$
  such that the pullback of $f$ to $U$ has a section $g : U \to Y$.
  Construct the pullback
  \[
    \begin{tikzcd}
      Y \times_X Z \ar[r] \ar[d] & Y \ar[d, "f"] \\
      Z \ar[r, "j"] & X
    \end{tikzcd}
  \]
  and then the (closed-by-proper) pushout
  \[
    \begin{tikzcd}
      Y \times_X Z \ar[r] \ar[d] & Y \ar[d, "q"] \\
      Z \ar[r, "k"] & Y'
    \end{tikzcd}
  \]
  obtaining an induced map $f' : Y' \to X$.
  Then $f'$ has a section $h : X \to Y'$
  such that $hj = k$.
\end{proposition}

\begin{proof}
  Intuitively, $Y'$ is obtained from $Y$
  by collapsing the fiber over each $x \in Z$ to a single point.
  The idea of the proof is to show that
  the composition $U \xrightarrow{g} Y \xrightarrow{q} Y'$
  extends continuously to all of $X$
  by sending each $x \in Z$ to the unique point of $Y'$ representing its fiber.

  Accordingly, we construct $h : X \to Y'$ as follows:
  \begin{itemize}
  \item
    On $U$, we define $h$ by the composition
    $U \xrightarrow{g} Y \xrightarrow{q} Y'$.
  \item
    On $Z$, we define $h$ by $Z \xrightarrow{k} Y'$.
  \end{itemize}
  This function $h$ is clearly definable and a section of $f'$
  (indeed, the unique section extending $g$).
  We need to prove it is continuous.
  By \cref{continuous-iff-curve},
  it suffices to prove that
  for any curve $\gamma : [0, 1] \to X$,
  the composition $h \circ \gamma : [0, 1] \to Y'$
  is continuous at $0$.
  We divide into cases depending on whether $\gamma(0)$ belongs to $Z$ or to $U$
  and whether $\gamma(t)$ belongs to $Z$ or to $U$ for small positive $t$
  (one of these cases necessarily holds by \cref{continuous-on-partition}).
  \begin{itemize}
  \item
    $\gamma(t)$ belongs to $Z$ or to $U$ for sufficiently small $t$,
    including $0$.
    Then the claim holds because $h$ is continuous
    when restricted to $Z$ or to $U$.
  \item
    $\gamma(0)$ belongs to $U$ but
    $\gamma(t)$ belongs to $Z$ for sufficiently small positive $t$.
    This case is impossible because $U$ is open.
  \item
    $\gamma(0)$ belongs to $Z$ but
    $\gamma(t)$ belongs to $U$ for sufficiently small positive $t$,
    say $t \le t_0$.
    In this case, we will show that $h \circ \gamma : [0, t_0] \to Y'$
    can be expressed as the image under $q : Y \to Y'$
    of a curve $\gamma' : [0, t_0] \to Y$;
    then $h \circ \gamma$ will be continuous.
    Define $\gamma'(t) = g(\gamma(t))$ for $0 < t \le t_0$.
    Then the square
    \[
      \begin{tikzcd}
        (0, t_0] \ar[r, "\gamma'"] \ar[d] & Y \ar[d, "f"] \\
        {[0, t_0]} \ar[r, "\gamma"] & X
      \end{tikzcd}
    \]
    commutes so, since $f$ is proper,
    $\gamma'$ extends (uniquely) to a map $\gamma' : [0, t_0] \to Y$.
    Now $q \circ \gamma' = h \circ \gamma$
    because, for $0 < t \le t_0$,
    \[
      q(\gamma'(t)) = q(g(\gamma(t))) = h(\gamma(t))
      \quad
      \mbox{(since $\gamma(t) \in U$)}
    \]
    while $q(\gamma'(0))$ equals $h(\gamma(0)) = k(\gamma(0))$
    because both are the unique point of $Y'$ lying above $\gamma(0) \in Z$.
    \qedhere
  \end{itemize}
\end{proof}

\begin{proof}[Proof of \cref{proper-sheaves-squares}]
  We have already shown one direction (\cref{proper-squares-of-sheaf}).
  Suppose, conversely, that $F$ is a presheaf (on $\DD$ or $\PP$)
  satisfying the conditions:
  \begin{enumerate}
  \item[(a)]
    $F(\emptyset) = *$.
  \item[(b)]
    $F$ sends each pushout square
    \[
      \begin{tikzcd}
        A \ar[r, "f"] \ar[d, "j"'] & A' \ar[d] \\
        B \ar[r] & B'
      \end{tikzcd}
    \]
    with $j$ a closed embedding and $f$ proper
    to a pullback square.
  \end{enumerate}
  Then (by (a) and (b) with $A = \emptyset$)
  $F$ takes finite coproducts to products,
  so it satisfies the sheaf condition for finite coproducts.
  It suffices then to show that
  $F$ satisfies the sheaf condition for a single surjective proper map.

  Let $\mathcal{C}$ denote the class of maps
  for which $F$ satisfies the sheaf condition.
  Then $\mathcal{C}$ enjoys the following properties.
  \begin{itemize}
  \item
    $\mathcal{C}$ is closed under composition
    \cite[C2.1.7]{Elephant2}.
  \item
    $\mathcal{C}$ contains any morphism which admits a section.
    (Such a morphism generates the sieve consisting of
    all maps into its codomain.)
  \item
    Suppose that
    \[
      \begin{tikzcd}
        A \ar[r, "f"] \ar[d, "j"'] & A' \ar[d, "j'"] \\
        B \ar[r, "g"] & B'
      \end{tikzcd}
    \]
    is a pushout square with $j$ a closed inclusion
    and $f$ a proper map that belongs to $\mathcal{C}$.
    Then we claim that $g$ also belongs to $\mathcal{C}$.
    Consider the diagram
    \[
      \begin{tikzcd}
        FB' \ar[r, "g^*"] \ar[d, "(j')^*"']
        & FB \ar[r, shift left, "\pi_1^*"] \ar[r, shift right, "\pi_2^*"']
        \ar[d, "j^*"]
        & F(B \times_{B'} B) \ar[d] \\
        FA' \ar[r, "f^*"]
        & FA \ar[r, shift left, "\pi_1^*"] \ar[r, shift right, "\pi_2^*"']
        & F(A \times_{A'} A)
      \end{tikzcd}
    \]
    in which the square on the left is a pullback, by condition (b).
    Since $F$ satisfies the sheaf condition for $f$ by the hypothesis on $f$,
    the bottom row is an equalizer.
    Then by a standard diagram chase, the top row is also an equalizer.
  \end{itemize}
  Now, we prove that every surjective proper map $f : Y \to X$
  belongs to $\mathcal{C}$
  by induction on the length of a splitting sequence for $f$
  (which exists by \cref{splitting-seq-of-surjective}).
  If the splitting sequence has length $-1$
  then $X$ is empty and $f$ is an isomorphism.
  So suppose the splitting sequence has length $n \ge 0$.
  Then it ends with a closed subset $Z \subset X$
  for which the pullback $f_1$ of $f$ to $Z$
  (again a surjective proper map)
  has a splitting sequence of length $n-1$
  and the pullback of $f$ to $U = X - Z$ has a section.
  By the inductive hypothesis, $f_1$ belongs to $\mathcal{C}$.
  By \cref{section-of-splitting},
  $f : Y \to X$ can be expressed as
  the composition of the pushout of $f_1$
  along the closed embedding $Y \times_X Z \subset Y$
  and a map $f'$ which admits a section.
  By the properties of $\mathcal{C}$ listed above,
  $f$ then also belongs to $\mathcal{C}$,
  completing the induction.
\end{proof}

\begin{remark}
  We have described the squares appearing in condition (b)
  as closed-by-proper pushouts
  in order to emphasize their relationships
  to effective quotients, cofibration categories,
  and the pushouts used to build up
  the realization of a simplicial set.
  However, there is a second, equivalent description
  of these squares
  which is precisely analogous to the definition of
  the proper cdh-structure of \cite{Voe10_unstable_motivic_homotopy}.
  We adjust our notation to match \cite{Voe10_unstable_motivic_homotopy}.
\end{remark}

\begin{proposition}
  \label{closed-by-proper-iff-blowup}
  A square in $\DD$
  \[
    \begin{tikzcd}
      B \ar[r, "i"] \ar[d, "e"'] & Y \ar[d, "f"] \\
      A \ar[r, "j"] & X
    \end{tikzcd}
  \]
  is a closed-by-proper pushout
  if and only if it is a pullback square
  such that $f$ is proper, $j$ is a closed embedding
  and the pullback of $f$ to $X - j(A)$
  is an isomorphism.
\end{proposition}

\begin{proof}
  The ``only if'' direction follows from \cref{closed-by-proper-pushout},
  the last part holding because it does in $\Set$
  and a proper bijection is an isomorphism.
  For the ``if'' direction,
  $e$ is proper and $i$ a closed embedding
  because they are pullbacks of $f$ and $j$ respectively.
  Thus, we may form a pushout
  \[
    \begin{tikzcd}
      B \ar[r, "i"] \ar[d, "e"'] & Y \ar[d, "e'"'] \ar[rdd, "f", bend left] \\
      A \ar[r, "i'"] \ar[rrd, "j"', bend right] & X' \ar[rd, "f'"'] \\
      & & X \ar[ul, dashed, bend right, "h"']
    \end{tikzcd}
  \]
  and we must prove that the induced map $f' : X' \to X$ is an isomorphism.
  Since $f$ is proper, \cref{section-of-splitting} implies that
  there is a map $h : X \to X'$ (as shown by the dotted arrow)
  such that $f'h = \id_L$ and $hj = i'$.
  We need to show that $hf'$ equals the identity of $X'$,
  which can be checked pointwise on~$X'$.
  Let $p \in X'$.

  If $f'(p) \in j(A)$, we claim first that $p \in i'(A)$.
  The other possibility is that $p = e'(y)$ for some $y \in Y$,
  but then $f(y) = f'(p) \in j(A)$ and so as the original square was a pullback,
  there exists $b \in B$ with $i(b) = y$.
  Then $p = e'(i(b)) = i'(e(b))$ is again in the image of $i'$.
  So, in either case, we can write $p = i'(a)$ for some $a$,
  and then $h(f'(p)) = h(j(a)) = i'(a) = p$.

  Now suppose $f'(p) \notin j(A)$, and let $q = h(f'(p)) \in X'$.
  Then $f'(q) = (f'h)(f'(p)) = f'(p)$ and so
  both $f'(p)$ and $f'(q)$ lie in $U := X - j(A)$.
  Then $p$ and $q$ must be images of points $y$ and $z$ of $Y$ respectively
  which satisfy $f(y) = f(z) \in U$.
  Since $f$ is an isomorphism when pulled back to $U$,
  it follows that $y = z$ and therefore $p = q = h(f'(p))$.
\end{proof}

\begin{remark}
  In $\DD$,
  open coverings are also covers for the closed topology
  (\cref{closed-finer-than-open})
  and hence also the proper topology.
  Moreover,
  since the implicit function theorem is valid in an o-minimal structure
  \cite[7.2.11]{vdD98},
  one would assume that
  the analogue of the \'etale maps in $\DD$
  are just the local homeomorphisms,
  which are also covers for the proper topology.
  Combining these facts with the observation that
  the ``completely decomposed'' condition is vacuous in $\DD$,
  we might regard the proper topology on $\DD$
  as analogous to pretty much any of the commonly-used variants
  of the h-topology in algebraic geometry.
\end{remark}

%% file: tex/p_vs_d.tex
\section{\texorpdfstring
{The relationship between $\Shv(\DD)$ and $\Shv(\PP)$}
{The relationship between Shv(D) and Shv(P)}}
\label{p_vs_d}

\begin{nul}
  Recall that we regard $\DD$ and $\PP$ as sites
  by equipping them with the proper topology
  unless specified otherwise.
  The inclusion functor $u : \PP \to \DD$ has the following properties:
  \begin{itemize}
  \item
    $u$ is fully faithful.
  \item
    $u$ preserves finite limits.
  \item
    $u$ preserves covering families:
    if $(L_i \to K)_{i \in I}$ is a covering family in $\PP$
    then its image $(uL_i \to uK)_{i \in I}$ is also a covering family in $\DD$.
  \item
    For $K \in \PP$, any covering family $(L_i \to uK)_{i \in I}$ in $\DD$
    is refined by the image under $u$ of a covering family in $\PP$.
    (This property is called ``cocontinuity'' in \cite{SGA4}.)
  \end{itemize}
  (The last two properties can be checked on the generating covering families.
  For the last one, if $K \in \PP$ and $L_i \to uK$ is a proper map,
  then $L_i$ also belongs to $\PP$.)

  These properties have the following consequences.
\end{nul}

\begin{proposition}
  The adjoint triple $u_! \dashv u^* \dashv u_*$
  of functors between $\PSh(\PP)$ and $\PSh(\DD)$
  induces an adjoint triple
  \[
    u_! : \Shv(\PP) \to \Shv(\DD)
  \]
  \[
    u^* : \Shv(\DD) \to \Shv(\PP)
  \]
  \[
    u_* : \Shv(\PP) \to \Shv(\DD)
  \]
  in which $u_!$ is fully faithful and preserves finite limits.
\end{proposition}

\begin{proof}
  Standard.
\end{proof}

\begin{nul}
  We will frequently identify $\Shv(\PP)$
  with the essential image of
  the fully faithful functor $u_! : \Shv(\PP) \to \Shv(\DD)$.
  Then $\Shv(\PP)$ is a coreflective subcategory of $\Shv(\DD)$,
  with coreflector $u^*$
  and associated comonad $u_! u^*$ on $\Shv(\DD)$.
  The subcategory $\Shv(\PP)$ is
  closed under arbitrary colimits and finite limits;
  it does not seem to be closed under infinite products.
  An object of $\Shv(\DD)$ belongs to $\Shv(\PP)$
  if and only if it can be written as a colimit of
  objects of the form $\yo K$ with $K \in \PP$.
\end{nul}

\begin{nul}
  The full subcategory $\PP \subset \DD$
  also has the following ``closure'' property
  with respect to the Grothendieck topologies.
  \begin{itemize}
  \item
    Suppose $(L_i \to K)_{i \in I}$ is a covering family of $K \in \DD$
    by objects $L_i$ which belong to $\PP$.
    Then $K$ also belongs to $\PP$.
  \end{itemize}
  The proof is simple:
  Because the generating covering families of $\DD$ are finite,
  some finite subfamily of $(L_i \to K)_{i \in I}$ must already cover $K$.
  Then the definable set $K$ is
  a finite union of images of closed and bounded sets,
  so $K$ is also closed and bounded.

  This has the following consequence for the categories of sheaves.
\end{nul}

\begin{proposition}
  The square of fully faithful functors
  \[
    \begin{tikzcd}
      \PP \ar[r, "\yo"] \ar[d, "u"'] &
      \Shv(\PP) \ar[d, "u_!"] \\
      \DD \ar[r, "\yo"] &
      \Shv(\DD)
    \end{tikzcd}
  \]
  is a (pseudo)pullback.
  In other words, an object of $\Shv(\DD)$
  is representable by an object of $\PP$
  if and only if
  it belongs to $\Shv(\PP)$
  and is representable by an object of $\DD$.
\end{proposition}

\begin{proof}
  The ``only if'' direction is obvious,
  so suppose $K \in \DD$ is such that $\yo K \in \Shv(\DD)$
  belongs to $\Shv(\PP)$.
  Then $\yo K$ can be written as a colimit of
  objects of the form $\yo L_i$ with $L_i \in \PP$,
  so in particular there is a cover $\coprod_{i \in I} \yo L_i \to \yo K$.
  Using the ``closure'' property noted above
  we conclude that $K$ must belong to $\PP$.
\end{proof}

\begin{nul}
  \label{r_loc}
  Let $M \in \DD$.
  We describe the coreflection $u_! u^* \yo M$.
  Let $P_M$ denote the family of definable subsets $L$ of $M$
  which are objects of $\PP$
  (that is, closed and bounded in the ambient space $R^n$ containing $M$),
  partially ordered by inclusion.
  Each $L \in P_M$ is equipped with an inclusion map $L \to M$,
  which induces a morphism $\yo L \to \yo M$.
  We set $M_\loc = \colim_{L \in P_M} \yo L$;
  then $M_\loc$ is equipped with an induced map $M_\loc \to \yo M$.
  We claim $M_\loc \to \yo M$ is the counit map $u_! u^* \yo M \to \yo M$.
  Indeed, $M_\loc$ belongs to $\Shv(\PP)$ by construction,
  and if $K \in \PP$ then any map $f : K \to M$ in $\DD$
  factors through the inclusion of the image $f(K) \subset M$,
  an object of $\PP$.
  Because the poset $P_M$ is filtered
  we conclude that the induced map
  \[
    \Hom(\yo K, M_\loc) = \colim_{L \in P_M} \Hom(K, L)
    \to \Hom(K, M) = \Hom(\yo K, \yo M)
  \]
  is an isomorphism.

  This description of $M_\loc = u_! u^* \yo M$
  as a filtered colimit of representables
  has the advantage of being canonical,
  but the poset $P_M$ is quite large.
  For geometric purposes it can be useful to
  replace $P_M$ by a smaller cofinal subset.
  We will consider just the case $M = R$.
  In this case, $P_R$ has a cofinal subset
  consisting of the intervals $[-r, r]$
  as $r$ ranges over all positive elements of $R$,
  and so $R_\loc = \colim_{r > 0} \yo [-r, r]$.
  To further reduce the size of the indexing diagram,
  we define the \emph{cofinality} of the field $R$
  to be the smallest cardinality of an unbounded subset of $R$.
  Then, if $R$ has cofinality $\kappa$,
  we can express $R_\loc$ as
  the filtered colimit of a diagram of size $\kappa$
  consisting of representables.
  (A more complicated construction applies to any $M \in \DD$
  which is ``locally complete'',
  meaning that each $x \in M$ has a neighborhood in $M$ which is an object of $\PP$.
  Compare \cite[Example~I.2.6]{LSA}.)
  Conversely, suppose $R_\loc = \colim_{i \in I} \yo K_i$ is
  a construction of $R_\loc$ as a colimit of representables with $K_i \in \PP$.
  Then each map $\yo K_i \to R_\loc \to \yo R$ induces
  a definable continuous map $K_i \to R$.
  Each of these maps has bounded image, and they must be jointly surjective;
  this is only possible if $|I| \ge \kappa$.
  In particular, $R_\loc$ can be written as
  the colimit of a countable sequence of closed embeddings
  $[-r_0, r_0] \subset [-r_1, r_1] \subset \cdots$
  if and only if the field $R$ is \emph{sequential} (has countable cofinality).

  The functor $u_! u^*$ has a right adjoint $u_* u^*$,
  so its value on a general object~$X$ of~$\DD$
  can be computed by writing $X$ as a colimit of representables
  and applying the above construction to each of these representables.
\end{nul}

\begin{remark}
  We have borrowed the notation $M_\loc$ from \cite{LSA},
  suggesting the word ``locally''.
  However, the adverb ``locally'' must be understood as meaning
  ``on definably compact subsets'',
  not in the sense involving open neighborhoods.
  In the nonarchimedean case these notions are different.
  For example the function which takes value $1$ on infinitesimals
  and $0$ elsewhere is locally definable (even constant) in the latter sense
  but not definable on the interval $[-1, 1]$.
\end{remark}

\begin{warning}
  Let $\Phi : (R, \sS) \to (R', \sS')$ be a morphism of o-minimal structures,
  and write $\PP'$ and $\DD'$ for the categories associated to the structure $(R', \sS')$.
  Then $\Phi$ induces base change functors
  \[
    \Phi_! : \Shv(\PP) \to \Shv(\PP')
    \quad\mbox{and}\quad
    \Phi_! : \Shv(\DD) \to \Shv(\DD').
  \]
  The square
  \[
    \begin{tikzcd}
      \Shv(\PP) \ar[r, "u_!"] \ar[d, "\Phi_!"'] &
      \Shv(\DD) \ar[d, "\Phi_!"] \\
      \Shv(\PP') \ar[r, "u_!"] &
      \Shv(\DD')
    \end{tikzcd}
  \]
  commutes (up to isomorphism), but the square
  \[
    \begin{tikzcd}
      \Shv(\PP) \ar[d, "\Phi_!"'] &
      \Shv(\DD) \ar[l, "u^*"'] \ar[d, "\Phi_!"] \\
      \Shv(\PP') &
      \Shv(\DD') \ar[l, "u^*"']
    \end{tikzcd}
  \]
  generally does not.
  For example, when $R'$ is a field of larger cofinality than $R$,
  this should be clear from the description of $R_\loc = u_! u^* \yo R$ above.
  (Compare \cite[Example~I.2.11]{LSA}.)
\end{warning}

\begin{example}
  In order to clarify the difference between $R_\loc$ and $\yo R$,
  we describe the Hom sets between the four possible pairs of these objects.
  Of course, a map $f : \yo R \to \yo R$ corresponds to a map $f : R \to R$ in $\DD$,
  that is, a continuous definable function from $R$ to itself.
  A map $f : \yo R \to R_\loc = \colim_{r > 0} {[-r, r]}$ corresponds to
  a continuous definable function $f : R \to R$ \emph{with bounded image},
  since $\Hom(\yo R, -)$ commutes with filtered colimits.
  On the other hand, a map $f : R_\loc \to \yo R$ corresponds to
  a function $f : R \to R$ which is \emph{continuous and definable when restricted to each $[-r, r]$}.
  The continuity condition is equivalent to $f$ being continuous on all of $R$,
  but the definability condition is strictly weaker, at least when $R$ is sequential.
  For example, if $r_1 < r_2 < \cdots$ is an unbounded sequence in $R$,
  then $f(x) = \min_i |x - r_i|$ is definable on every bounded interval
  but not globally definable (since $f(x) = 0$ has an infinite discrete set of solutions).
  We do not know whether such examples can be constructed
  when $R$ is not sequential.
  Finally, the maps $R_\loc \to R_\loc$ are the same as the maps $R_\loc \to R$
  because $\Shv(\PP)$ is a coreflective subcategory of $\Shv(\DD)$
  (or because every continuous definable function on a closed interval is bounded).

  The implication is that if we work in $\Shv(\PP)$,
  the best available approximation to the object $\yo R \in \Shv(\DD)$ is $u^* \yo R$,
  which we can identify with~$R_\loc$.
  Consequently, in $\Shv(\PP)$ we have no way to describe
  the continuous and globally definable functions $f : R \to R$.
  Instead, we must make do with the functions which are
  continuous and definable on bounded intervals.
  Depending on the application, this might be either an advantage or a disadvantage.
\end{example}

\begin{example}
  \label{exp-loc}
  The standard universal covering map
  $f : \RR \to S^1 \subset \CC = \RR^2$
  given by $f(t) = \exp (2\pi i t)$
  cannot be definable in any o-minimal structure,
  because $\{\,t \in \RR \mid f(t) = 1\,\} = \ZZ \subset \RR$
  is an infinite discrete set.
  We will have the same problem for any other universal covering map
  $g : \RR \to S^1$ as well,
  since $\{\,t \in \RR \mid g(t) = 1\,\}$ must be discrete and infinite.

  We can sidestep this problem by working in $\Shv(\PP)$,
  or by working in $\Shv(\DD)$ but with the object $\RR_\loc$.
  For example, there is an o-minimal structure $\RR_\an$
  which contains all restrictions of analytic functions to compact intervals
  (see \cref{real-examples} for details).
  Then the function $f(t) = \exp (2\pi i t)$
  does determine a morphism $f : \RR_\loc \to \yo S^1$.
  More generally, let $G$ be a classical Lie group,
  considered as a group object of $\DD$
  by embedding it in some $\RR^{N \times N}$ in the standard way,
  and let $\mathfrak{g} \subset \RR^{N \times N}$ be its Lie algebra.
  Then the exponential map $\exp : \mathfrak{g} \to G$
  becomes definable in this structure
  when viewed as a morphism from~$\mathfrak{g}_\loc$ to~$G_\loc$.

  If we want to work in a smaller o-minimal structure such as $\Qalgre_\sa$,
  we can still choose a semialgebraic parametrization
  $g : [0, 1]/\{0 \sim 1\} \to S^1$
  and extend it by periodicity to all of $\Qalgre$.
  On each interval $[-N, N]$ of $\Qalgre$,
  the extended function is definable
  (by a formula whose complexity grows with $N$)
  and so we again obtain a universal covering map in the category $\Shv(\PP)$.
  (However, this universal covering map cannot be a group homomorphism
  when the trigonometric functions are not definable on compact intervals.)
\end{example}

\begin{speculation}
  \label{proper-homotopy-theory}
  We saw in \cref{ho-p-eq-ho-d} that
  every object of $\DD$ is homotopy equivalent to an object of $\PP$
  and so working in $\Shv(\DD)$ seems to confer no particular advantage
  for the purpose of the classical homotopy theory of spaces,
  which in any case is normally presented by
  building up spaces from pieces that are closed and bounded
  (be they closed balls in the case of CW complexes,
  or simplices in the case of simplicial sets,
  or cubes, \ldots).

  However, we speculate that the extra information contained in $\Shv(\DD)$
  may be useful for some form of \emph{proper homotopy theory}.
  Heuristically, the difference between
  a representable object of $\Shv(\DD)$ and its restriction to $\Shv(\PP)$
  is that the former retains information about
  definable paths that ``run off to infinity''.
  On $\PP$, there is a unique numerical invariant (the Euler characteristic)
  which is homotopy invariant, additive on pushouts by closed embeddings,
  and takes the value $1$ on a point.
  However, to a general definable set $K \in \DD$
  we can associate \emph{two} numerical invariants:
  \begin{itemize}
  \item
    The (o-minimal, \cite[\S4.2]{vdD98}) Euler characteristic $\chi(K)$,
    which agrees with the ordinary Euler characteristic on polytopes
    and is additive with respect to arbitrary partitions
    into finitely many definable sets.
    It is not homotopy invariant in general:
    if $K = [0, 1)$ then $\chi(K) = \chi([0, 1]) - \chi(\{1\}) = 1 - 1 = 0$.
    It is however a proper homotopy invariant
    at least when restricted to locally closed definable sets.
  \item
    The homotopy invariant Euler characteristic $\chi^h(K)$,
    defined as the Euler characteristic of
    a polytope homotopy equivalent to $K$
    (so that for $K = [0, 1)$, we have $\chi^h(K) = 1$).
  \end{itemize}
  These two invariants suggest that
  $\Shv(\DD)$ should naturally be a setting
  for some more refined homotopy theory
  than the classical homotopy theory of spaces.
\end{speculation}

%% file: tex/modelcat.tex
\section{Model category structures}
\label{modelcat}

\begin{nul*}
  In this section,
  we will construct ``Serre--Quillen-style'' model category structures
  on the categories $\Shv(\PP)$ and $\Shv(\DD)$.
  Our approach is as follows.
  \begin{itemize}
  \item
    We introduce the geometric realization functor
    $|{-}| : \sSet \to \Shv(\PP)$
    and describe its relationship to
    the geometric realization of finite simplicial sets
    discussed in \cref{geometric-realization}.
    The geometric realization functor is left adjoint to a functor
    $\Sing : \Shv(\PP) \to \sSet$.
  \item
    Following the original approach of Quillen~\cite{Qui67},
    we obtain a transferred model category structure on $\Shv(\PP)$
    in which the fibrations and weak equivalences are detected by $\Sing$.
    We further transfer this model category structure to $\Shv(\DD)$
    along the adjunction
    $u_! : \Shv(\PP) \rightleftarrows \Shv(\DD) : u^*$.
  \item
    To show that
    $|{-}| : \sSet \rightleftarrows \Shv(\PP) : \Sing$
    is a Quillen equivalence,
    we check that the unit $K \to \Sing |K|$ is a weak equivalence
    for each $K \in \sSet$
    by reducing to the case of $K$ finite,
    which was treated in \cref{homotopy-d-p}.
  \end{itemize}

  We expect that it should also be possible
  to construct these model category structures
  without making use of the theory of simplicial sets,
  by instead imitating
  a self-contained construction of the model structure on $\Top$,
  as found for example in \cite{Hir19}.
  However, since we intend to show that
  the geometric realization--$\Sing$ adjunction
  is a Quillen equivalence anyways,
  here we will just present the easier construction
  of these model category structures
  by transferring the Kan--Quillen model category structure on $\sSet$.
\end{nul*}

\begin{nul}
  We begin with geometric realization.
  Recall from \cref{geometric-realization}
  the o-minimal simplex functor $|\Delta^\bullet|_\PP : \Delta \to \PP$.
  By composing with the Yoneda embedding we obtain a functor
  \[
    |\Delta^\bullet|_{\Shv(\PP)} :
    \Delta \xrightarrow{|\Delta^\bullet|_\PP} \PP
    \xrightarrow{\yo} \Shv(\PP)
  \]
  and we define $\Sing : \Shv(\PP) \to \sSet$ by the formula
  \[
    (\Sing X)_n = \Hom(|\Delta^n|_{\Shv(\PP)}, X)
    = X(|\Delta^n|_\PP).
  \]
  Since $\Shv(\PP)$ is a cocomplete category,
  the functor $\Sing$ certainly has a left adjoint,
  which we provisionally denote by $|{-}|_{\Shv(\PP)}$.

  Now we have the following fact,
  which is the main point of the category $\Shv(\PP)$.
\end{nul}

\begin{proposition}
  \label{geometric-realization-sset}
  The square
  \[
    \begin{tikzcd}
      \sSetfin \ar[r, "|{-}|_\PP"] \ar[d] & \PP \ar[d, "\yo"] \\
      \sSet \ar[r, "|{-}|_{\Shv(\PP)}"] & \Shv(\PP)
    \end{tikzcd}
  \]
  commutes (up to a canonical natural isomorphism).
\end{proposition}

\begin{proof}
  By \cref{proper-sheaves-squares},
  the functor $\yo : \PP \to \Shv(\PP)$
  preserves the initial object and closed-by-proper pushouts.
  By \cref{finite-realizations}, then,
  $\yo$ sends the $\PP$-valued realization of a finite simplicial set of $K$
  to a $(\yo \circ |\Delta^\bullet|_\PP)$-realization of $K$,
  which can be canonically identified with $|K|_{\Shv(\PP)}$.
\end{proof}

\begin{remark}
  If we had instead equipped $\PP$ with the closed topology $\taucl$
  (\cref{closed-topology}),
  \cref{geometric-realization-sset} would not hold;
  we would only obtain a natural isomorphism
  $\yo |K|_\PP \cong |K|_{\Shv(\PP)}$
  for finite simplicial \emph{complexes} $K$,
  as discussed in \cref{quotient-not-closed-cover}.
  Likewise, \cref{geometric-realization-sset} depends indirectly on
  the definability of the multiplication map,
  as explained in \cref{multiplication-required}.
\end{remark}

\begin{nul}
  We may enlarge the square of \cref{geometric-realization-sset}
  to a larger commutative diagram
  \[
    \begin{tikzcd}
      \sSetfin \ar[r, "|{-}|_\PP"] \ar[d] &
      \PP \ar[r, "u"] \ar[d, "\yo"] &
      \DD \ar[d, "\yo"] \\
      \sSet \ar[r, "|{-}|_{\Shv(\PP)}"] &
      \Shv(\PP) \ar[r, "u_!"] &
      \Shv(\DD)
    \end{tikzcd}
  \]
  and, by \cref{finite-realizations},
  the $\DD$-valued geometric realization $|{-}|_\DD : \sSetfin \to \DD$
  is canonically naturally isomorphic to
  the composition $u \circ |{-}|_\PP$.
  Since all these notions of geometric realization
  agree up to the inclusions of full subcategories,
  \emph{we henceforth drop the subscripts on all these functors},
  writing $|K|$ for an object of
  $\PP$, $\DD$, $\Shv(\PP)$, or $\Shv(\DD)$,
  to be disambiguated by the context.
\end{nul}

\begin{proposition}
  \label{geom-real-lex}
  Geometric realization preserves finite limits.
\end{proposition}

\begin{proof}
  By \cite[section~VIII.8]{MacLane_Moerdijk},
  $\sSet$ is the classifying topos of
  the theory of intervals equipped with distinct bottom and top elements.
  We claim that the object $I = \yo [0, 1]$ of $\Shv(\PP)$,
  equipped with its internal order $\yo Q \rightrightarrows \yo [0, 1]$
  where $Q = \{\,(x, y) \in [0, 1]^2 \mid x \le y\,\}$,
  is a model of this theory.
  The only axiom which requires comment is
  the totality axiom $\forall x, y \, (x \le y \vee y \le x)$;
  it holds because $Q$ and
  its transpose $\{\,(x, y) \in [0, 1]^2 \mid y \le x\,\}$
  form a covering family of $[0, 1]^2$ in $\PP$.
  The functor $\Delta^\bullet_I : \Delta \to \Shv(\PP)$
  associated to this interval object $I$
  is easily seen to coincide with
  the o-minimal simplex functor $\yo \circ |\Delta^\bullet|_\PP$
  we have defined directly,
  and so by \cite[Proposition~VIII.8.4]{MacLane_Moerdijk}
  the associated geometric realization functor is left exact.
\end{proof}

\begin{nul}
  In particular, the promise made in \cref{finite-realization-finite-limits}
  is now fulfilled:
  the $\PP$-valued realization of finite simplicial sets
  preserves finite limits,
  because the inclusion $\sSetfin \subset \sSet$
  and the Yoneda embedding $\yo : \PP \to \Shv(\PP)$
  preserve and reflect finite limits.
\end{nul}

\begin{proposition}
  \label{sing-fibrant}
  For any object $X \in \Shv(\PP)$,
  the simplicial set $\Sing X$ is a Kan complex.
\end{proposition}

\begin{proof}
  It suffices to check that the map $X \to *$
  has the right lifting property with respect to
  the geometric realization of each horn inclusion $|\Lambda^n_i| \to |\Delta^n|$
  ($0 < n$, $0 \le i \le n$),
  and in view of \cref{geometric-realization-sset},
  we already checked in \cref{horn-inclusion-retract}
  that these geometric realizations admit retracts.
\end{proof}

\begin{nul*}
  We are now in a position to apply standard results
  on the existence of transferred model structures.
  For our purposes, the following criterion will suffice.
\end{nul*}

\begin{proposition}
  \label{monoidal-transfer}
  Let $F : M \rightleftarrows N : U$ be
  a monoidal adjunction (i.e., with $F$ strong monoidal)
  between closed symmetric monoidal locally presentable categories
  and suppose $M$ is equipped with
  a monoidal combinatorial model category structure with cofibrant unit.
  Assume that $UX$ is fibrant in $M$ for every object~$X$ of~$N$.
  Then there is a transferred model structure on $N$
  in which a map $f$ of $N$ is a weak equivalence (respectively, fibration)
  if and only if $Uf$ is a weak equivalence (respectively, fibration) of $M$.
\end{proposition}

\begin{remark}
  \label{transferred-features}
  The transferred model structure also has the following properties
  which are consequences of the ones described above.
  \begin{itemize}
  \item
    The original adjunction $F : M \rightleftarrows N : U$
    becomes a Quillen adjunction
    for the original model structure on $M$
    and the transferred model structure on $N$.
    Furthermore, a left adjoint $G : N \to P$ to another model category
    is a left Quillen functor if and only if
    the composition $GF : M \to P$ is one.
  \item
    Generating (acyclic) cofibrations for $N$
    can be obtained by applying $F$ to those for $M$.
    In particular, $N$ is again combinatorial
    and the transferred model structure on $N$ is automatically monoidal
    because $F$ commutes with pushout products.
  \item
    The monoidal adjunction $F : M \rightleftarrows N : U$
    gives $N$ the structure of a category
    enriched, tensored, and cotensored over $M$.
    In particular, the tensor action of $M$ on $N$
    is given by the formula $K \otimes X = FK \otimes_N X$
    for $K \in M$, $X \in N$.
    This structure is compatible with the transferred model structure,
    making $N$ into an $M$-enriched model category.
  \end{itemize}
\end{remark}

\begin{proof}
  In view of the hypothesis that $UX$ is fibrant for every $X$,
  the existence of the transferred model structure
  will follow by a standard argument
  once we verify the following condition:
  \begin{itemize}
  \item
    For every object $X$ of $N$,
    there exists a ``path factorization''
    $X \xrightarrow{r} PX \xrightarrow{q} X \times X$
    of the diagonal map of $X$
    such that $Ur$ is a weak equivalence and $Uq$ is a fibration.
  \end{itemize}
  (Compare \cite[Theorem~2.2.1]{HKRS} which gives an explicit criterion
  for the ``acyclicity'' condition in the dual situation
  of transfer along a left adjoint.)
  Choose a cylinder object $1 \amalg 1 \xrightarrow{i} I \xrightarrow{p} 1$
  for the unit object $1$ of $M$.
  Then
  \[
    X \cong [F1, X]_N \xrightarrow{(Fp)^*} [FI, X]_N
    \xrightarrow{(Fi)^*} [F(1 \amalg 1), X]_N \cong X \times X
  \]
  is the required path factorization of $X$,
  by a standard adjunction argument
  using the fact that $1 \in M$ is cofibrant and $UX$ is fibrant.
\end{proof}

\begin{theorem}
  \label{q-model-structure}
  There is a model category structure on $\Shv(\PP)$
  in which a map $f : X \to Y$ is a weak equivalence (respectively, fibration)
  if and only if $\Sing f$ is a weak equivalence (respectively, Kan fibration)
  of simplicial sets.
  Furthermore, the model category $\Shv(\PP)$ is cartesian monoidal.
\end{theorem}

\begin{proof}
  We apply the above criterion to the adjunction
  $|{-}| : \sSet \rightleftarrows \Shv(\PP) : \Sing$.
  We checked in \cref{geom-real-lex} that
  geometric realization commutes with finite products
  and in \cref{sing-fibrant} that
  $\Sing X$ is fibrant for every object~$X$ of~$\Shv(\PP)$.
\end{proof}

\begin{remark}
  The above argument is essentially the same as
  Quillen's original construction of the model category structure on $\Top$,
  modulo some technicalities which do not arise in our setting
  (the $\Top$-valued geometric realization is monoidal
  only on finite simplicial sets,
  and $\Top$ is not locally presentable
  so a little more work is required to construct factorizations).
\end{remark}

\begin{nul}
  We refer to the model category structure
  constructed in \cref{q-model-structure}
  as the \emph{Kan--Quillen model structure} on $\Shv(\PP)$,
  or as the $q$-model structure for short.
  By its construction, the geometric realization--$\Sing$ adjunction
  \[
    |{-}| : \sSet \rightleftarrows \Shv(\PP) : \Sing
  \]
  is a Quillen adjunction.
  In fact, this adjunction is a Quillen equivalence
  and so $\Shv(\PP)$ models the homotopy theory of spaces.
  We will prove this at the end of this section.
\end{nul}

\begin{nul}
  \label{shv-d-model-structure}
  We can apply \cref{monoidal-transfer} again to the adjunction
  $u_! : \Shv(\PP) \rightleftarrows \Shv(\DD) : u^*$
  to produce a $q$-model structure on $\Shv(\DD)$.
  Indeed, $u_!$ preserves finite limits by \cref{geom-real-lex},
  and the condition on fibrant objects is automatic
  because every object of $\Shv(\PP)$ is fibrant.
  The description of the weak equivalences and fibrations of $\Shv(\DD)$
  is the same as that for $\Shv(\PP)$
  except that $\Sing$ is now the functor $\Sing : \Shv(\DD) \to \sSet$.

  From the perspective of this model structure,
  the category $\Shv(\DD)$ is somewhat uninteresting
  because (by the small object argument) all of its cofibrant objects belong to
  the coreflective subcategory $\Shv(\PP)$.
  In particular,
  the adjunction $u_! \dashv u^*$ is a Quillen equivalence
  simply because $\Shv(\PP)$ and $\Shv(\DD)$
  have the same cofibrant objects and the same weak equivalences between them.
  The functor $u^*$ also has a right adjoint $u_*$,
  and $u^*$ preserves the generating (acyclic) cofibrations of $\Shv(\DD)$
  because $u^* u_!$ is naturally isomorphic to the identity,
  so $u^* \dashv u_*$ is also a Quillen adjunction
  (hence a Quillen equivalence).
  By the construction of the model structure on $\Shv(\DD)$,
  the functor $u^*$ preserves and reflects weak equivalences.

  We expect that $\Shv(\DD)$ will support
  other model category structures with more cofibrant objects:
  models for either the homotopy theory of spaces,
  or some more refined homotopy theory
  as discussed in \cref{proper-homotopy-theory}.

  In the remainder of this section
  we focus mainly on the model category $\Shv(\PP)$,
  leaving the consequences for $\Shv(\DD)$ to the reader.
\end{nul}

\begin{nul}[Simplicial structure]
  The monoidal adjunction
  $|{-}| : \sSet \rightleftarrows \Shv(\PP) : \Sing$
  gives $\Shv(\PP)$ the structure of
  a category enriched, tensored and cotensored over simplicial sets,
  with formulas given respectively by
  \[
    \Map(X, Y) = \Sing (Y^X), \qquad
    K \otimes X = |K| \times X, \qquad
    Y^K = Y^{|K|},
  \]
  and the $q$-model category structure is simplicial.
\end{nul}

\begin{nul}[Cofibrations]
  The $q$-model structure has a set of generating cofibrations
  $I = \{\,|\partial \Delta^n| \to |\Delta^n| \mid n \ge 0\,\}$.
  The map $|\partial \Delta^n| \to |\Delta^n|$ is definably homeomorphic to
  the boundary inclusion $S^{n-1} \to D^n$
  where
  \[
    D^n := \{\,x \in R^n \mid ||x||_2 \le 1\,\}, \quad
    S^{n-1} := \partial D^n := \{\,x \in R^n \mid ||x||_2 = 1\,\}
  \]
  (note that $S^{-1} = \emptyset$).
  Indeed, if we replace the $\ell_2$ norm by the $\ell_\infty$ norm
  (so that the $n$-disc becomes a cube $[-1, 1]^n$)
  then the homeomorphism can be chosen to be
  a continuous piecewise-linear map with rational coefficients.

  By \cref{triangulation},
  any closed embedding $A \to B$ in $\PP$
  can be expressed as the geometric realization of
  a monomorphism between simplicial sets.
  We thus obtain the following result.
\end{nul}

\begin{proposition}
  \label{yo-mono-cof}
  If $i : A \to B$ is a monomorphism of $\PP$
  then $\yo i : \yo A \to \yo B$ is a cofibration of $\Shv(\PP)$.
  In particular, $\yo B$ is cofibrant for any $B \in \PP$.
\end{proposition}

\begin{remark}
  This observation yields
  the following characterization of the $q$-model structure
  which does not explicitly mention any specific geometric objects
  (such as spheres or simplices):
  it is the unique model category structure on $\Shv(\PP)$
  whose cofibrations are generated by
  (the images under $\yo$ of) the monomorphisms of $\PP$
  and in which every object is fibrant.

  Informally, we may think of $\Shv(\PP)$ as
  a relative of ``shape-based'' model categories
  (such as simplicial sets or cubical sets)
  in which we probe a space by its maps from \emph{all} reasonable spaces---%
  where the ``reasonable spaces'' are
  the closed and bounded sets definable in the chosen o-minimal structure.
  By restricting to sheaves for the proper topology,
  we preserve certain well-behaved colimits of $\PP$
  and ensure that every object is fibrant.

  The functor $\Sing$ is given simply by restriction
  to the (non-full) subcategory~$\Delta \subset \PP$ containing the simplices $|\Delta^n|_\PP$
  and the vertex-preserving, order-preserving affine maps between them.
  Note that $\Delta$ has no nontrivial sieves
  whose images in $\PP$ form covering families.
  Comparison functors to other categories such as cubical sets
  can be obtained by restriction to
  correspondingly chosen subcategories of $\PP$.
\end{remark}

\begin{nul}[Non-cofibrations]
  Not every object of $\Shv(\PP)$ is cofibrant.
  We give the following examples without proof;
  the theory required to verify these claims
  will be developed in \cite{paper3}.
  \begin{itemize}
  \item
    Nonseparated objects
    (i.e., non-quasiseparated in the sense of \cite{SGA4})
    cannot be cofibrant.
    For example, let $G$ be a group object of $\PP$ of positive dimension
    (such as $\mathrm{SO}(2) = S^1$)
    and let $\flat G$ denote the constant sheaf
    on the underlying discrete group.
    Then $\flat G$ acts on $\yo G$
    and the quotient $\yo G / \flat G$ is a nonseparated sheaf.
    The subobject classifier $\Omega$ is also a nonseparated sheaf.
  \item
    The object $R_\loc = u^* \yo R$ described in \cref{r_loc}
    is cofibrant if and only if the field $R$ is sequential.
    The ``if'' direction is easy:
    if $0 < r_0 < r_1 < r_2 < \cdots$ is an unbounded sequence in $R$
    then $\emptyset \to R_\loc$ is a transfinite composition of
    the cofibrations $\emptyset \to [-r_0, r_0] \to [-r_1, r_1] \to \cdots$.
    The more difficult converse is \cite[IV.9.10]{WSA},
    up to the equivalence
    between weak polytopes and cofibrant objects of $\Shv(\PP)$
    which we will explain in~\cite{paper3}.
    Informally, the problem with non-sequential $R$ is that
    if we apply a similar procedure to a longer sequence
    $0 < r_0 < r_1 < \cdots < r_{\omega} < r_{\omega + 1} < \cdots$,
    after $\omega$ steps we will have constructed $\bigcup_{i < \omega} {[-r_i, r_i]}$,
    which is generally not even a definable set, let alone an object of $\PP$.
    Then there is no reason for
    the inclusion $\bigcup_{i < \omega} {[-r_i, r_i]} \to [-r_\omega, r_\omega]$,
    which should be the next ``cell'' to attach,
    to be a cofibration.
  \end{itemize}
\end{nul}

\begin{nul}[Acyclic cofibrations and fibrations]
  The $q$-model structure has a set of generating acyclic cofibrations
  $J = \{\,|\Lambda^n_i| \to |\Delta^n| \mid n > 0, 0 \le i \le n\,\}$.
  For fixed $n$, all the geometric realizations $|\Lambda^n_i| \to |\Delta^n|$
  are definably homeomorphic to the map
  $D^{n-1} \times \{0\} \to D^{n-1} \times I$,
  where $I$ denotes the unit interval $\yo [0, 1]$,
  as well as to the ``open cylinder'' inclusion map
  \[
    w_{n-1} : (\partial D^{n-1} \times I) \cup (D^{n-1} \times \{0\})
    \to D^{n-1} \times I.
  \]
  Hence, we may equivalently take the maps
  $\{\,D^n \times \{0\} \to D^n \times I \mid n \ge 0\,\}$
  or $\{\,w_n \mid n \ge 0\,\}$ as generating acyclic cofibrations.
  The fibrations, then, are those maps $f : X \to Y$
  with the right lifting property with respect to the inclusions
  $D^n \times \{0\} \to D^n \times I$;
  these are the Serre fibrations in our context.

  Again,
  because any closed inclusion $A \subset B$ in $\PP$ can be triangulated,
  the induced map $(A \times I) \cup (B \times \{0\}) \to B \times I$
  is always an acyclic cofibration---%
  a fact which also follows from
  the pushout-product axiom for $\Shv(\PP)$.
\end{nul}

\begin{nul}[Homotopies]
  If $f_0$, $f_1 : X \to Y$ are two maps
  we define a \emph{homotopy} from $f_0$ to $f_1$
  to be a map $H : X \times I \to Y$
  with $H(-, 0) = f_0$ and $H(-, 1) = f_1$.
  This agrees with the model category-theoretic notion of \emph{right} homotopy
  because a map $H : X \times I \to Y$
  corresponds to a map $\hat H : X \to Y^I$,
  and $Y \to Y^I \to Y \times Y$ is always a path object for $Y$
  (because $* \amalg * = \{0, 1\} \to I \to *$ is a cylinder object for $*$,
  and every object $Y$ is fibrant).
  When $X$ is cofibrant,
  $X \times I$ is a cylinder object for $X$
  and $H$ is also a left homotopy.
  However, it so happens that (as in $\Top$)
  the relation of right homotopy is an equivalence relation
  regardless of whether $X$ is cofibrant,
  because homotopies can be concatenated and reversed.

  In the usual way we can also define \emph{homotopy rel $A$}
  for a map $j : A \to X$
  (normally a cofibration, in which case
  it agrees with the general model category-theoretic notion).
  Taking $A = *$, we obtain a notion of pointed homotopy
  between basepoint-preserving maps between pointed spaces,
  which is the notion of right homotopy
  in the slice model category $\Shv(\PP)_{*/}$.
  A base point $x : * \to X$ of an object $X \in \Shv(\PP)$
  can be identified with an element~$x$ of~$X(*)$.
  Note that if $X = \yo B$ for $B \in \PP$,
  then any map $x : * \to X$ is a cofibration.
\end{nul}

\begin{nul}[Weak equivalences]
  By Whitehead's theorem in a model category,
  the weak equivalences between cofibrant objects of $\Shv(\PP)$
  are precisely the homotopy equivalences
  with respect to the above notion of homotopy.

  For general objects,
  we can describe the weak equivalences in the classical way.
  We choose, for each $n \ge 0$, a point $*$ of $S^n$ to serve as basepoint.
\end{nul}

\begin{definition}
  Let $X$ be an object of $\Shv(\PP)$.
  We define $\pi_0 X$ to be the set $[*, X]$
  of homotopy classes of maps from $*$ to $X$,
  and for $x \in X(*)$ and $n \ge 0$,
  we define $\pi_n(X, x)$ to be the set $[(S^n, *), (X, x)]$
  of pointed homotopy classes of pointed maps.

  $\pi_0$ is evidently a functor,
  and $\pi_n(X, x)$ is functorial in the sense that
  a map $f : X \to Y$ induces a map $\pi_n(X, x) \to \pi_n(Y, fx)$
  for $x \in X(*)$.
\end{definition}

\begin{proposition}
  \label{weq-iff-pi}
  A map $f : X \to Y$ of $\Shv(\PP)$ is a weak equivalence
  if and only if the induced maps $\pi_0 X \xrightarrow{f_*} \pi_0 Y$
  and $\pi_n(X, x) \xrightarrow{f_*} \pi_n(Y, fx)$
  are isomorphisms for each $x \in X(*)$ and each $n \ge 0$.
\end{proposition}

\begin{proof}
  By definition, $f$ is a weak equivalence if and only if
  $\Sing f$ is a weak equivalence of simplicial sets.
  The simplicial sets $\Sing X$ and $\Sing Y$ are Kan complexes,
  so we can check whether $\Sing f$ is a weak equivalence
  using simplicial homotopy groups.
  The pointed space $(S^n, *)$ is definably homeomorphic to
  the geometric realization of
  the pointed simplicial set $(\Delta^n / \partial \Delta^n, *)$,
  so pointed maps $(S^n, *) \to (X, x)$ correspond under adjunction
  to maps $\Delta^n \to \Sing X$
  which carry the subspace $\partial \Delta^n$ to
  the $0$-simplex $\Delta^0 \xrightarrow{\Sing x} \Sing X$
  corresponding to $* = |\Delta^0| \xrightarrow{x} X$.
  Furthermore, $S^n \times I$ is the geometric realization of
  $(\Delta^n / \partial \Delta^n) \times \Delta^1$
  and so the pointed homotopy relation on pointed maps $(S^n, *) \to (X, x)$
  corresponds to simplicial homotopy rel $\partial \Delta^n$.
  We conclude that $\pi_n(X, x)$ may be identified with
  the simplicial homotopy group $\pi_n(\Sing X, \Sing x)$.
  Similar but simpler arguments apply to $\pi_0$.
  The claim then follows from the characterization of
  weak equivalences between Kan complexes by simplicial homotopy groups,
  together with the fact that the $0$-simplices of $\Sing X$
  are in bijection with $X(*)$.
\end{proof}

\begin{nul*}
  In some arguments the following criterion is more useful.
  We will see an example in the next section.
  It can also be used to prove the above statement directly
  (i.e., without recourse to simplicial sets)
  by a standard ``obstruction-theoretic'' argument.
\end{nul*}

\begin{proposition}
  \label{weq-iff-lift}
  For a map $f : X \to Y$ of $\Shv(\PP)$,
  the following conditions are equivalent:
  \begin{itemize}
  \item[(a)]
    $f$ is a weak equivalence.
  \item[(b)]
    For any monomorphism $j : K \to L$ of $\PP$
    and square
    \[
      \begin{tikzcd}
        \yo K \ar[r, "h"] \ar[d, "\yo j"'] & X \ar[d, "f"] \\
        \yo L \ar[r, "k"] & Y
      \end{tikzcd}
    \]
    there is a map $l : \yo L \to X$
    such that $lj = h$
    and $fl$ is homotopic rel $\yo K$ to $k$.
  \item[(c)]
    The same condition as (b),
    but with $j$ restricted to the monomorphisms $S^{n-1} \to D^n$ ($n \ge 0$).
  \end{itemize}
\end{proposition}

\begin{proof}
  This follows from \cite[Proposition~3.5]{Isaev},
  since every object of $\Shv(\PP)$ is fibrant
  and we can take all monomorphisms of $\PP$,
  or just the inclusions $S^{n-1} \to D^n$,
  as generating cofibrations for $\Shv(\PP)$.
\end{proof}

\begin{nul}[Properness]
  The model category $\Shv(\PP)$ is right proper
  because all of its objects are fibrant.
  We will show in the next section
  that it is left proper as well (\cref{left-proper}).
\end{nul}

\begin{nul}[Functoriality]
  Suppose $\Phi : (R, \sS) \to (R', \sS')$ is a morphism of o-minimal structures.
  Then there is a diagram of categories
  \[
    \begin{tikzcd}
      \Delta \ar[r] \ar[rd] &
      \PP(R, \sS) \ar[r, "u"] \ar[d, "\Phi"] &
      \DD(R, \sS) \ar[d, "\Phi"] \\
      & \PP(R', \sS') \ar[r, "u"] &
      \DD(R', \sS')
    \end{tikzcd}
  \]
  which commutes because
  $\Phi$ preserves the structure $0$, $1$, $\le$
  used to define the o-minimal simplex functor $|\Delta^\bullet|_\PP$.
  The diagram above induces a diagram of left adjoints
  \[
    \begin{tikzcd}
      \sSet \ar[r, "|{-}|"] \ar[rd, "|{-}|"'] &
      \Shv(\PP(R, \sS)) \ar[r, "u_!"] \ar[d, "\Phi_!"] &
      \Shv(\DD(R, \sS)) \ar[d, "\Phi_!"] \\
      & \Shv(\PP(R', \sS')) \ar[r, "u_!"] &
      \Shv(\DD(R', \sS'))
    \end{tikzcd}
  \]
  in which the vertical functors $\Phi_!$ are also left Quillen functors;
  this follows from the universal property of transferred model structures
  (\cref{transferred-features}).
\end{nul}

\begin{nul}[Topological realization]
  Suppose $R = \RR$.
  Then the topological realization functor
  $|{-}|_\tp : \Shv(\DD) \to \Top$
  described in \cref{topological-realization}
  is a left Quillen functor.
  This follows from the fact that the composition
  \[
    \sSet \xrightarrow{|{-}|} \Shv(\DD)
    \xrightarrow{|{-}|_\tp} \Top
  \]
  is the classical geometric realization,
  a left Quillen functor.
\end{nul}

\begin{nul}[Universal property]
  As described in \cref{cc-of-model},
  we may regard the model category $\Shv(\PP)$
  also as a cofibration category.
  Recall that $\PP$ has the structure of a cofibration category
  in which the cofibrations are the closed embeddings
  and the weak equivalences are the homotopy equivalences
  defined using the interval object $[0, 1]$.
  The functor $\yo : \PP \to \Shv(\PP)$
  preserves the initial object and pushouts of closed embeddings
  by \cref{proper-sheaves-squares}.
  Two parallel maps $f_0$, $f_1 : A \to B$ in $\PP$
  are homotopic if and only if their images under $\yo$ are.
  Therefore $\yo$ takes the weak equivalences of $\PP$
  to homotopy equivalences and, in particular, weak equivalences
  of $\Shv(\PP)$.
  Finally, $\yo$ preserves cofibrations by \cref{yo-mono-cof},
  so we conclude that $\yo : \PP \to \Shv(\PP)$ is an exact functor.

  We now claim that $\yo$ exhibits $\Shv(\PP)$
  as a ``model category cocompletion'' of the cofibration category $\PP$
  in the following sense.
\end{nul}

\begin{proposition}
  \label{model-category-universal}
  The functor $\yo : \PP \to \Shv(\PP)$
  is the universal exact functor from $\PP$ to a model category.
  More precisely, given any model category $M$,
  composition with $\yo : \PP \to \Shv(\PP)$
  induces an equivalence from
  the category of all left Quillen functors from $\Shv(\PP)$ to $M$
  to the category of all exact functors from $\PP$ to
  the underlying cofibration category of $M$.
\end{proposition}

\begin{proof}
  Let $S$ denote the following set of colimit cocones of $\PP$:
  \begin{itemize}
  \item
    the initial object $\emptyset$,
    as a colimit of the empty diagram;
  \item
    closed-by-proper pushouts, or equivalently,
    pushouts of cofibrations.
  \end{itemize}
  By definition a functor $F : \PP \to M$ to a cofibration category $M$ is exact
  if and only if it satisfies the conditions:
  \begin{enumerate}
  \item
    $F$ preserves the colimit cocones of $S$.
  \item
    $F$ preserves cofibrations and acyclic cofibrations.
  \end{enumerate}
  By \cref{proper-sheaves-squares},
  $\Shv(\PP)$ is the full subcategory of $\PSh(\PP)$
  consisting of the presheaves which send the elements of $S$ to limit cones.
  Then, by the theory of sketches,
  $\yo : \PP \to \Shv(\PP)$ is the free cocompletion
  preserving the colimit cocones $S$.
  In other words, for any cocomplete $M$,
  giving a functor $F : \PP \to M$ satisfying condition (\textit{i})
  is the same as giving a colimit-preserving functor $\hat F : \Shv(\PP) \to M$,
  where $F$ is related to $\hat F$ by $F = \hat F \circ \yo$.
  (See \cite[Theorem~11.5]{Fiore96enrichmentand} for a precise statement.)

  Now, let $M$ be a model category.
  It remains to check that, in this correspondence,
  the functor $F$ satisfies condition (\textit{ii})
  if and only if the associated colimit-preserving functor $\hat F$
  preserves all cofibrations and acyclic cofibrations of $\Shv(\PP)$.
  This is clear because
  the images of the (acyclic) cofibrations under $\yo : \PP \to \Shv(\PP)$
  are generating (acyclic) cofibrations for $\Shv(\PP)$.
\end{proof}

\begin{nul*}
  Finally, as promised earlier,
  we prove that the geometric realization--$\Sing$ adjunction
  is a Quillen equivalence.
  Our proof relies on
  the identification of the homotopy theory of $\PP \subset \Shv(\PP)$
  with the homotopy theory of finite CW complexes (\cref{ho-p}),
  and therefore on
  the model-theoretic methods of Delfs and Knebusch \cite{LSA}
  as well as the ``normal triangulations'' of Baro and Otero \cite{BO10}.
\end{nul*}

\begin{theorem}
  \label{quillen-eq}
  The Quillen adjunction
  $|{-}| : \sSet \rightleftarrows \Shv(\PP) : \Sing$
  is a Quillen equivalence.
\end{theorem}

\begin{proof}
  As every object of $\sSet$ is cofibrant
  and every object of $\Shv(\PP)$ is fibrant,
  we must check that:
  \begin{enumerate}
  \item
    for each $K \in \sSet$,
    the unit map $\eta_K : K \to \Sing |K|$ is a weak equivalence
    of simplicial sets;
  \item
    for each $X \in \Shv(\PP)$,
    the counit map $\varepsilon_X : \mathopen|\Sing X\mathclose| \to X$ is a weak equivalence
    in $\Shv(\PP)$.
  \end{enumerate}
  In (\textit{ii}), by definition, $\varepsilon_X$ is a weak equivalence
  if and only if $\Sing(\varepsilon_X) : \Sing {\mathopen|\Sing X\mathclose|} \to \Sing X$
  is a weak equivalence in $\sSet$.
  This map is left inverse to
  $\eta_{\Sing X} : \Sing X \to \Sing {\mathopen|\Sing X\mathclose|}$,
  so it suffices to verify (\textit{i}).

  Let $K$ be a simplicial set.
  We may express $K$ as the filtered colimit
  of its finite simplicial subsets.
  The functor $\Sing$ commutes with filtered colimits,
  because $(\Sing X)_n = \Hom_{\Shv(\PP)}(\yo |\Delta^n|, X)$
  and the proper topology is finitary.
  Furthermore, a filtered colimit of weak equivalences in $\sSet$
  is again a weak equivalence \kerodon{00XE}.
  Therefore it suffices to verify the claim
  for finite simplicial sets.
  This we did in \cref{eta-weq-of-finite}.
\end{proof}

\begin{remark}
  In subsequent work, we will show that
  the full subcategory of $\Shv(\PP)$
  consisting of the cofibrant objects
  is equivalent to the category of \emph{weak polytopes}
  of \cite{WSA} (in the semialgebraic case)
  and \cite{Pi08} (in the case of a general o-minimal structure).
  Using this identification,
  we could also deduce \cref{quillen-eq}
  from the ``two main theorems on homotopy sets'' of \cite{WSA,Pi08}.
  Those theorems ultimately rely on the same inputs
  used to prove \cref{ho-p}:
  normal triangulations to reduce to the semialgebraic case
  and then model-theoretic arguments to reduce to the case $R = \RR$.

  When $R$ is archimedean,
  \cref{quillen-eq} can be proven
  by techniques of classical algebraic topology
  (for example, by imitating the proof in \kerodon{012Y}).
  The difficulty in the nonarchimedean case
  is that an open cover of a simplex
  might not be refined by $N$-fold barycentric subdivision for any $N$,
  so arguments based on the $\Ex^\infty$ functor
  or ``small simplices'' arguments in homology
  become unavailable.
  It is however possible to prove \cref{quillen-eq}
  (or more precisely \cref{ho-p}, on which it relies)
  using a variation of the o-minimal triangulation theorem
  which gives better control over the homotopy type
  of the produced simplicial complex;
  this will appear in forthcoming work.
\end{remark}

%% file: tex/justleftproper.tex
\section{Left properness}
\label{left-properness}

\begin{nul*}
  The aim of this section is to prove
  the last remaining claim about the model category $\Shv(\PP)$,
  namely that it is left proper.
  This means that for every cofibration $j : A \to A'$
  and weak equivalence $w : A \to B$,
  in the pushout square
  \[
    \begin{tikzcd}
      A \ar[r, "w"] \ar[d, "j"'] & B \ar[d, "j'"] \\
      A' \ar[r, "w'"] & B'
    \end{tikzcd}
  \]
  the induced map $w'$ is also a weak equivalence.
\end{nul*}

\begin{nul*}
  A typical approach to proving that the model category $\Top$ is left proper
  is as follows.
  First, reduce by general arguments to the case that
  $j$ attaches a single cell.
  Then use the van~Kampen and Mayer--Vietoris theorems
  to prove that $w'$ induces a weak equivalence on fundamental groupoids
  and on singular homology with local coefficients.
  Finally, use the Hurewicz theorem to deduce that
  $w'$ is a weak homotopy equivalence.
  In $\Shv(\PP)$, we do not yet have these tools available.
  Furthermore, it is not clear how one could reduce
  the left properness of $\Shv(\PP)$ to that of $\Top$,
  especially as we only have a direct comparison functor $\Shv(\PP) \to \Top$
  when the field $R$ is a subfield of the real numbers.

  There is however a more direct approach to left properness of $\Top$
  using the following result on gluing weak equivalences,
  apparently due to Gray
  \cite[Proposition 16.24]{Gray}\rlap{.}\footnote{
    We give a simplified statement which is enough for our purposes.
    In the original statement,
    the cover $Y = Y_1 \cup Y_2$ only needs to be excisive,
    meaning the interiors of $Y_1$ and $Y_2$ cover $Y$,
    and the cover $X = X_1 \cup X_2$ need not be
    the one pulled back along $f$,
    merely an excisive cover of $X$ satisfying $X_i \subset f^{-1}(Y_i)$.}
  Let $f : X \to Y$ be a map of topological spaces,
  let $Y = Y_1 \cup Y_2$ be an open cover of $Y$,
  and set $X_1 = f^{-1}(Y_1)$ and $X_2 = f^{-1}(Y_2)$.
  Then if the restrictions $f_1 : X_1 \to Y_1$, $f_2 : X_2 \to Y_2$
  and $f_{12} : X_1 \cap X_2 \to Y_1 \cap Y_2$
  are all weak homotopy equivalences,
  then so is $f$.

  To deduce that $\Top$ is left proper,
  we again reduce to the case that $j : A \to A \cup e^n = A'$
  attaches a single cell.
  Pick an interior point $p$ in $e^n$,
  and apply the gluing theorem to $w' : A' \to B'$
  and the cover of $B'$ by $B' - \{p\}$ and the interior of $e^n$.
  The restriction of $w'$ to a map $A' - \{p\} \to B' - \{p\}$
  is homotopy equivalent to the original map $w$,
  hence a weak equivalence,
  and the other two restrictions of $w'$ involved
  are homeomorphisms.
\end{nul*}

\begin{nul*}
  The proof of Gray's gluing theorem goes roughly as follows.
  To say that $f : X \to Y$ is a weak homotopy equivalence
  means that in any square
  \[
    \begin{tikzcd}
      \partial I^n \ar[r, "h"] \ar[d] & X \ar[d, "f"] \\
      I^n \ar[r, "k"] & Y
    \end{tikzcd}
  \]
  there is a lift $l : I^n \to X$
  making the top triangle commute strictly
  and the bottom triangle commute up to homotopy rel $\partial I^n$.
  Pull back the open cover $Y = Y_1 \cup Y_2$ along $k$
  to obtain an open cover of $I^n$.
  By a Lebesgue number argument,
  we can subdivide $I^n$ into small cubes
  finely enough that each small cube is mapped into
  either $Y_1$ or $Y_2$.
  Then we can realize $I^n$ as a CW complex
  which is the union of two subcomplexes $V_1$ and $V_2$,
  with the property that
  $k$ maps $V_1$ into $Y_1$ and $V_2$ into $Y_2$.
  Now using the hypotheses on $f_{12}$, $f_1$, and $f_2$,
  one constructs the desired lift first on $V_1 \cap V_2$,
  then extends it separately to both $V_1$ and $V_2$.

  In our setting, a Lebesgue number argument is not available
  since the field $R$ may be nonarchimedean.
  However, the triangulation theorem
  (or more specifically, its corollary that
  every monomorphism of $\PP$ has the homotopy extension property)
  is a more than adequate replacement.
  In particular, the analogue of an excisive cover in $\Shv(\PP)$
  is just the built-in notion of a cover that it has as a topos.
  We postpone a full treatment of Gray-style gluing theorems
  to future work,
  and here present just enough to prove that
  $\Shv(\PP)$ is left proper.
\end{nul*}

\begin{remark}
  Again, our arguments will depend only on
  a few formal properties of the category $\PP$
  together with its ``interval'' object $[0, 1]$.
  However, in the absence of other known suitable examples,
  it seems premature to fix any specific definition
  of ``good category with interval object''
  to which our arguments would apply.
\end{remark}

\begin{convention}
  By now the explicit notation $\yo$ for the Yoneda embedding of $\PP$ in $\Shv(\PP)$
  has fulfilled any disambiguatory purpose it once had.
  Henceforth, we simply regard $\PP$ as a full subcategory of $\Shv(\PP)$.
  This subcategory is closed under finite limits, finite coproducts
  and quotients by equivalence relations,
  including pushouts by monomorphisms.
  It is not closed under ``bad'' pushouts,
  such as those of \cref{bad-colimits}.
\end{convention}

\begin{lemma}
  Let $Y \in \Shv(\PP)$ be covered by subobjects $Y_1 \subset Y$ and $Y_2 \subset Y$,
  let $L \in \PP$ and let $k : L \to Y$ be a map.
  Then there is a cover of $L$ by closed subsets $L_1 \subset L$ and $L_2 \subset L$
  such that $k$ maps $L_1$ into $Y_1$ and $L_2$ into $Y_2$.
\end{lemma}

\begin{proof}
  Pulling back the cover $\{Y_1 \to Y, Y_2 \to Y\}$ of $Y$ along $k : L \to Y$
  produces a cover of $L$,
  which can be refined by a generating covering family
  $\{A_j \xrightarrow{p_j} L\}_{j \in J}$
  (with $J$ a finite set, $A_j \in \PP$, and the $p_j$ jointly surjective).
  Then $kp_j : A_j \to Y$ factors through either $Y_1$ or $Y_2$ for each $j$.
  By forming coproducts,
  we may assume that this family consists of
  just two maps $\{A_1 \xrightarrow{p_1} L, A_2 \xrightarrow{p_2} L\}$
  such that $kp_1$ factors through $Y_1$ and $kp_2$ factors through $Y_2$.
  However, the maps $p_1$ and $p_2$ need not be injective.

  For $i \in \{1, 2\}$, let $L_i$ denote the image of $p_i : A_i \to L$;
  then $L_i \subset L$ is closed, and $L = L_1 \cup L_2$.
  The map $p_i : A_i \to L_i$ is surjective in $\Shv(\PP)$
  because the family $\{p_i\}$ is a generating covering family for the proper topology.
  Then, the fact that $kp_i : A_i \to Y$ factors through the subobject $Y_i$ of $Y$
  means that the composition $L_i \subset L \xrightarrow{k} Y$
  also factors through $Y_i$.
\end{proof}

\begin{lemma}
  \label{mono-mono-weq-pushout}
  Let
  \[
    \begin{tikzcd}
      X_0 \ar[r, "i"] \ar[d, "w_0"'] & X \ar[d, "w"] \\
      Y_0 \ar[r] & Y
    \end{tikzcd}
  \]
  be a pushout square in $\Shv(\PP)$
  in which $w_0$ and $i$ are both monomorphisms.
  If $w_0$ is a weak equivalence, then so is $w$.
\end{lemma}

\begin{proof}
  By \cref{weq-iff-lift}, it suffices to show that
  for any monomorphism $j : K \to L$ of $\PP$
  and any square
  \[
    \begin{tikzcd}
      K \ar[r, "h"] \ar[d, "j"'] & X \ar[d, "w"] \\
      L \ar[r, "k"] & Y
    \end{tikzcd}
  \]
  there is a lift $l : L \to X$ making the top triangle commute strictly
  and the bottom triangle commute up to homotopy rel $K$.
  We apply the preceding lemma to the map $k : L \to Y$
  and the cover $Y = Y_0 \cup X$,
  obtaining a cover of $L$ by closed subsets $L_1$ and $L_2$
  such that $k$ maps $L_1$ into $Y_0$ and $L_2$ into $X$.

  We claim that without loss of generality,
  we may assume that $L_2 = K$.
  First, since $K \subset L$ is a closed subset
  that $k$ maps into $X$,
  we may enlarge $L_2$ to $L_2 \cup K$.
  So, assume that $K \subset L_2$.
  Then, since $L_2 \subset L \xrightarrow{k} Y$
  factors through $X \subset Y$,
  we may also enlarge $K$ to $L_2$
  and replace $h : K \to X$ by this factorization $h' : L_2 \to X$.
  A solution $l : L \to X$ to the new lifting problem
  is also a solution to the original one:
  the homotopy $H : wl \simeq k$ rel $L_2$
  is also a homotopy rel $K$.
  Hence, we may assume $L_2 = K$.

  Now we rename $L_1$ to $L_0$ and set $K_0 = K \cap L_0$,
  obtaining a cube
  \[
    \begin{tikzcd}[row sep=small, column sep=small]
      K_0 \ar[rrrr, "h_0", pos=0.6] \ar[dd, "j_0"'] \ar[rd] & & & & X_0 \ar[dd, "w_0", pos=0.8] \ar[rd] \\
      & K \ar[rrrr, crossing over, "h", pos=0.4] & & & & X \ar[dd, "w"] \\
      L_0 \ar[rrrr, "k_0", pos=0.6] \ar[rd] \ar[rrrruu, loosely dashed, bend left = 20, "l_0"] & & & & Y_0 \ar[rd] \\
      & L \ar[rrrr, "k", pos=0.4] & & & & Y
      \ar[from=2-2, to=4-2, crossing over, "j"', pos=0.2]
    \end{tikzcd}
  \]
  in which each of the left and right squares is both a pushout and a pullback.
  By hypothesis $w_0$ is a weak equivalence,
  so in the back face we can find a lift $l_0 : L_0 \to X_0$ with $l_0 j_0 = h_0$
  and a homotopy $H_0 : w_0 l_0 \simeq k_0$ rel $K_0$.
  We define $l : L \to X$ to be $h$ on $K$ and $l_0$ on $L_0$,
  and $H : L \times [0, 1] \cong L_0 \times [0, 1] \cup K \times [0, 1] \to Y$
  to be $H_0$ on $L_0 \times [0, 1]$ and
  the constant homotopy at $kj : K \to Y$ on $K \times [0, 1]$.
  This makes sense because $H_0$ was a homotopy rel $K_0$,
  so equals the constant homotopy at $k_0 j_0$ on $K_0 \times I$.
  Then $H$ is a homotopy rel $K$ from $wl$ to $k$.
\end{proof}

\begin{nul*}
  The rest of the proof of left properness follows classical lines.
  We leave some of the detailed verifications to the reader.
\end{nul*}

\begin{nul}
  Let $f : A \to B$ be a map of $\Shv(\PP)$.
  The \emph{mapping cylinder} of $f$ is the object $M_f$
  constructed as the pushout below.
  (Since the interval $[0, 1]$ has an automorphism
  swapping its endpoints $0$ and $1$,
  the orientation of the mapping cylinder is not important.)
  \[
    \begin{tikzcd}
      A \ar[r, "f"] \ar[d, "\id \times 1"'] & B \ar[d] \\
      A \times [0, 1] \ar[r] & M_f
    \end{tikzcd}
  \]
  There is a map $g : M_f \to B \times [0, 1]$
  given by $f \times \id$ on $A \times [0, 1]$ and by $\id \times 1$ on $B$.
  The map $f$ has the \emph{homotopy extension property}
  if this map $g$ admits a retraction.
  In this case $g$ must be a monomorphism
  and so the original map $f$ must be one as well,
  since there are pullback squares as shown below.
  \[
    \begin{tikzcd}
      A \ar[r, "f"] \ar[d] & B \ar[r] \ar[d] & \{0\} \ar[d] \\
      M_f \ar[r, "g"] & B \times [0, 1] \ar[r, "\pi_2"] & {[0, 1]}
    \end{tikzcd}
  \]
  For any $B \in \Shv(\PP)$,
  the unique map $f : \emptyset \to B$ has the homotopy extension property,
  since then $g : M_f = B \times \{1\} \to B \times [0, 1]$ is the inclusion.
\end{nul}

\begin{nul}
  By an adjointness argument,
  a map $f : A \to B$ has the homotopy extension property
  if and only if it has the left lifting property
  with respect to $\ev_1 : Z^{[0, 1]} \to Z$ for every $Z \in \Shv(\PP)$.
  Then, since $\Shv(\PP)$ is a cartesian monoidal model category,
  we deduce that every cofibration has the homotopy extension property.
\end{nul}

\begin{lemma}
  $\Shv(\PP)$ admits the structure of a cofibration category
  in which the cofibrations are the maps with the homotopy extension property
  and the weak equivalences are the homotopy equivalences.
\end{lemma}

\begin{proof}
  We apply \cref{cc-of-icat} with the cylinder functor $IX = X \times [0, 1]$
  and the class of all maps with the homotopy extension property
  as the cofibrations.
  Since these maps are determined by a left lifting property,
  the only nontrivial condition is (I4):
  if $f : A \to B$ has the homotopy extension property
  then so does the induced ``relative cylinder'' map
  \[
    f' : A \times [0, 1] \amalg_{A \times \{0, 1\}} B \times \{0, 1\} \to B \times [0, 1].
  \]
  By an adjointness argument, $f'$ having the homotopy extension property
  amounts to saying that
  $f$ has the left lifting property with respect to the restriction map
  $Z^{[0, 1]^2} \to Z^K$ for every $Z \in \Shv(\PP)$,
  where $K \subset [0, 1]^2$ denotes the subset
  $K = [0, 1] \times \{0, 1\} \cup \{1\} \times [0, 1]$.
  The inclusion $K \subset [0, 1]^2$ is isomorphic in $\PP$
  to the inclusion $\{1\} \times [0, 1] \subset [0, 1]^2$
  by a direct geometric construction,
  and so the restriction map $Z^{[0, 1]^2} \to Z^K$
  is isomorphic to $\ev_1 : (Z^{[0, 1]})^{[0, 1]} \to (Z^{[0, 1]})$.
  Hence the claim follows by applying the homotopy extension property of $f$
  to the object $Z^{[0, 1]}$.
  (See \cite{SV02} for a more systematic treatment of this kind of argument.)
\end{proof}

\begin{nul}
  In particular, we conclude:
  \begin{enumerate}
  \item
    Any map $f : A \to B$ admits a factorization into
    a map with the homotopy extension property
    followed by a homotopy equivalence.
    (In fact, the usual mapping cylinder factorization of $f$
    is such a factorization,
    although we will not make use of this.)
  \item
    A pushout of a homotopy equivalence
    by a map with the homotopy extension property
    is a homotopy equivalence.
    This is because every object is cofibrant
    in the cofibration category structure we just constructed.
  \end{enumerate}

  We do not know whether this cofibration category structure
  extends to a ``Str\o m-style'' model category structure on $\Shv(\PP)$.
  We continue to use the terms cofibration and weak equivalence
  to refer to the ``Serre--Quillen-style'' model category structure
  constructed in the previous section.

  Now, we can finish the proof that $\Shv(\PP)$ is left proper.
\end{nul}

\begin{proposition}
  \label{left-proper}
  The model structure of \cref{q-model-structure} is left proper.
\end{proposition}

\begin{proof}
  Let $j : A \to A'$ be a cofibration and $w : A \to B$ a weak equivalence,
  and form the pushout square below.
  We must show that $w' : A' \to B'$ is a weak equivalence.
  \[
    \begin{tikzcd}
      A \ar[r, "w"] \ar[d, "j"'] & B \ar[d] \\
      A' \ar[r, "w'"] & B'
    \end{tikzcd}
  \]
  Factor $w$ into a map $i : A \to C$ with the homotopy extension property
  followed by a homotopy equivalence $q : C \to B$,
  and form two pushout squares as shown below.
  \[
    \begin{tikzcd}
      A \ar[r, "i"] \ar[d, "j"'] & C \ar[r, "q"] \ar[d] & B \ar[d] \\
      A' \ar[r, "i'"] & C' \ar[r, "q'"] & B'
    \end{tikzcd}
  \]
  Since $j$ is a cofibration it has the homotopy extension property,
  and so both $i$ and $j$ are monomorphisms.
  Since $q$ is a homotopy equivalence it is also a weak equivalence,
  so by two-out-of-three, $i$ is also a weak equivalence.
  Now $i'$ is a weak equivalence by \cref{mono-mono-weq-pushout},
  and $q'$ is a pushout of a homotopy equivalence
  by (a pushout of) a map with the homotopy extension property,
  hence a homotopy equivalence and so also a weak equivalence.
\end{proof}

\begin{nul}
  Using \cref{shv-d-model-structure},
  we deduce that $\Shv(\DD)$ is also left proper, as follows.
  Given a pushout of a weak equivalence by a cofibration in $\Shv(\DD)$,
  apply $u^*$ to it;
  the result is a pushout of a weak equivalence by a cofibration in $\Shv(\PP)$,
  hence a weak equivalence.
  Because $u^*$ reflects weak equivalences,
  we deduce that the original map in $\Shv(\DD)$
  was also a weak equivalence.
\end{nul}

%% file: tex/related.tex
\section{Related work and questions}
\label{related}

\begin{nul*}
  The semialgebraic sets,
  or more generally, the sets defined in a fixed o-minimal structure
  form a ``category of spaces'' which is strikingly pleasant,
  but limited to spaces with the homotopy type of a finite CW complex.
  The idea of enlarging the class of spaces
  by gluing definable sets together in certain ways
  is certainly not a new one.
  Here, we describe some existing work in this direction
  and the relations to the theory presented here.
\end{nul*}

\begin{nul}
  The first such notion we are aware of
  is the \emph{locally semialgebraic spaces}
  of Delfs and Knebusch \cite{LSA}.
  Modeled on the theory of schemes in algebraic geometry,
  a locally semialgebraic space consists of
  a kind of ringed ``generalized topological space''
  which admits an open cover by affine semialgebraic spaces.
  Extensions to the o-minimal setting
  can be found in \cite{BO10b} and \cite{Pi08}.

  This greatly expands the class of possible spaces,
  for example by allowing spaces with infinitely many path components,
  but is not particularly convenient for homotopy theory
  because the possible kinds of gluing are too restricted.
  In particular, every point of a locally semialgebraic space
  must have an open neighborhood which looks like
  a semialgebraic subset of $R^n$ for some $n$.
  Consequently,
  many infinite CW complexes with tame attaching maps
  (such as the wedge of countably many circles at a single point,
  or infinite-dimensional projective spaces)
  cannot be realized as locally semialgebraic spaces
  (although one could realize them up to homotopy equivalence
  by something like a mapping telescope construction).
\end{nul}

\begin{nul}
  In order to include such infinite CW complexes,
  Knebusch introduced \emph{weakly semialgebraic spaces} \cite{WSA}.
  A weakly semialgebraic space is
  a kind of generalized topological space equipped with a sheaf of rings
  which admits an \emph{exhaustion}:
  a carefully controlled directed system of closed semialgebraic subsets.
  When these closed semialgebraic subsets can be chosen to be polytopes,
  the weakly semialgebraic space is called a \emph{weak polytope}.
  This theory was extended to a general o-minimal structure
  by Pi\k{e}kosz \cite{Pi08}.

  Weakly semialgebraic spaces and especially weak polytopes
  are significantly more convenient than locally semialgebraic spaces
  for the purposes of homotopy theory.
  The main point is that one glues along closed subsets
  rather than open ones;
  for instance, the category of weak polytopes
  admits sequential colimits of closed embeddings
  and pushouts of closed embeddings along arbitrary maps.
  In addition, closed embeddings between weak polytopes
  have the homotopy extension property.
  Knebusch proved that
  the homotopy category of weak polytopes
  is equivalent to the classical homotopy category of CW complexes.

  In subsequent work \cite{paper3}, we will show how
  the success of the theory of weak polytopes
  is ``explained'' by the existence of the model category $\Shv(\PP)$:
  \begin{itemize}
  \item
    The cofibrant objects of $\Shv(\PP)$
    form a category equivalent to the category of weak polytopes.
  \item
    The cofibrations and the weak equivalences between cofibrant objects of $\Shv(\PP)$
    correspond to the closed embeddings and the homotopy equivalences between weak polytopes,
    respectively.
  \end{itemize}
  It is remarkable that Knebusch effectively defined
  the cofibrant objects of $\Shv(\PP)$
  without having the full model category at hand
  and by rather different techniques.

  The advantage of the full model category $\Shv(\PP)$
  is that it contains (for example) loop spaces
  with the correct ``point-set level'' universal property,
  which are usually not cofibrant.
  This limitation of weak polytopes was noted in the introduction of \cite{WSA},
  which works around this deficiency to some extent
  by using what are effectively cofibrant replacements
  for the missing loop spaces.
\end{nul}

\begin{nul}
  Marra and Menni recently introduced the \emph{PL topos} \cite{MarraMenni}.
  This is the category of sheaves on the site $\mathsf{P}$
  whose objects are polyhedra and morphisms are piecewise-linear maps
  and whose covering families are generated by
  finite covers by closed subsets.
  This site would be the same as our $\PP(R, \sS)$ with the closed topology
  for $\sS$ the semilinear o-minimal structure on an ordered field $R$
  (\cref{omin-examples}).
  However, multiplication as a function from $R \times R$ to $R$
  is not definable in this structure
  and so it does not fit into our setup.
  In particular, the category $\mathsf{P}$
  does not admit pushouts by closed inclusions
  (\cref{multiplication-required})
  and so we cannot expect to have an analogue of \cref{proper-sheaves-squares}.

  While working in an o-minimal structure with definable multiplication
  confers several convenient formal properties,
  it would also be of interest to formulate a version of the theory
  that applies to the PL topos.
  For example, one could try to make the PL topos into a model category
  in a way which is intermediate between simplicial sets
  and the model categories we have introduced here.
  Such a model structure should also be directly related to
  the classical homotopy theory of simplicial complexes.
\end{nul}

\begin{nul}
  Our work is obviously influenced by
  the methods of motivic homotopy theory generally
  and the cdh-topology in particular,
  and the model category $\Shv(\PP)$ could be seen as
  a model for ``unstable o-minimal motivic homotopy theory''.
  The biggest difference relative to the algebric setting
  is the fact that we have enough definable maps
  (for example, retractions $|\Delta^n| \to |\Lambda^n_i|$)
  to get a good homotopy theory
  without relying on an external model for spaces such as simplicial sets.

  In contrast to the case of algebraic geometry
  (but like the case of differential geometry, see \cite[Proposition~8.3]{Dugger})
  the ``motivic homotopy theory'' that we obtain
  is just the ordinary homotopy theory of spaces.
  This is essentially because
  any polytope has the definable homotopy type (even homeomorphism type)
  of a simplicial complex, by the triangulation theorem.
  A ``linear'' version of this was noted in \cite{Huber_semialg_mot}.

  We expect there to also be a direct relationship to motivic homotopy theory
  in the form of ``realization'' or comparison functors
  defined on the motivic homotopy theory of varieties over a field of the form $R[i]$
  for a real closed field $R$.
\end{nul}

\begin{nul}
  We end with some open questions.
  \begin{itemize}
  \item
    We have constructed a $q$-model structure on $\Shv(\PP)$.
    Is there also an ``$h$-model structure''
    (corresponding to the Str\o m model structure on $\Top$)
    in which the weak equivalences are the homotopy equivalences
    and the cofibrations are
    the closed embeddings with the homotopy extension property?
    If so, the associated mixed model structure \cite{Cole} on $\Shv(\PP)$
    would be another model category for spaces with every object fibrant,
    but many more cofibrant objects.
  \item
    Suppose $X \in \PP$ and $Y$ is a cofibrant object of $\Shv(\PP)$.
    Is the exponential $Y^X$ homotopy equivalent to a cofibrant object?
    For $\Top$ this is a classical theorem of Milnor \cite{Milnor59}.
  \item
    Suppose we work in a constructive background theory,
    and interpret ``o-minimal structure'' in a suitable sense
    as described in \cref{constructive}.
    How much of the theory developed here
    can be proved in this setting?
  \item
    Can we construct the $q$-model structure
    in a ``self-contained'' fashion,
    without reference to simplicial sets?
    This is possibly helpful for a constructive treatment
    but may also be of independent interest.
  \item
    The category $\PP$ is a pretopos with an interval object $[0, 1]$
    with respect to which every monomorphism has the homotopy extension property.
    Furthermore, the topology $\tauad$ on $\PP$ described in \cref{topologies-agree}
    agrees with the coherent topology.
    As discussed in \cref{pretopoi-for-spaces},
    most of the theory of $\Shv(\PP)$ depends only on these properties.
    Are there any other nontrivial examples of such categories?
  \end{itemize}
\end{nul}

%% file: tex/cofibration.tex
\section{Cofibration categories}
\label{cofibration-categories}

\begin{nul*}
  Model categories are a framework
  for ``large'' homotopy theories that admit all limits and colimits,
  like the homotopy theory of all spaces.
  However, we are also interested in ``small'' homotopy theories
  like the homotopy theory of \emph{finite} spaces
  (finite simplicial sets, or finite CW complexes).
  Cofibration categories model finitely cocomplete homotopy theories
  and functors between them that preserve finite colimits.

  In this appendix we collect facts about cofibration categories
  which will be used in the main text.
  Our main aim is \cref{cc-unit-weq} which gives conditions
  under which the unit of a relative adjunction between cofibration categories
  is a weak equivalence.
\end{nul*}

\begin{definition}
  \label{cofibration-category}
  A \emph{cofibration category} is a category $C$
  equipped with two classes of maps
  called \emph{cofibrations} and \emph{weak equivalences}
  satisfying the axioms below.
  \begin{enumerate}
  \item
    $C$ has an initial object $\emptyset$,
    and $\id_\emptyset$ is a cofibration.
  \end{enumerate}
  We call an object $A$ \emph{cofibrant}
  if the unique map $\emptyset \to A$ is a cofibration.
  \begin{enumerate}[resume]
  \item
    The weak equivalences of $C$ satisfy the two-out-of-three axiom,
    and isomorphisms are weak equivalences.
  \item
    \label{pushout-cofibration}
    The cofibrations of $C$ are closed under composition.
    The pushout of a cofibration $j : A \to B$ with cofibrant domain
    along a map $f : A \to A'$ to another cofibrant object exists,
    and any such pushout is again a cofibration,
    which is also a weak equivalence if $j$ is one.
  \item
    Every map with cofibrant domain
    factors as a cofibration followed by a weak equivalence.
  \item
    $C$ satisfies the following conditions
    which are equivalent given the preceding axioms
    \cite[Theorem~7.2.7]{RB}:
    \begin{itemize}
    \item
      A retract of a weak equivalence is a weak equivalence.
    \item
      The weak equivalences satisfy the two-out-of-six axiom:
      if $gf$ and $hg$ are weak equivalences then so is $g$
      (hence also $f$ and $h$).
    \end{itemize}
  \end{enumerate}
  Condition (\textit{iii}) implies that
  any isomorphism between cofibrant objects is a cofibration
  (because it may be written as a pushout of $\id_\emptyset$).

  We call a map of a cofibration category an \emph{acyclic cofibration}
  if it is both a cofibration and a weak equivalence.
  If $C$ and $D$ are cofibration categories,
  then a functor $F : C \to D$ is \emph{exact}
  if $F$ preserves cofibrations, acyclic cofibrations,
  the initial object, and pushouts of cofibrations.
\end{definition}

\begin{nul}
  There are a number of different notions of cofibration category
  in the literature.
  The one above is from \cite{Cis10}
  (where it is called a ``categorie d\'erivable \`a droite'')
  and \cite{RB}
  (where it is called a ``precofibration category''),
  except that we have built the saturation condition (\textit{v})
  into the definition of a cofibration category:
  this is the condition that ensures that
  the weak equivalences of $C$ are precisely
  the maps that become isomorphisms in the homotopy category.
  Cofibration categories that appear in practice are saturated\rlap{.}%
  \footnote{
    If you know an example of a non-saturated cofibration category,
    please describe it at \cite{NonSatCofCat}!
    }

  We are mainly interested in cofibration categories
  in which every object is cofibrant,
  but allowing the possibility of non-cofibrant objects
  will simplify some of the statements below.
\end{nul}

\begin{example}
  \label{cc-of-model}
  A model category $M$ has an ``underlying'' cofibration category
  with the same underlying category, cofibrations, and weak equivalences.
  In fact, this data determines the model category $M$
  and so we may think of a model category as
  a (very) special kind of cofibration category.
  If $M$ and $N$ are model categories viewed as cofibration categories
  then a functor $F : M \to N$ is a left Quillen functor
  if and only if it is both exact and a left adjoint.
\end{example}

\begin{example}
  If $C$ is a cofibration category
  we write $C^\cof$ for the full subcategory of $C$ on the cofibrant objects.
  Then $C^\cof$ is again a cofibration category:
  the main point is that it inherits pushouts of cofibrations from $C$.
  Note that normally $C^\cof$ will not have all pushouts, even if $C$ does.
\end{example}

\begin{example}
  Write $\sSetfin$ for the category of finite simplicial sets.
  We equip $\sSetfin$ with the structure of a cofibration category
  by restricting the corresponding parts of
  the model category structure on $\sSet$.
  That is, a morphism $f : A \to B$ of $\sSetfin$ is:
  \begin{itemize}
  \item a cofibration if it is a monomorphism;
  \item a weak equivalence if it is a weak equivalence in $\sSet$.
  \end{itemize}
  The axioms for a cofibration category follow from
  the fact that $\sSet$ is a model category
  with the exception of the last one, the factorization axiom.
  For any map $f : A \to B$,
  the mapping cylinder factorization
  $A \to A \times \Delta^1 \amalg_A B \xrightarrow{\sim} B$
  expresses $A$ as a cofibration followed by a weak equivalence,
  with the intermediate object $A \times \Delta^1 \amalg_A B$
  again a finite simplicial set.
\end{example}

\begin{example}
  Write $\Topfin$ for the full subcategory of $\Top$
  on those objects which are homeomorphic to finite cell complexes.
  We equip $\Topfin$ with the structure of a cofibration category
  by declaring a morphism $f : A \to B$ to be
  \begin{itemize}
  \item
    a cofibration if it is (homeomorphic to)
    the inclusion of a finite relative cell complex;
  \item
    a weak equivalence if it is a homotopy equivalence.
  \end{itemize}
  By Whitehead's theorem, the homotopy equivalences in $\Topfin$
  are also those maps which are weak equivalences in $\Top$.
  Again, the cofibration category axioms for $\Topfin$
  are easily checked using the fact that $\Top$ is a model category
  (under the standard Serre--Quillen model category structure)
  together with the usual mapping cylinder construction.

  This example can also be constructed by applying
  the following result of Baues.
\end{example}

\begin{definition}
  \label{cc-fibrant}
  An object $X$ of a cofibration category is \emph{fibrant}
  if for any acyclic cofibration $f : A \to B$ with cofibrant domain,
  the induced map $\Hom(B, X) \xrightarrow{- \circ f} \Hom(A, X)$
  is surjective.
\end{definition}

\begin{example}
  Every object of $\Topfin$ is fibrant.
  On the other hand, $\sSetfin$ has very few fibrant objects.
  In particular, if $X$ is a fibrant finite simplicial set
  then each homotopy group of $X$ is finite,
  so $S^1 = \Delta^1/\partial \Delta^1$
  cannot have a fibrant approximation in $\sSetfin$.

  Suppose $M$ is a model category, $X$ is an object of $M$,
  and either of the following two hypotheses holds.
  \begin{enumerate}
  \item
    $M$ admits a class of generating acyclic cofibrations
    with cofibrant domains.
  \item
    $X$ is cofibrant.
  \end{enumerate}
  Then $X$ is fibrant in $M$ as a cofibration category
  (\cref{cc-of-model})
  if and only if $X$ is fibrant in the usual sense
  as an object of the model category $M$.
  (Under the second hypothesis, to check that $X \to *$
  has the right lifting property
  with respect to an acyclic cofibration $i : A \to B$,
  we first push forward $i$ to $X$,
  making its domain cofibrant.)
\end{example}

\begin{proposition}[{\cite[Theorem I.3.3]{Baues}}]
  \label{cc-of-icat}
  Let $C$ be an \emph{$I$-category}, that is,
  a category equipped with a class of cofibrations
  and a ``cylinder'' functor $I : C \to C$
  and natural transformations
  $p : I \to \id_C$, $i_0 : \id_C \to I$ and $i_1 : \id_C \to I$
  satisfying the following axioms:
  \begin{enumerate}[label=(I\arabic*)]
  \item
    $pi_0 = pi_1 = \id$.
  \item
    $C$ has an initial object and pushouts of cofibrations,
    and these colimits are preserved by the functor $I : C \to C$.
    Pushouts of cofibrations are again cofibrations.
  \item
    The cofibrations are closed under composition
    and for every object $X$, the map $\emptyset \to X$ is a cofibration.
    Moreover, cofibrations have the \emph{homotopy extension property}:
    for each cofibration $f : A \to B$, object $X$,
    and $\varepsilon \in \{0, 1\}$,
    any horizontal morphism
    \[
      \begin{tikzcd}
        IA \amalg_{A,\varepsilon} B \ar[r] \ar[d] & X \\
        IB \ar[ru, dashed]
      \end{tikzcd}
    \]
    admits a lift as shown by the dotted arrow.
    Here the vertical map is built from $i_\varepsilon : \id_C \to I$.
  \item
    For each cofibration $f : A \to B$,
    the induced ``relative cylinder'' map $B \amalg_{A,0} IA \amalg_{A,1} B \to IB$
    is a cofibration.
  \item
    For each object $X$, there exists an ``interchange'' map $T : IIX \to IIX$
    such that $Ti_\varepsilon = Ii_\varepsilon$
    and $T(Ii_\varepsilon) = i_\varepsilon$
    for $\varepsilon \in \{0, 1\}$.
  \end{enumerate}
  Then
  $C$ admits the structure of a cofibration category
  with the same cofibrations,
  and homotopy equivalences (defined using the cylinder functor $I$)
  as the weak equivalences.
  Moreover, every object of $C$ is both cofibrant and fibrant.
\end{proposition}

\begin{definition}
  \label{ho-cat}
  The \emph{homotopy category} $\Ho C$ of a category $C$ with weak equivalences $\mathcal{W}$
  (such as a cofibration category or model category)
  is the category $C[\mathcal{W}^{-1}]$ obtained by formally inverting
  the weak equivalences of $C$.
\end{definition}

\begin{nul}
  We recall the following basic facts about homotopy categories.
  Functors that preserve weak equivalences
  induce functors between homotopy categories.
  For a cofibration category $C$,
  the inclusion $C^\cof \to C$ induces
  an equivalence of categories $\Ho C^\cof \to \Ho C$ \cite[Theorem 6.1.6]{RB}.
  In particular, any object~$X$ has a cofibrant approximation $\tilde X$,
  which can be obtained by factoring $\emptyset \to X$
  into a cofibration $\emptyset \to \tilde X$
  followed by a weak equivalence $\tilde X \to X$.
  An exact functor $F : C \to D$ between cofibration categories
  does not necessarily preserve all weak equivalences,
  but it does preserve cofibrant objects
  and weak equivalences between cofibrant objects
  by ``Ken Brown's lemma''.
  We define the \emph{(left) derived functor} of $F$
  to be the functor $\lder F$ (defined up to unique natural isomorphism)
  which fits in the square below.
  \[
    \begin{tikzcd}
      \Ho C^\cof \ar[r, "\Ho F^\cof"] \ar[d, "\sim"'] &
      \Ho D^\cof \ar[d, "\sim"] \\
      \Ho C \ar[r, "\lder F"'] & \Ho D
    \end{tikzcd}
  \]
  When every object of $C$ is cofibrant,
  we simply take $\lder F = \Ho F$.

  In general, the morphisms of $\Ho C$ are equivalence classes of
  zigzags in which the backwards maps are weak equivalences.
  However, in certain cases
  we can give a more explicit description of the Hom-sets of $\Ho C$.
\end{nul}

\begin{definition}[{\cite[\S6.3]{RB}}]
  Two maps $f_0 : A \to B$ and $f_1 : A \to B$
  between cofibrant objects $A$ and $B$ of a cofibration category
  are \emph{strictly left homotopic}
  if there exists
  a factorization $A \amalg A \to I \to A$ of the fold map $A \amalg A \to A$
  into a cofibration followed by a weak equivalence
  together with a map $H : I \to B$
  such that the diagram below commutes.
  \[
    \begin{tikzcd}
      A \amalg A \ar[r, "{\langle f_0, f_1 \rangle}"] \ar[d] & B \\
      I \ar[ru, "H"']
    \end{tikzcd}
  \]
  The maps $f_0$ and $f_1$ are \emph{left homotopic}
  if they become strictly left homotopic
  after postcomposition with some acyclic cofibration $B \to B'$.
\end{definition}

\begin{nul}
  Left homotopic maps between cofibrant objects
  become equal in the homotopy category.
  In fact, the converse also holds \cite[Theorem 6.3.1]{RB}.
  However, not every element of $\Hom_{\Ho C}(A, B)$
  arises from a morphism in $C$ from $A$ to $B$.
  When $A$ and $B$ are cofibrant,
  each element of $\Hom_{\Ho C}(A, B)$ may be represented by a zigzag
  of the form $[A \to B' \xleftarrow{\sim} B]$.
  in which the map $B \xrightarrow{\sim} B'$
  may be chosen to be an acyclic cofibration
  \cite[Theorem 6.4.5]{RB}.

  Unlike in a model category,
  we generally cannot choose a single object~$\widehat B$
  so that every element of $\Hom_{\Ho C}(A, B)$
  is represented by a morphism from $A$ to $\widehat B$.
  (Consider the fact that the Hom-sets of $\sSetfin$ are finite,
  while $\Hom_{\Ho \sSetfin}(\partial \Delta^2, \partial \Delta^2) =
  \Hom_{\Ho \sSet}(\partial \Delta^2, \partial \Delta^2) = \ZZ$, as we will see later.)
  However, if $B$ happens to be fibrant,
  then $\Hom_{\Ho C}(A, B)$ is given by the familiar description.
\end{nul}

\begin{proposition}
  \label{cc-ho-cof-fib}
  Suppose $A$ is cofibrant and $X$ is both cofibrant and fibrant.
  Then strict left homotopy is an equivalence relation $\sim^{s\ell}$
  on $\Hom_C(A, X)$
  and the induced map $\Hom_C(A, X)/{\sim^{s\ell}} \to \Hom_{\Ho C}(A, X)$
  is a bijection.
  Furthermore, strict left homotopy may be detected
  using any fixed cylinder object for $A$.
\end{proposition}

\begin{proof}
  Under these hypotheses, any acyclic cofibration $X \to X'$ has a retraction,
  so left homotopy and strict left homotopy of maps into $X$ agree.
  Then the claim follows from
  \cite[Proposition 7.3.2, Theorem 6.3.1 and Lemma 6.3.2]{RB}.
\end{proof}

\begin{example}
  \label{ho-topfin-ff}
  For any objects $A$ and $B$ of $\Topfin$,
  the natural map $\Hom_{\Ho \Topfin}(A, B) \to \Hom_{\Ho \Top}(A, B)$
  is a bijection:
  both sides can be computed as homotopy classes of maps from $A$ to $B$
  in the classical sense.
\end{example}

\begin{proposition}
  \label{ho-ssetfin-ff}
  The inclusion $\sSetfin \to \sSet$
  induces a fully faithful functor $\Ho \sSetfin \to \Ho \sSet$.
\end{proposition}

\begin{proof}
  This follows from \cite[Theorem 4.6]{BaSc2} and \cite[Proposition 6.1]{BaSc1}.
  We give a quick sketch of the proof specialied to $\sSet$.
  Let $B$ be a finite simplicial set.
  While $B$ generally does not have a finite fibrant replacement,
  we can find a sequence $B = B_0 \to B_1 \to B_2 \to \cdots$
  of finite simplicial sets and acyclic cofibrations
  such that $\widehat B = \colim_{i \in \NN} B_i$ is fibrant in $\sSet$
  (for example by taking $B_i = \Ex^i B$).
  Then if $A$ is a finite simplicial set,
  every element $h$ of $\Hom_{\Ho \sSet}(A, B)$
  is represented by an element of $\Hom_\sSet(A, \widehat B)$
  and therefore by an element $g$ of $\Hom_\sSet(A, B_i)$ for some $i \in \NN$.
  We send $h$ to the element $[A \xrightarrow{g} B_i \xleftarrow{\sim} B]$
  of $\Hom_{\Ho \sSetfin}(A, B)$.
  Then one calculates that
  this map does not depend on the choice of $i$
  or on the choice of representative $g$,
  and provides an inverse to the natural map
  $\Hom_{\Ho \sSetfin}(A, B) \to \Hom_{\Ho \sSet}(A, B)$.
\end{proof}

\begin{example}
  \label{ssetfin-equiv-topfin}
  Geometric realization $|{-}| : \sSet \to \Top$ is a left Quillen functor
  and it restricts to an exact functor $|{-}| : \sSetfin \to \Topfin$.
  We claim that the latter functor
  induces an equivalence of homotopy categories.
  Indeed, we have a commutative square
  \[
    \begin{tikzcd}
      \Ho \sSetfin \ar[r, "\Ho |{-}|"] \ar[d] & \Ho \Topfin \ar[d] \\
      \Ho \sSet \ar[r, "\Ho |{-}|"] & \Ho \Top
    \end{tikzcd}
  \]
  in which the vertical morphisms are fully faithful functors
  by \cref{ho-ssetfin-ff,ho-topfin-ff}
  and the bottom functor is an equivalence.
  Thus, it suffices to show that the top functor is essentially surjective.
  This is a classical fact:
  every finite cell complex
  is homotopy equivalent to the geometric realization of a finite simplicial set.
\end{example}

\begin{nul*}
  Next we describe a very useful criterion for detecting when
  an exact functor induces an equivalence of homotopy categories,
  due to Cisinski.
  This criterion is the reason we included the saturation axiom
  in the definition of a cofibration category.
\end{nul*}

\begin{definition}[{\cite[3.6]{Cis10}}]
  An exact functor $F : C \to D$ between cofibration categories
  satisfies the \emph{approximation property}
  if:
  \begin{enumerate}[label=(AP\arabic*)]
  \item
    When restricted to cofibrant objects, $F$ reflects weak equivalences.
  \item
    Suppose $A$ is a cofibrant object of $C$
    and $f : FA \to B$ is a morphism to a cofibrant object of $D$.
    Then there exists a morphism $u : A \to A'$ with $A'$ cofibrant
    and a diagram
    \[
      \begin{tikzcd}
        FA \ar[r, "f"] \ar[d, "Fu"'] & B \ar[d, "\sim"] \\
        FA' \ar[r, "\sim"] & B'
      \end{tikzcd}
    \]
    with $B'$ cofibrant and the marked morphisms weak equivalences in $D$.
  \end{enumerate}
\end{definition}

\begin{proposition}[{\cite[Th\'eor\`eme 3.19]{Cis10}}]
  Let $F : C \to D$ be an exact functor between cofibration categories.
  Then $\lder F : \Ho C \to \Ho D$ is an equivalence of categories
  if and only if $F$ satisfies the approximation property.
\end{proposition}

\begin{nul}
  Let $C$ be a cofibration category and $Z$ a cofibrant object.
  Then the slice category $C_{Z/}$ has an induced cofibration category structure
  in which a morphism is a cofibration or a weak equivalence
  if and only if its underlying map is one in $C$.
  Note that normally $C_{Z/}$ will not have all objects cofibrant,
  even if $C$ does.

  When $C$ has a cofibrant terminal object $*$, we write $C_*$ for $C_{*/}$,
  the category of pointed objects of $C$.
\end{nul}

\begin{proposition}
  \label{cc-slice-equiv}
  Let $F : C \to D$ be an exact functor between cofibration categories
  which induces an equivalence of homotopy categories,
  and let $Z$ be a cofibrant object of $C$.
  Then the induced functor $F_{Z/} : C_{Z/} \to D_{FZ/}$ is also exact
  and induces an equivalence of homotopy categories.
\end{proposition}

\begin{proof}
  The verification that $F_{Z/}$ is exact is routine.
  To show that it induces an equivalence of homotopy categories,
  we check the approximation property.
  Here, for an object $X$ or map $f$ of $C_{Z/}$,
  we write $\underline{X}$ or $\underline{f}$
  for the underlying object or map of $C$.
  \begin{enumerate}[label=(AP\arabic*)]
  \item
    Suppose $f : A \to B$ is a map between cofibrant objects of $C_{Z/}$
    such that $F_{Z/}(f)$ is a weak equivalence.
    Then $\underline{A}$ and $\underline{B}$ are cofibrant objects of $C$
    (because $Z$ is cofibrant,
    and the maps $Z \to \underline{A}$ and $Z \to \underline{B}$
    are cofibrations)
    and $\underline{F_{Z/}(f)} = F\underline{f}$ is a weak equivalence.
    Hence by (AP1) for $F$,
    $\underline{f}$ and so also $f$ are weak equivalences.
  \item
    Similarly, suppose given cofibrant objects $A$ of $C_{Z/}$ and $B$ of $D_{FZ/}$
    and a morphism $f : F_{Z/}A \to B$.
    Unwinding the definitions, this corresponds to a commutative diagram
    \[
      \begin{tikzcd}[column sep=tiny]
        & FZ \ar[ld, "Fa"'] \ar[rd, "b"] \\
        F \underline{A} \ar[rr, "\underline{f}"] && \underline{B}
      \end{tikzcd}
    \]
    in $D$, in which the maps $a$ and $b$ are cofibrations.
    In particular, $\underline{A}$ and $\underline{B}$ are cofibrant
    because $Z$ is.
    So we may apply (AP2) for $F$ and enlarge the diagram to
    \[
      \begin{tikzcd}[column sep=tiny]
        & FZ \ar[ld, "Fa"'] \ar[rd, "b"] \\
        F \underline{A} \ar[rr, "\underline{f}"] \ar[d, "F\underline{u}"'] &&
        \underline{B} \ar[d, "\sim"] \\
        F \underline{A'} \ar[rr, "\sim"] && \underline{B'}
      \end{tikzcd}
    \]
    for cofibrant objects $\underline{A'}$ and $\underline{B'}$
    and a morphism $\underline{u} : \underline{A} \to \underline{A'}$ in $C$.
    Using the factorization axiom as described in \cite[Scholie 3.7]{Cis10},
    we may moreover assume that $\underline{u}$
    and the map $\underline{B} \to \underline{B'}$
    are cofibrations.
    We may regard $\underline{u}$ and $\underline{B'}$
    as underlying a morphism $u : A \to A'$ of $C_{Z/}$
    and an object $B'$ of $D_{FZ/}$ respectively,
    with the structural maps determined by the above diagram.
    This translates back into a diagram
    \[
      \begin{tikzcd}
        F_{Z/}A \ar[r, "f"] \ar[d, "F_{Z/}u"'] & B \ar[d, "\sim"] \\
        F_{Z/}A' \ar[r, "\sim"] & B'
      \end{tikzcd}
    \]
    of the form required for (AP2).
    The objects $A'$ and $B'$ are cofibrant because
    their structural maps $Z \to \underline{A'}$ and $FZ \to \underline{B'}$
    are cofibrations.
    \qedhere
  \end{enumerate}
\end{proof}

\begin{remark}
  The natural functor $\Ho(C_{Z/}) \to \Ho(C)_{Z/}$ is generally not faithful.
  For example, the morphisms of $\Ho(\Top_{*/})$
  are basepoint-preserving homotopy classes of basepoint-preserving maps,
  while the morphisms of $\Ho(\Top)_{*/}$
  are \emph{free} homotopy classes of maps that preserve the basepoint component.
\end{remark}

\begin{proposition}
  \label{slice-fibrant}
  Suppose $C$ is a cofibration category,
  $Z$ is a cofibrant object of $C$
  and $X$ is an object of $C_{Z/}$
  whose underlying object $\underline{X}$ is fibrant in $C$.
  Then $X$ is fibrant in $C_{Z/}$.
  In particular, if every object of $C$ is fibrant,
  the same is true of $C_{Z/}$.
\end{proposition}

\begin{proof}
  We have to check that
  for every cofibration $f : A \to B$ in $C_{Z/}$ with cofibrant domain,
  any extension problem of the form below admits a lift.
  \[
    \begin{tikzcd}
      A \ar[r] \ar[d, "f"'] & X \\
      B \ar[ru, dashed]
    \end{tikzcd}
  \]
  Unwinding the definitions, this is equivalent to finding a lift in the diagram
  \[
    \begin{tikzcd}[column sep=tiny]
      & Z \ar[ld] \ar[rd] \\
      \underline{A} \ar[rr] \ar[d, "\underline{f}"'] & & \underline{X} \\
      \underline{B} \ar[rru, dashed]
    \end{tikzcd}
  \]
  and $\underline{f}$ is also a cofibration with cofibrant domain,
  so such a lift exists because $\underline{X}$ is fibrant.
\end{proof}

\begin{proposition}
  \label{ho-ssetfinptd-ff}
  The inclusion $\sSetfin_* \to \sSet_*$ induces
  a fully faithful functor $\Ho \sSetfin_* \to \Ho \sSet_*$.
\end{proposition}

\begin{proof}
  Consider the following diagram,
  where all the functors are induced by either inclusions
  or geometric realization.
  \[
    \begin{tikzcd}
      & \Ho {(\Topfin_*)^\cof} \ar[d, "(1)"] \\
      \Ho \sSetfin_* \ar[r, "(2)"] \ar[d, "(3)"'] &
      \Ho \Topfin_* \ar[d, "(4)"] \\
      \Ho \sSet_* \ar[r, "(5)"] & \Ho \Top_*
    \end{tikzcd}
  \]
  We want to prove that the functor labeled (3) is fully faithful.
  The horizontal functors are equivalences,
  using \cref{cc-slice-equiv,ssetfin-equiv-topfin}.
  Functor (1) is also an equivalence,
  and then functor (4) is fully faithful
  by the same argument as \cref{ho-topfin-ff}
  (every object of $\Topfin_*$ is fibrant).
\end{proof}

\begin{nul}
  Finally, we prove a specialized result
  regarding the unit of a relative adjunction
  between cofibration categories.
  We apply it in the main text to promote
  an exact functor inducing an equivalence on homotopy categories
  to a Quillen equivalence between model categories.

  We begin with a technical lemma.
  It corresponds to the following fact about model categories:
  if $F : C \rightleftarrows D : G$ is a Quillen adjunction
  and $\lder F : \Ho C \to \Ho D$ is fully faithful,
  then the ``derived unit'' of the adjunction is a weak equivalence.
  In the setting of cofibration categories
  we cannot generally expect an exact functor to have a right adjoint,
  so we assume instead that $F$ is part of
  an adjunction \emph{relative} to a functor $J : C \to C'$ (\cref{relative-adjoint}).
  The conclusion is therefore weaker:
  we can only deduce that the unit map
  ``looks like a weak equivalence''
  to objects of $C'$ in the image of the functor $J$.
  (Even for $C = \sSetfin$, $C' = \sSet$, $J$ the inclusion,
  this is not by itself sufficient to conclude that
  the unit map is actually a weak equivalence.)
\end{nul}

\begin{lemma}
  Let $C$, $C'$ and $D$ be cofibration categories
  and $F : C \to D$, $U : D \to C'$ an adjunction relative to $J : C \to C'$
  (so there is an isomorphism
  $\Hom_D(FA, X) \xrightarrow{\sim} \Hom_{C'}(JA, UX)$
  natural in $A \in C$ and $X \in D$)
  satisfying the following conditions:
  \begin{enumerate}
  \item
    All objects of $C$ and of $C'$ are cofibrant.
  \item
    All objects of $D$ are fibrant,
    as are their images under $U$.
  \item
    $J$ and $F$ are exact.
  \item
    $\lder J : \Ho C \to \Ho C'$ and $\lder F : \Ho C \to \Ho D$
    are fully faithful.
  \end{enumerate}
  Then for any objects $A$ and $K$ of $C$,
  the unit map $\eta_K : JK \to UFK$
  (corresponding under the relative adjunction to $\id_{FK} : FK \to FK$)
  induces an isomorphism
  $(\eta_K)_* : \Hom_{\Ho C'}(JA, JK) \to \Hom_{\Ho C'}(JA, UFK)$.
\end{lemma}

\begin{proof}
  First, let $A$ be an object of $C$ and $X$ a cofibrant object of $D$.
  In the diagram
  \[
    \begin{tikzcd}
      \Hom_D(FA, X) \ar[r, "\sim"] \ar[d] & \Hom_{C'}(JA, UX) \ar[d] \\
      \Hom_{\Ho D}(FA, X) \ar[r, dashed, "\sim"] & \Hom_{\Ho C'}(JA, UX)
    \end{tikzcd}
  \]
  each vertical map is the quotient by the strict left homotopy relation,
  because all objects involved are cofibrant and $X$ and $UX$ are fibrant.
  Moreover, choosing a cylinder object
  $A \amalg A \to I \to A$ for $A$ in $C$,
  we may detect these strict left homotopy relations
  on the images of this cylinder object under $F$ and $J$.
  Then an adjunction argument shows that
  these two strict left homotopy relations correspond,
  so the adjunction isomorphism descends to the homotopy category
  as shown by the dotted arrow.

  Now let $A$ and $K$ be objects of $C$.
  We next claim the square
  \[
    \begin{tikzcd}
      \Hom_{\Ho C}(A, K) \ar[r] \ar[d] & \Hom_{\Ho D}(FA, FK) \ar[d, dashed] \\
      \Hom_{\Ho C'}(JA, JK) \ar[r, "(\eta_K)_*"] & \Hom_{\Ho C'}(JA, UFK)
    \end{tikzcd}
    \tag{$*$}
  \]
  commutes,
  where the dotted arrow was constructed above, taking $X = FK$.
  We can represent a general element of $\Hom_{\Ho C}(A, K)$
  in the form $[A \to K' \xleftarrow{\sim} K]$
  with the map $K \to K'$ an acyclic cofibration.
  Its image in $\Hom_{\Ho C'}(JA, UFK)$ under the two solid arrows
  is represented by the zigzag
  $[JA \to JK' \xleftarrow{\sim} JK \xrightarrow{\eta_K} UFK]$,
  while its image under the top map is $[FA \to FK' \xleftarrow{\sim} FK]$.
  Since $FK \xrightarrow{\sim} FK'$ is an acyclic cofibration
  and $FK$ is fibrant,
  we can choose a retraction $r : FK' \to FK$,
  so that the composition $FK \to FK' \xrightarrow{r} FK$ is the identity.
  Then the homotopy class $[FA \to FK' \xleftarrow{\sim} FK]$
  is represented by a map $FA \to FK' \xrightarrow{r} FK$ in $D$
  which (by naturality) the adjunction isomorphism takes to the composition
  $JA \to UFK' \xrightarrow{Ur} UFK$ in $C'$.
  Hence, in the diagram
  \[
    \begin{tikzcd}
      JA \ar[rd] & JK \ar[r, "\eta_K"] \ar[d, "\sim"'] & UFK \ar[d] \\
      & JK' \ar[r, "\eta_{K'}"] & UFK' \ar[u, bend right, shift right, pos=.4, "Ur"']
    \end{tikzcd}
  \]
  we need to check that the zigzags
  \[
    \begin{tikzcd}
      JA \ar[rd] & JK \ar[r, "\eta_K"] \ar[d, "\sim"'] & UFK \\
      & JK'
    \end{tikzcd}
  \]
  and
  \[
    \begin{tikzcd}
      JA \ar[rd] & & UFK \\
      & JK' \ar[r, "\eta_{K'}"] & UFK' \ar[u, bend right, shift right, pos=.4, "Ur"']
    \end{tikzcd}
  \]
  become equal in $\Ho C'$.
  Indeed, we compute
  \begin{align*}
    & [JA \to JK' \xrightarrow{\eta_{K'}} UFK' \xrightarrow{Ur} UFK] \\
    {}={} & [JA \to JK' \xleftarrow{\sim} JK \xrightarrow{\sim} JK'
            \xrightarrow{\eta_{K'}} UFK' \xrightarrow{Ur} UFK] \\
    {}={} & [JA \to JK' \xleftarrow{\sim} JK \xrightarrow{\eta_K} UFK
            \to UFK' \xrightarrow{Ur} UFK] \\
    {}={} & [JA \to JK' \xleftarrow{\sim} JK \xrightarrow{\eta_K} UFK].
  \end{align*}

  Finally, we conclude that $(\eta_K)_*$ is an isomorphism
  because the other three maps in ($*$) are
  (the top and left by the hypotheses on $J$ and $F$,
  and the dotted one by the first part of the proof).
\end{proof}

\begin{nul}
  Now we specialize to the following situation.
  Let $J : \sSetfin \to \sSet$ be the inclusion,
  let $D$ be a cofibration category
  and suppose $|{-}| : \sSetfin \to D$ and $\Sing : D \to \sSet$
  is an adjunction relative to $J$
  satisfying the following conditions:
  \begin{enumerate}
  \item
    Every object of $D$ is fibrant.
  \item
    $|{-}|$ is exact and induces an equivalence of homotopy categories.
  \end{enumerate}
  Then $\Sing X$ is also fibrant for every object $X$ of $D$,
  because $\sSet$ has generating acyclic cofibrations
  which are the images under $J$ of maps between cofibrant objects of $\sSetfin$.
  Furthermore, we saw earlier that $J$ induces
  a fully faithful functor $\lder J = \Ho J : \Ho \sSetfin \to \Ho \sSet$.
  Therefore, the conditions for the lemma are satisfied.

  Next, we pass to the pointed setting.
  The inclusion $J_{*/} : \sSetfin_* \to \sSet_*$ is also an exact functor
  which induces a fully faithful functor on homotopy categories
  (\cref{ho-ssetfinptd-ff}),
  and $|{-}|$ induces an exact functor
  $|{-}|_{*/} : \sSetfin_* \to D_{|*|/}$
  which again induces an equivalence of homotopy categories
  (\cref{cc-slice-equiv}).
  Define $\Sing' : D_{|{*}|/} \to \sSet_*$ by sending
  $(|{*}| \xrightarrow{x} X)$ to
  the pointed simplicial set given by the composition
  $({*} \xrightarrow{\eta_*} \Sing |{*}| \xrightarrow{\Sing x} \Sing X)$.
  One easily checks that $|{-}|_{*/}$ and $\Sing'$
  form an adjunction relative to $J_{*/}$
  whose unit map on an object $K$ of $\sSetfin_*$
  has underlying map given by
  the unit map of the original relative adjunction
  on the underlying object of $K$.
  Finally, by \cref{slice-fibrant}
  every object of $D_{|{*}|/}$ is fibrant
  and their images under $\Sing'$ are also fibrant.
  Thus, the lemma also applies to this pointed situation.
\end{nul}

\begin{proposition}
  \label{cc-unit-weq}
  In the above situation,
  the unit map $\eta_K : K \to \Sing |K|$ is a weak equivalence
  for every object $K$ of $\sSetfin$.
\end{proposition}

\begin{proof}
  We check that $\eta_K$ induces an isomorphism
  on $\pi_0$ and on $\pi_n$ for every $n \ge 1$ and choice of basepoint of $K$.
  Here, by $\pi_n$ of a general pointed simplicial set
  we really mean the simplicial homotopy groups
  of a functorial fibrant replacement in $\sSet_*$,
  so that $\pi_n(-)$ is naturally isomorphic to
  $\Hom_{\Ho \sSet_*}((\Delta^n/\partial \Delta^n, *), -)$.

  The map on $\pi_0$ induced by $\eta_K$ is isomorphic to
  \[
    (\eta_K)_* : \Hom_{\Ho \sSet}(*, K) \to \Hom_{\Ho \sSet}(*, \Sing |K|)
  \]
  and so it is an isomorphism by the lemma.
  For $\pi_n$, let $k \in K_0$ be a vertex of~$K$.
  We must show that the induced map
  $\pi_n(K, k) \to \pi_n(\Sing |K|, \eta_K(k))$
  is an isomorphism.
  Let $K_* = (K, k) \in \sSetfin_*$.
  Then the pointed map $(K, k) \to (\Sing |K|, \eta_K(k))$
  is the unit map $K_* \to \Sing' |K_*|_{*/}$ of
  the relative adjunction between $|{-}|_{*/}$ and $\Sing'$.
  By the lemma again,
  $\Hom_{\Ho \sSet_*}((\Delta^n/\partial \Delta^n, *), -)$
  takes this map to an isomorphism.
\end{proof}